\documentclass[12pt]{amsart}
\usepackage{amscd, amsfonts, amssymb,pictex}
\usepackage{amscd}
\usepackage{epsfig}
\usepackage{amssymb}
\usepackage[mathscr]{eucal}
\usepackage{verbatim}
\usepackage{fullpage}
\usepackage{latexsym}
\usepackage{lscape}

\newcommand\bb{\mathfrak b}
\newcommand\cc{\mathfrak c}
\renewcommand\gg{\mathfrak g}
\newcommand\hh{\mathfrak h}
\newcommand\ttt{\mathfrak t}

\newcommand\mm{\mathfrak m}
\newcommand\aaa{\mathfrak a}
\newcommand\nn{\mathfrak n}
\newcommand\uu{\mathfrak u}
\newcommand\CO{\mathcal O}
\newcommand\pp{\mathfrak p}
\newcommand\qq{\mathfrak q}

\newcommand{\cha}{\mathop{\mathrm{char}}\nolimits}
\newcommand{\htt}{\mathop{\mathrm{ht}}\nolimits}
\newcommand{\fib}{\mathop{G\!\ast_H\!Y}\nolimits}
\newcommand{\fibp}{\mathop{G\!\ast_P\!Y}\nolimits}
\newcommand{\fibl}{\mathop{G\!\ast_P\!\gg_{\geqslant 2}}\nolimits}
\newcommand{\rk}{\mathop{\mathrm{rk}}\nolimits}
\newcommand{\lv}{\mathop{\mathrm{lv}}\nolimits}
\newcommand{\supp}{\mathop{\mathrm{supp}}\nolimits}

\newcommand\inverse{{^{-1}}}

\newcommand\B{\mathcal B}
\newcommand\D{\mathcal D}

\newcommand\N{\mathcal N}
\newcommand\U{\mathcal U}

\newcommand\ZZ{\mathbb Z}

\DeclareMathOperator{\ad}{ad} \DeclareMathOperator{\Ad}{Ad}
 
\DeclareMathOperator{\codim}{codim} 
\DeclareMathOperator{\SL}{SL} 
\DeclareMathOperator{\SO}{SO} \DeclareMathOperator{\SP}{Sp}
\DeclareMathOperator{\Hom}{Hom} \DeclareMathOperator{\Lie}{Lie}

\DeclareMathOperator{\rank}{rank}

\numberwithin{equation}{section}

\newtheorem{thm}[equation]{Theorem}

\newtheorem{lem}[equation]{Lemma}
\newtheorem{cor}[equation]{Corollary}
\newtheorem{prop}[equation]{Proposition}

\theoremstyle{definition}
\newtheorem{defn}[equation]{Definition}
\newtheorem{exmp}[equation]{Example}
\newtheorem{exmps}[equation]{Examples}
\theoremstyle{remark}
\newtheorem{rem}[equation]{Remark}

\theoremstyle{remark}
\newtheorem{rems}[equation]{Remarks}


\title[Spherical Nilpotent Orbits] % in Positive Characteristic]
{Spherical Nilpotent Orbits in Positive Characteristic}

\author[R.\  Fowler]
{Russell Fowler}
\address
{School of Mathematics, University of Birmingham,
Birmingham B15 2TT, UK}  %United Kingdom}
\author[G. R\"ohrle]{Gerhard R\"ohrle}
\address%[G.~R\"{o}hrle]
{Fakult\"at f\"ur Mathematik,
Ruhr-Universit\"at Bochum,
D-44780 Bochum, Germany}
\email{gerhard.roehrle@rub.de}
\urladdr{http://www.ruhr-uni-bochum.de/ffm/Lehrstuehle/Lehrstuhl-VI/roehrle.html}

\begin{document}

\begin{abstract}
Let $G$ be a connected reductive linear algebraic group
defined over an algebraically closed field of characteristic $p$.
Assume that $p$ is good for $G$.
In this note we classify all the
spherical nilpotent $G$-orbits in the Lie algebra of $G$. 
The classification is the same as in the
characteristic zero case obtained  by D.I.\ Panyushev in 1994, \cite{Pa3}:
for $e$ a nilpotent element in the Lie algebra of $G$, 
the $G$-orbit $G\cdot e$ is spherical if and only if 
the height of $e$ is at most $3$. 
\end{abstract}

\maketitle

\tableofcontents

\section{Introduction}
\label{s:intro}
Let $G$ be a connected reductive linear algebraic group
defined over an algebraically closed field $k$
of characteristic $p > 0$. With the exception of Subsection
\ref{sub:bad}, we assume throughout that $p$ is \emph{good} for $G$
(see Subsection \ref{sub:not} for a definition).

A \emph{spherical} $G$-variety $X$ is an (irreducible) algebraic $G$-variety 
on which a Borel subgroup $B$ of $G$ acts with a dense orbit.
Homogeneous spherical 
$G$-varieties $G/H$, for $H$ a closed subgroup of $G$, 
are of particular interest. They include flag varieties (when $H$ is a 
parabolic subgroup of $G$) as well as  
symmetric spaces (when $H$ is the fixed point subgroup of an involutive
automorphism of $G$). 
We refer the reader to \cite{Bri1.5} and \cite{Bri1}
for more information on spherical varieties and for their representation
theoretic significance.
These varieties enjoy a remarkable property:
a Borel subgroup of $G$ acts on a spherical $G$-variety only 
with a finite number of orbits. 
This fundamental result is due to M.~Brion \cite{Bri2} and
\'E.~B.~Vinberg \cite{Vi1} independently 
in characteristic $0$, and to F.~Knop \cite[2.6]{Kn} 
in arbitrary characteristic.

Let $\gg = \Lie G$ be the Lie algebra of $G$. 
%and let  $\N$ be the nilpotent cone of $\gg$. 
The aim of this note is to classify
the spherical nilpotent $G$-orbits in $\gg$.
In case $k$ is of characteristic zero, this classification 
was  obtained by D.I.\ Panyushev in 1994 in \cite{Pa3}.
The classification is the same 
in case the characteristic of $k$ is good for $G$:
for $e \in \gg$ nilpotent, 
$G\cdot e$ is spherical if and only if 
the height of $e$ is at most $3$ (Theorem \ref{thm 8.40}). 
The height of $e$ is the highest degree in the grading of $\gg$
afforded by a cocharacter of $G$ associated to $e$ 
(Definition \ref{defn:height}).

The methods employed by Panyushev in \cite{Pa3} 
do not apply in positive characteristic, e.g.\ parts of the argument
are based on the concept of ``stabilizers in general position''; it is unknown 
whether these exist generically in positive characteristic.
Thus a different approach is needed to address the question in 
this case.

\smallskip
We briefly sketch the contents of the paper.
In Section \ref{s:prelim} we collect the preliminary results we require. 
In particular, we discuss the concepts of complexity and sphericity, 
and more specifically the question of complexity of homogeneous spaces.
In Subsection \ref{sub:kempf} we recall the basic results of
Kempf--Rousseau Theory and in Subsection \ref{sub:assco} we 
recall the fundamental concepts of 
associated cocharacters for nilpotent elements from 
\cite[\S 5]{Ja1} and \cite{Pr1}.
There we also recall the grading of $\gg$ afforded by a cocharacter
associated to a given nilpotent element and 
define the notion of the height of a nilpotent element as the 
highest occurring degree of such a grading, Definition \ref{defn:height}.
The complexity of fibre bundles is discussed in 
Subsection \ref{sub:Fibbun} which is crucial for the sequel.
In particular, in Theorem \ref{thm 5.10} we show that the complexity
of a fixed nilpotent orbit $G\cdot e$ is given by 
the complexity of a smaller reductive group acting on a linear space.
Precisely, let  $\lambda$ be a cocharacter 
of $G$ that is associated to $e$. Then $P_\lambda$
is the destabilizing parabolic subgroup $P(e)$ defined by $e$, 
in the sense of Geometric Invariant Theory.
Moreover, $L = C_G(\lambda(k^*))$ is a  Levi subgroup of $P(e)$.
We show in Theorem \ref{thm 5.10} that the complexity
of $G\cdot e$ equals the complexity of the action of
$L$ on the subalgebra $\bigoplus_{i \geqslant 2}\gg(i,\lambda)$ 
of $\gg$ where the grading $\gg = \bigoplus_{i \in \ZZ}\gg(i,\lambda)$ 
is afforded by $\lambda$. 
In Subsection \ref{sub:sec 6.2} we 
recall the concept of a weighted Dynkin diagram associated to 
a nilpotent orbit from \cite[\S 5]{Ca}. There we also
present the classification of the 
parabolic subgroups $P$ of a simple algebraic group $G$ 
admitting a dense action of a Borel subgroup of a 
Levi subgroup of $P$ on the unipotent radical of $P$ from 
\cite[Thm.\ 4.1]{Br1}. Here we also remind the reader of 
the classification of the parabolic subgroups of $G$ 
with an abelian unipotent radical.

In Section \ref{sect:classification} we give the classification of the 
spherical nilpotent orbits in good characteristic:
a nilpotent element $e$ in $\gg$ is spherical if and only if the height of $e$
is at most $3$ (Theorem \ref{thm 8.40}).
In Subsections \ref{sub:ht2} and \ref{sub:ht4} we show that 
orbits of height $2$ are spherical and orbits of height at least $4$ are not,
respectively.
The subsequent subsections deal with the cases of height $3$ nilpotent 
classes. For classical groups these only occur for the orthogonal groups.
For the exceptional groups the height $3$ cases are handled
in Subsection \ref{sub:ex} with the aid of a computer programme 
of S.M.\ Goodwin.

In Section \ref{sec:appl} we discuss some further results and some 
applications of the classification.
In Subsection \ref{sub:dist} we discuss the spherical nilpotent orbits 
that are distinguished and in Subsection \ref{sub:orth} we extend a result
of Panyushev in characteristic zero to good positive characteristic:
a characterization of the spherical nilpotent orbits in terms of 
pairwise orthogonal simple roots, see Theorem \ref{thm 20.2}. 

In Subsection \ref{sub:ideals}
we discuss generalizations of results from  \cite{PaRo1} and \cite{PaRo2}
to positive characteristic.
%More precisely, in Proposition \ref{prop 20.4} and Theorem \ref{thm 20.03} 
%we describe the maximal ideals $\cc$ of a Borel subalgebra $\bb$ of $\gg$
%with the property that the $G$-saturation $G\cdot \cc$ is a spherical 
%variety. Moreover, 
In Theorem \ref{thm 20.04} we show that if $\aaa$ is an abelian ideal of $\bb$,
then $G\cdot \aaa$ is a spherical variety. 
In Subsection \ref{sub:geom}
we describe a geometric characterization of spherical orbits 
in simple algebraic groups from \cite{CaCaCo} and \cite{Carno}. 
Finally, in Subsection \ref{sub:bad} we very briefly touch on the
issue of spherical nilpotent orbits in bad characteristic.

Thanks to the fact that a Springer isomorphism between the unipotent variety
of $G$ and the nilpotent variety of $\gg$ affords 
a bijection between the unipotent $G$-classes in $G$ and the nilpotent
$G$-orbits in $\gg$ (cf.\ \cite{serre}), 
there is an analogous classification of the spherical 
unipotent conjugacy classes in $G$.

For results on algebraic groups we refer the reader to Borel's book
\cite{Bor}  and for information on nilpotent classes we cite 
Jantzen's monograph \cite{Ja1}.

\section{Preliminaries} 
\label{s:prelim}

\subsection{Notation}
\label{sub:not} 
Let $H$ be a linear algebraic 
group defined over an algebraically closed field $k$. We denote
the Lie algebra of $H$ by $\Lie H$ or by $\hh$.
We write $H^{\circ}$ for the
identity component of $H$ and $Z(H)$
for the centre of $H$.
The derived subgroup of $H$ is denoted by $\D H$ and
we write $\rank H$ for the dimension of a maximal torus of $H$.
The unipotent radical of $H$ is denoted by $R_u(H)$. We say that
$H$ is reductive provided $H^\circ$ is reductive. 
Let $K$ be a subgroup of $H$. We write
$C_H(K) = \{ h \in H \mid hxh\inverse = x \ \text{for all}\   x \in K\}$ 
for the centralizer of $K$ in $H$. 

Suppose $H$ acts morphically on an algebraic variety $X$.
Then we say that $X$ is an $H$-variety.
Let $x \in X$. Then $H \cdot x$ denotes the $H$-orbit of $x$ in $X$
and $C_H(x) = \{h \in H \mid h\cdot x = x \  \text{for all}\  h \in H\}$ 
is the stabilizer of $x$ in $H$.

For $e \in \hh$ we denote the centralizers of $e$ in $H$ and $\hh$
by $C_H(e) = \{h\in H \mid \Ad(h) e = e\}$ and $\cc_\hh(e) = \{x \in
\hh \mid [x,e] = 0\}$, respectively.
For $S$ a subset of $H$ we write 
$\cc_\hh(S) =\{x \in \hh \mid \Ad(s) x = x \ \text{for all}\  s \in S\}$
for the centralizer of $S$ in $\hh$.

Suppose $G$ is a connected reductive algebraic group.
By $\N$ we denote the nilpotent cone of $\gg$.
Let $T$ be a maximal torus of $G$. Let $\Psi = \Psi(G,T)$ denote
the set of roots of $G$ with respect to $T$. Fix a Borel subgroup
$B$ of $G$ containing $T$ and let
$\Pi = \Pi(G, T)$   % = \{\alpha_1, \ldots , \alpha_r\}
be the set of simple roots of $\Psi$ defined by $B$. Then $\Psi^+
= \Psi(B,T)$ is the set of positive roots of $G$ with respect to $B$. 
For $I \subset \Pi$, we denote by $P_I$ and  $L_I$ the \emph{standard}
parabolic and \emph{standard} Levi subgroups of $G$ defined by $I$,
respectively; see \cite[\S 2]{Ca}.

For $\beta \in \Psi^+$ write $\beta = \sum_{\alpha \in \Pi} c_{\alpha\beta}
\alpha$ with $c_{\alpha\beta} \in \mathbb N_0$. A prime $p$ is
said to be \emph{good} for $G$ if it does not divide
$c_{\alpha\beta}$ for any $\alpha$ and $\beta$, \cite[Defn.\ 4.1]{SpSt}.
Let $U = R_u(B)$ and set $\uu = \Lie U$.
For a $T$-stable Lie subalgebra $\mm$ of $\uu$ we write
$\Psi(\mm)=\{\beta \in \Psi^+ \mid \gg_{\beta}\subseteq \mm\}$
for the set of roots of $\mm$ (with respect to $T$).

For every root $\beta \in \Psi$ we choose a
generator $e_{\beta}$ for the corresponding root
space $\gg_{\beta}$ of $\gg$. Any element $e \in\uu$ can be
uniquely written as 
$e = \sum_{\beta \in \Psi^+}c_{\beta}e_{\beta}$, 
where $c_{\beta} \in k$. 
The \emph{support} of $e$ is defined as
$\supp(e)=\{\beta \in \Psi^+ \mid c_{\beta}\neq 0\}$.

The variety of all Borel subgroups of $G$ is denoted by  
$\B$. Note that $\B$ is a single conjugacy class $\B =\{B^g \mid g \in G\}$. 
Also note the isomorphism $\B \cong G/B$.

%Since $B$ is self-normalizing,  $\B$ is identified with $G/B$.

Let $Y(G) = \Hom(k^*,G)$ denote the set of \emph{cocharacters} 
(one-parameter subgroups) 
of $G$, likewise for a closed subgroup $H$ of $G$, we
set $Y(H) = \Hom(k^*,H)$ for the set of cocharacters of $H$.
For $\lambda\in Y(G)$ and $g\in G$ we define
$g\cdot \lambda \in Y(G)$ by $(g\cdot \lambda)(t) =
g\lambda(t)g^{-1}$ for $t \in k^*$; 
this gives a left action of $G$ on $Y(G)$.
For $\mu \in Y(G)$ we write 
$C_G(\mu)$ for the centralizer of $\mu$ under this action 
of $G$ which coincides with $C_G(\mu(k^*))$.

By a Levi subgroup of $G$ we mean a Levi subgroup of a parabolic
subgroup of $G$. The Levi subgroups of $G$ are precisely
the subgroups of $G$ which are of the form 
$C_G(S)$ where $S$ is a torus of $G$, \cite[Thm.\  20.4]{Bor}. 
Note that for $S$ a torus of $G$
the group $C_G(S)$ is connected, \cite[Cor.\ 11.12]{Bor}.

\subsection{Complexity}
\label{sub:Comp}
Suppose the linear algebraic group $H$ acts morphically on the 
(irreducible) algebraic variety $X$.
Let $B$ be a Borel subgroup of $H$. Recall that the 
\emph{complexity of}  $X$ (with
respect to the $H$-action on $X$) is defined as
 \[
\kappa_H(X) := \min_{x\in X}\codim_{X} B\cdot x,
\]
see also \cite{Bri1}, \cite{Kn},  \cite{LuVu}, \cite{Pa3}, and \cite{Vi1}.

Since the Borel subgroups of $H$ are conjugate in $H$ 
(\cite[Thm.\  21.3]{Hu}),
the complexity of the variety $X$ is well-defined.

Since a Borel subgroup of $H$ is connected, 
we have $\kappa_H(X)=\kappa_{H^{\circ}}(X)$. Thus for considering the 
complexity of an $H$-action, 
we may assume that $H$ is connected.

Concerning basic properties of complexity,
we refer the reader to \cite[\S 9]{Vi1}. 

We return to the general situation of a linear algebraic group $H$
acting on an algebraic variety $X$. For a Borel subgroup $B$ of $H$, we define 
\[
\Gamma_{X}(B) := \{x\in X \mid \codim_{X} B\cdot x = \kappa_H(X)\}\subseteq X. 
\]
Then we set
\[ 
\Gamma_{X} := \bigcup_{B \in \mathcal{B}}\Gamma_{X}(B) \subseteq X.
\]

For $x \in X$, we define 
\[ 
\Lambda_H(x) := \{B \in \mathcal{B} \mid \codim_{X} B\cdot x =\kappa_H(X)\}
\subseteq \B.
\]
\begin{rem}
\label{rem1.20} 
The following statements are 
immediate from the definitions.
\begin{itemize}
\item[(i)] If $H$ acts transitively on $X$, then $\Gamma_{X}=X$.
\item[(ii)] $B\in \Lambda_H(x)$ if and only if $x \in \Gamma_{X}(B)$.
\item[(iii)] $\Lambda_H(x)= \varnothing$ if and only if $x \notin
\Gamma_{X}$.
\end{itemize}
\end{rem}

The complexity of a reducible variety can easily be determined from 
the complexities of its irreducible
components: Since a Borel subgroup $B$ of $G$ is connected, it
stabilizes each irreducible component of $X$, 
cf.\ \cite[Prop.\ 8.2(d)]{Hu}. 
Let $x \in \Gamma_{X}(B)$ 
and choose an irreducible
component $X'$ of $X$ such that $x \in X'$. Then
$\kappa_G(X)=\kappa_G(X')+\codim_{X}X'$.
Therefore, from now on we may assume that $X$ is irreducible.

Next we recall the upper semi-continuity of
dimension, e.g.\ see \cite[Prop.\ 4.4]{Hu}.

\begin{prop}
\label{prop 1.10} 
Let $\varphi : X\to Y$ be a dominant morphism
of irreducible varieties. For $x \in X$, let
$\varepsilon_{\varphi}(x)$ be the maximal dimension of any
component of $\varphi^{-1}(\varphi(x))$ passing through $x$. Then
$\{x \in X \mid\: \varepsilon_{\varphi}(x)\geqslant n \}$ is
closed in $X$, for all $n \in \mathbb{Z}$.
\end{prop}

\begin{cor}
\label{cor 1.10} 
Let $X$ be an $H$-variety. The set $\{x \in X \mid \dim H\cdot
x \leqslant n \}$ is closed in $X$ for all $n \in \mathbb{Z}$. In
particular, the union of all $H$-orbits of maximal dimension in
$X$ is an open subset of $X$.
\end{cor}

\begin{lem}
\label{lem 1.40}
For every $B \in \mathcal{B}$, we have
$\Gamma_{X}(B)$ is a non-empty open subset of $X$.
\end{lem}

\begin{proof}
Note that $\Gamma_{X}(B)$ is the union of $B$-orbits of maximal dimension.
Thus, by Corollary \ref{cor 1.10}, 
$\Gamma_{X}(B)$ is open in $X$. 
%If $x \in \Gamma_{X}(B)$, 
%then $b\cdot x \in \Gamma_{X}(B)$ for all $b \in B$. 
%Thus $\Gamma_{X}(B)$ is the union of $B$-orbits, but all these
%orbits have maximal dimension, and each $B$-orbit in $X$ of maximal dimension  
%is contained in $\Gamma_{X}(B)$. Thus, by Corollary \ref{cor 1.10}, 
%$\Gamma_{X}(B)$ is open in $X$. 
\end{proof}

\begin{cor} 
\label{cor 1.20}
$\Gamma_{X}$ is open in $X$.
\end{cor}

Next we need an easy but useful lemma; the proof is elementary.

\begin{lem}
\label{lem 1.50}
Let $\varphi : X \to Y$ be an $H$-equivariant 
dominant morphism of irreducible $H$-varieties. For $x \in X$ set
$F_{\varphi(x)}=\varphi^{-1}(\varphi(x))$. Then $F_{\varphi(x)}$
is $C_H(\varphi(x))$-stable.
\end{lem}

%\begin{proof} 
%Let $g \in C_G(\varphi(x))$ and $y \in F_{\varphi(x)}$. 
%Since $\varphi$ is $G$-equivariant, we have
%$\varphi(g\cdot y)=g\cdot \varphi(y)=g\cdot \varphi(x) =
%\varphi(x)$. Hence $g\cdot y \in F_{\varphi(x)}$. 
%%Therefore, $F_{\varphi(x)}$ is $C_G(\varphi(x))$-stable. 
%\end{proof}

Before we can prove the main result of this subsection we need
another preliminary result, see  \cite[Thm.\  4.3]{Hu}.

\begin{thm}
\label{thm 1.20}
Let $\varphi: X \to Y$ be a
dominant morphism of irreducible varieties. Set $r = \dim X - \dim Y$. 
Then there is a non-empty open subset $V$ of $Y$ such
that $V\subseteq \varphi(X)$ and if $Y'\subseteq Y$ is closed,
irreducible and meets $V$ and $Z$ is a component of
$\varphi^{-1}(Y')$ which meets $\varphi^{-1}(V)$, then $\dim
Z = \dim Y' +r$. In particular, if $v \in V$, then $\dim
\varphi^{-1}(v) = r$.
\end{thm}

%\begin{proof} Since $Y \subset V$, we have $\varphi^{-1}(Y)\subset \varphi^{-1}(V)$ and thus every
%component of $\varphi^{-1}(Y)$ meets $\varphi^{-1}(V)$ and so has
%dimension $r$. \end{proof}

For the remainder of this section let $G$ be connected reductive.
Let $\varphi : X \to Y$ be a $G$-equivariant 
dominant morphism of irreducible $G$-varieties. 
Then $\kappa_G(Y)\leqslant \kappa_G(X)$, \cite[\S 9]{Vi1}.
%double check!! 
%Let $c\in \mathbb{Z}_{\geqslant 0}$
%(which, of course, depends on $G,X,Y$, the $G$-action on $X$
%and $Y$ and the morphism $X\to Y$) 
%such that $\kappa_G(X)=\kappa_G(Y)+c$. 
%It is natural to ask how to determine the constant $c$.
In the main result of this subsection
we give an interpretation for the difference 
$\kappa_G(X) - \kappa_G(Y)$ in terms of the complexity of
a smaller subgroup acting on a fibre of $\varphi$. 

\begin{thm}
\label{thm 1.10} 
Let $\varphi : X \to Y$ be a $G$-equivariant 
dominant morphism of irreducible $G$-varieties. For $x \in X$ set
$F_{\varphi(x)}=\varphi^{-1}(\varphi(x))$. Then for every $B \in
\mathcal{B}$ there exists $x \in \Gamma_{X}(B)$ such that
for $H = C_B(\varphi(x))^{\circ}$ we have
\[
\kappa_G(X)=\kappa_G(Y)+\kappa_H(Z),
\] 
where $Z$ is an irreducible component of $F_{\varphi(x)}$ passing through $x$.
\end{thm}

\begin{proof} 
Let $B\in \mathcal{B}$.
Let $V$ be a non-empty open subset of $Y$ which satisfies the
conditions in Theorem \ref{thm 1.20}. Since $Y$ is irreducible,
Lemma \ref{lem 1.40} implies that $\Gamma_{Y}(B)\cap V \neq
\varnothing$. For $y \in \Gamma_{Y}(B)\cap V$, 
Theorem \ref{thm 1.20} implies that any component of $\varphi^{-1}(y)$ 
has dimension $r = \dim X -\dim Y$, 
in particular, $\dim \varphi^{-1}(y) = r$.
Since $\varphi^{-1}(\Gamma_{Y}(B)\cap V)$ is open in $X$, 
we have 
$\varphi^{-1}(\Gamma_{Y}(B)\cap V)\cap \Gamma_{X}(B)\neq \varnothing$,
by Lemma \ref{lem 1.40}. 
Now choose $x\in \varphi^{-1}(\Gamma_{Y}(B)\cap V)\cap \Gamma_{X}(B)$. 
In particular, $\dim F_{\varphi(x)}=r$. Lemma \ref{lem 1.50} implies
that $F_{\varphi(x)}$ is $C_B(\varphi(x))$-stable. 
Clearly, $C_B(x)$ is the stabilizer
of $x$ in $C_B(\varphi(x))$. Thus we obtain

\begin{align*}
\codim_{F_{\varphi(x)}}C_B(\varphi(x))\cdot x
&= \dim F_{\varphi(x)} -\dim C_B(\varphi(x))\cdot x\\
&= r - \dim C_B(\varphi(x))+ \dim C_B(x)\\
&= \dim X -\dim Y - \dim C_B(\varphi(x))+\dim C_B(x)+\dim B-\dim B\\ 
&= \left( \dim X - \dim B+ \dim C_B(x) \right) \\
& \qquad \qquad -\left(\dim Y -\dim B+ \dim C_B(\varphi(x))\right)\\ 
&= \kappa_G(X)-\kappa_G(Y),
\end{align*}
where the last equality holds because $x \in \Gamma_{X}(B)$ and
$\varphi(x) \in \Gamma_{Y}(B)$. 

Let $Z$ be an irreducible
component of $F_{\varphi(x)}$ which passes through $x$. Theorem
\ref{thm 1.20} implies that $Z$ has the same dimension as
$F_{\varphi(x)}$. The connected group $H = C_B(\varphi(x))^{\circ}$ 
stabilizes $Z$. 
Note that for each $z \in Z$ we have $\varphi(z) = \varphi(x)$
and $C_B(z) = C_{C_B(\varphi(x))}(z)$ (observed for $z=x$ above).
Since $x \in \Gamma_{X}(B)$, 
$\dim C_B(x)$ is minimal among groups of the form $C_B(z)$ for $z \in Z$.
Therefore, because $C_B(z) = C_{C_B(\varphi(x))}(z)$,
we see that $\dim C_{C_B(\varphi(x))}(x)$ 
is minimal among groups of the form $C_{C_B(\varphi(z))}(z)$ for $z \in Z$.
We deduce that $x \in \Gamma_{Z}(H)$.
Consequently, 
\[
\kappa_{H}(Z)
= \dim Z -\dim C_B(\varphi(x))^{\circ}+\dim C_{C_B(\varphi(x))^{\circ}}(x)
= \codim_{F_{\varphi(x)}}C_B(\varphi(x))\cdot x.
\] 
The result follows. 
\end{proof}

\subsection{Spherical Varieties}
A $G$-variety $X$ is called \emph{spherical} if a
Borel subgroup of $G$ acts on $X$ with a dense orbit, that is
$\kappa_G(X) = 0$.
We recall some standard facts concerning spherical varieties, 
see \cite{Bri1}, \cite{Kn} and \cite{Pa3}.

First we recall  an important result due to \'{E}.B.\ Vinberg \cite{Vi1} 
and M.\ Brion \cite{Bri2}
independently in characteristic zero and F.\ Knop 
\cite[Cor.\ 2.6]{Kn} in arbitrary characteristic.
Let $B$ be a Borel subgroup of $G$.

\begin{thm}
\label{lem 1.60} 
A spherical $G$-variety consists only of a finite number of
$B$-orbits.
\end{thm}

We have an immediate corollary.
 
\begin{cor}
\label{cor 1.20a} 
%Suppose that $G$ acts on the irreducible variety $X$. 
The following are equivalent.
\begin{itemize}
\item[(i)] The $G$-variety $X$ is spherical.
\item[(ii)] There is an open $B$-orbit in $X$. 
\item[(iii)] The number of $B$-orbits in $X$ is finite. 
\end{itemize}
\end{cor}

\subsection{Homogeneous Spaces}
Let $H$ be a closed subgroup of $G$. 
Since $G/H$ is a $G$-variety, we may consider the complexity
$\kappa_G(G/H)$. Let $B$ be a Borel subgroup of $G$. The orbits of
$B$ on $G/H$ are in bijection with the $(B,H)$-double cosets of $G$. We have that
$\kappa_G(G/H) = \codim_{G/H}BgH/H$ for $gH \in \Gamma_{G/H}(B)$.
Clearly, $G$ acts transitively on $G/H$, so Remark \ref{rem1.20}(i)
implies that we can choose a Borel subgroup $B$ such that 
$B \in \Lambda_{G}(1H)$. Thus, for this choice of $B$, we have 
\begin{align}
\label{SubCom}
\notag \kappa_G(G/H) & = \codim_{G/H}B H/H = \dim G/H -\dim B H/H \\
& = \dim G/H -\dim B /B\cap H \\
\notag & = \dim G -\dim H - \dim B + \dim B \cap H. 
\end{align}

Following M.\ Kr\"{a}mer \cite{Kr}, a subgroup $H$ of $G$ is called
\emph{spherical} if $\kappa_G(G/H) = 0$. 

Since $\kappa_G(G/H) = \kappa_G(G/H^\circ)$, by \eqref{SubCom},
in considering the complexity of
homogeneous spaces $G/H$ we may assume that the
subgroup $H$ is connected.

\begin{comment}
\begin{cor}\label{cor 1.70} 
If $T$ is a torus, then
any subgroup of $T$ is spherical.
\end{cor}
\end{comment}

We have an easy lemma.

\begin{lem}
\label{lem 1.30} 
Let $G$ be connected reductive and let $H$ be a subgroup of $G$ 
which contains the unipotent radical of a Borel subgroup of $G$. 
Then $H$ is spherical. In particular, a parabolic subgroup of $G$ is
spherical.
\end{lem}

\begin{proof} 
Let $B$ be a Borel subgroup of
$G$ such that $U = R_u(B)\leqslant H$. Denote by $B^-$ the opposite
Borel subgroup to $B$, relative to some maximal torus of $B$, see
\cite[\S 26.2 Cor.\ C]{Hu}. The \emph{big cell} $B^-U$
is an open subset of $G$, \cite[Prop.\ 28.5]{Hu}. We 
have $B^-U \subseteq B^-H$, so $B^-H$ is a dense subset of
$G$. Thus, %by Corollary \ref{cor 1.20}, 
$G/H$ is spherical.
\end{proof}

\begin{rem}
\label{rem:simplespherical}
If both $G$ and $H$ are reductive, then $G/H$ is an affine
variety, see \cite[Thm.\  A]{Ri4}. This case has been studied
greatly. The classification of spherical reductive subgroups of
the simple algebraic groups in characteristic zero was obtained by M.\
Kr\"{a}mer \cite{Kr} and was shown to be the same in positive
characteristic by J.\ Brundan \cite{Br1}. M.\ Brion
\cite{Bri1.5} classifies all the spherical reductive subgroups of
an arbitrary reductive group in characteristic zero. 
In positive characteristic no such classification is known.
However, the classification of the reductive spherical subgroups in simple
algebraic groups in positive characteristic follows from
work of T.A.\ Springer \cite{Sp2} (see also G.\ Seitz \cite{seitz}), 
J.\ Brundan \cite{Br1} and R.\ Lawther \cite{Law}.

Important examples of reductive spherical subgroups are centralizers of
involutive automorphisms of $G$: Suppose that $\cha k \neq 2$ and let
$\theta$ be an involutive automorphism of $G$. Then the fixed
point subgroup $C_G(\theta)=\{ g \in G \mid \theta(g)=g\}$ of $G$
is spherical, see \cite[Cor.\ 4.3.1]{Sp2}.
\end{rem}

For more on the complexity and sphericity of homogeneous spaces
see \cite{Bri2}, \cite{LuVu} and \cite{Pa4}.

\begin{rem}
In order to compute the complexity of an orbit variety,
it suffices to determine the complexity of a homogeneous space.
For, suppose that $G$ acts on an algebraic variety $X$. Let $x\in
X$. Since $G$ is connected, the orbit $G\cdot x$ is irreducible.
The map $\pi_x : G/C_G(x) \to G\cdot x$,
by $\pi_x(gC_G(x))=g\cdot x$ is a bijective $G$-equivariant morphism, 
\cite[\S 2.1]{Ja1}. Thus, by applying Theorem \ref{thm 1.10}
to $\pi_x$, we have
\begin{equation}
\label{eq 1.10}
\kappa_G(G/C_G(x)) = \kappa_G(G\cdot x).
\end{equation}

The relevance of \eqref{eq 1.10} is that the left hand side is
easier to compute, since calculating $\kappa_G(G/C_G(x))$ only requires
the study of groups of the form $C_B(x)$, cf.\ \eqref{SubCom}, 
where $B$ is a Borel subgroup of $G$. 
\end{rem}

\subsection{Kempf--Rousseau Theory}
\label{sub:kempf}

Next we require some standard facts from Geometric Invariant Theory, see
\cite{Ke}, also see \cite[\S 2]{Pr1}, \cite[\S 7]{Ri3}.
Let $X$ be an affine variety and $\phi : k^* \to X$
be a morphism of algebraic varieties. We say that
$ \underset{t\to 0}{\lim}\,\phi(t)$ exists if there exists a
morphism $\widehat{\phi}:k\to X$ such that
$\widehat{\phi}|_{k^*}=\phi$. If such a limit exists, we set
$ \underset{t\to 0}{\lim}\,\phi(t)=\widehat{\phi}(0)$. Note,
that if such a morphism $\widehat{\phi}$ exists, it is necessarily unique.

Let $\lambda$ be a cocharacter of $G$. Define 
$P_{\lambda}=\{ x \in G\mid \underset{t\to 0}{\lim}\,
\lambda(t)x\lambda(t)^{-1} \text{ exists}\}$. 
Then $P_{\lambda}$ is a parabolic subgroup of $G$, 
the unipotent radical of $P_{\lambda}$ is given by 
$R_u(P_{\lambda}) = \{x \in G \mid \underset{t\to 0}{\lim}\,
\lambda(t)x\lambda(t)^{-1}=1\}$, and a Levi subgroup of
$P_{\lambda}$ is the centralizer $G_G(\lambda) = C_G(\lambda(k^*))$
of the image of $\lambda$ in $G$, 
\cite[\S 8.4]{Sp0.5}.

Let the connected reductive group $G$ act on the affine variety $X$
and suppose $x \in X$ is a point such that $G\cdot x$ is not closed in $X$.
Let $C$ denote the unique closed $G$-orbit in the 
closure of $G\cdot x$, cf.~\cite[Lem.~1.4]{Ri4}.
Set $\Lambda(x) := \{\lambda \in Y(G) \mid 
\underset{t\to 0}{\lim}\, \lambda(t)\cdot x  \textrm{ exists and lies in } C \}$. 
Then  there is a so-called \emph{optimal class} 
$\Omega(x) \subseteq \Lambda(x)$ of cocharacters associated to $x$.
%Roughly speaking, the optimal cocharacters  [add stuff here]
The following theorem is due to G.R.~Kempf, \cite[Thm.\ 3.4]{Ke} 
(see also \cite{rousseau}).

\begin{thm}
\label{thm 4.15}
Assume as above. Then we have the following:
\begin{itemize}
\item[(i)] 
$\Omega(x) \neq \varnothing$.
\item[(ii)] 
There exists a parabolic subgroup $P(x)$ of $G$ such that
$P(x) = P_\lambda$ for every $\lambda \in \Omega(x)$.
\item[(iii)] $\Omega(x)$ is a single $P(x)$-orbit.
\item[(iv)] For $g \in G$, we have $\Omega(g\cdot x) = g\cdot \Omega(x)$
and $P(g \cdot x) = gP(x)g^{-1}$. In particular, 
$C_G(x) \leqslant N_G(P(x)) = P(x)$.
\end{itemize}
\end{thm}

Frequently, $P(x)$ in Theorem \ref{thm 4.15}
is called the \emph{destabilizing}
parabolic subgroup of $G$ defined by $x \in X$.

\subsection{Associated Cocharacters}
\label{sub:assco}
In this subsection we closely follow A.\ Premet \cite{Pr1}; also see
\cite[\S 5]{Ja1}. We recall that $p$ is a good prime for $G$ throughout
this section.

Every cocharacter $\lambda \in Y(G)$ 
induces a grading of $\gg$:
\[
\gg = \bigoplus_{i \in \mathbb{Z}}\gg(i,\lambda),
\] 
where 
\[
\gg(i,\lambda)=\{x\in \gg\mid\Ad(\lambda(t))(x)=t^ix\text{ for all }t\in k^*\},
\]
see \cite[\S 5.1]{Ja1}. 
For $P_\lambda$ as in the the previous subsection,
we have the following equalities: 
$\Lie P_{\lambda}=\bigoplus_{i \geqslant 0}
\gg(i,\lambda); \Lie R_u(P_{\lambda})=\bigoplus_{i > 0}
\gg(i,\lambda)$; and $\Lie C_G(\lambda)=\gg(0,\lambda)$.
Frequently, we write $\gg(i)$ for $\gg(i,\lambda)$
once we have fixed a cocharacter $\lambda \in Y(G)$.

Let $H$ be a connected reductive subgroup of $G$.
A nilpotent element $ e \in \hh$ is called 
\emph{distinguished in $\hh$} provided 
each torus in $C_H(e)$ is contained in the centre of $H$, 
\cite[\S 4.1]{Ja1}.

\begin{comment}\begin{rem}\label{rem 4.05} Some
texts, for $\cha k=0$ or large positive characteristic, give
alternative definitions, for example \cite[\S 5.7]{Ca} states:
For $G$ simple, a nilpotent element $e \in \gg$ is distinguished
if it does not commute with any nonzero semisimple elements of
$\gg$. This can easily be seen to be equivalent, under the
assumptions on $\cha k$, to the definition above.
\end{rem}\end{comment}

The following characterization of distinguished nilpotent elements in 
the Lie algebra of a Levi subgroup of $G$ can be found in 
\cite[\S 4.6, \S 4.7]{Ja1}.

\begin{prop}
\label{prop4.20} 
Let $e \in \gg$ be nilpotent and let $L$ be a Levi subgroup of $G$. Then $e$
is distinguished in $\Lie L$ if and only if $L = C_G(S)$, where $S$
is a maximal torus of $C_G(e)$.
\end{prop} 

Next we recall the definition of an associated 
cocharacter, see \cite[\S 5.3]{Ja1}.

\begin{defn}
\label{def 4.30} 
A cocharacter $\lambda : k^* \to G$
is \emph{associated} to $e \in \N$ if
$e \in \gg(2,\lambda)$ and there exists a Levi subgroup $L$ of $G$
such that $e$ is distinguished in $\Lie L$, and 
$\lambda(k^*) \leqslant \D L$.
\end{defn}

\begin{rem}
\label{rem 4.06}
In view of Proposition \ref{prop4.20}, the last two conditions in  
Definition \ref{def 4.30} are equivalent to 
the existence of a maximal torus $S$ of
$C_G(e)$ such that $\lambda(k^*)\leqslant\D C_G(S)$. We will use
this fact frequently in the sequel.
\end{rem}

Let $e \in \N$.
In \cite[\S 2.4, Prop.\  2.5]{Pr1}, A.~Premet explicitly defines a 
cocharacter of $G$ which is associated to $e$.
Moreover, in  \cite[Thm.\  2.3]{Pr1}, Premet shows that each of these
associated cocharacters belongs to the optimal class $\Omega(e)$ 
determined by $e$.
Premet shows this under the so called \emph{standard hypotheses} on $G$, 
see \cite[\S 2.9]{Ja1}. These restrictions were subsequently removed by 
G.\ McNinch in \cite[Prop.\ 16]{Mc2} so that this fact holds
for any connected reductive group $G$ in good characteristic.
It thus follows from %\cite[Thm.\  2.3]{Pr1}, 
\cite[Prop.\ 16]{Mc2},
Theorem \ref{thm 4.15}(iv), and 
the fact that any two associated cocharacters are conjugate 
under $C_G(e)$, \cite[Lem.\ 5.3]{Ja1}, that 
all the cocharacters of $G$ associated to $e \in \N$ belong to the 
optimal class $\Omega(e)$ defined by $e$;
see also \cite[Prop.\ 18, Thm.\ 21]{Mc2}. 
This motivates and justifies 
the following notation which we use in the sequel.

\begin{defn}
\label{d:Gamma}
Let $e \in \gg$ be nilpotent. Then we denote the set of 
cocharacters of $G$ associated to $e$ by 
\[
\Omega_G^a(e)  := \{\lambda \in Y(G)\mid 
\lambda  \text{ is associated to } e  \} \subseteq \Omega(e).
\]
Further, if $H$ is a (connected) reductive subgroup of $G$ with $e \in \hh$ nilpotent 
we also write $\Omega_H^a(e)$ to denote the cocharacters of $H$ that are
associated to $e$.
\end{defn}

As indicated above, in good characteristic, associated cocharacters
are known to exist for any nilpotent element $e \in \gg$;
more precisely,  we have the following, \cite[\S 5.3]{Ja1}:

\begin{prop}
\label{prop 4.10} 
Suppose that $p$ is good for $G$.
Let $e \in \gg$ be nilpotent. Then $\Omega_G^a(e) \ne \varnothing$. 
Moreover, if $\lambda\in
\Omega_G^a(e)$ and $\mu \in Y(G)$, then $\mu \in \Omega_G^a(e)$ if and
only if $\mu$ and $\lambda$ are conjugate by an element of
$C_G(e)$.
\end{prop}

Fix a nilpotent element $e \in \gg$ and an associated cocharacter
$\lambda\in \Omega_G^a(e)$ of $G$. Set $P = P_{\lambda}$. 
By Theorem \ref{thm 4.15}(ii), $P$ only depends on $e$ 
and not on the choice of the
associated cocharacter $\lambda$. 
%We now consider how
%$P=C_G(\lambda)R_u(P)$ acts, via the adjoint action, on $\Lie
%P=\pp$. 
Note that $C_G(\lambda)$ stabilizes $\gg(i)$
for every  $i \in \mathbb{Z}$. For $n \in \mathbb{Z}_{\geqslant 0}$ we set
\[
\gg_{\geqslant n} = \bigoplus\limits_{i\geqslant n}\gg(i)\ \ \text{ and }\ \  
\gg_{>n} = \bigoplus\limits_{i>n}\gg(i).
\] 
Then we have 
\[
\gg_{\geqslant 0} = \Lie P \ \ \text{ and }\ \  
\gg_{>0} = \Lie R_u(P).
\] 
Also, $C_G(e) = C_P(e)$, by Theorem \ref{thm 4.15}(iv).

The next result is \cite[Prop.\ 5.9(c)]{Ja1}.

\begin{comment}\begin{proof} By definition, $C_G(\lambda)=\{x \in G \mid x
\lambda(t)x^{-1}=\lambda(t) \text{ for all } t \in k^*\}$. Let $Y
\in \gg(i)$ and $x \in C_G(\lambda)$. Now
$\Ad(\lambda(t))\Ad(x)(Y)=\Ad(xx^{-1}\lambda(t)x)(Y)=\Ad(x)\Ad(\lambda(t))(Y)
=t^i\Ad(x)(Y)$. Thus, $\Ad(x)(Y) \in \gg(i)$.
\end{proof}\end{comment}

\begin{prop}
\label{prop 4.30}
The $P$-orbit of $e$ in $\gg_{\geqslant 2}$ is dense in $\gg_{\geqslant 2}$.
\end{prop}

\begin{cor}
\label{cor 4.10} 
The $C_G(\lambda)$-orbit of $e$ in $\gg(2)$ is dense in $\gg(2)$.
\end{cor}

\begin{comment}\begin{proof} Lemma \ref{lem 4.45} states that $C_G(\lambda)$
stabilizes each $\gg(i)$. In particular, $C_G(\lambda)$ stabilizes
$\gg(2)$. It is known, see \cite[\S 5.10]{Ja1}, that
$(\Ad(R_u(P))(e)-e)\subseteq \gg_{\geqslant 3}$, thus
$\Ad(R_u(P))(e)\cap\gg(2)=\{e\}$. The result follows.
 \end{proof}

\begin{cor}\label{cor 4.11}If the characteristic of $k$ is very
good for $G$, then.\begin{enumerate}\item $[\pp,e]=\gg_{\geqslant
2}$;\item $[\gg_{>0},e]=\gg_{\geqslant 3}$;\item
$[\gg(0),e]=\gg(2)$.\end{enumerate}\end{cor} \begin{proof} The
assumption on $\cha k$ implies that $\dim [\pp,e]=\dim P\cdot
e=\dim \gg_{\geqslant 2}$, see \cite[\S 1.14]{Ca}, so
$[\pp,e]=\gg_{\geqslant 2}$. The other results follow from the
fact that $\pp=\gg(0)\bigoplus\gg_{>0}$ and $e \in \gg(2)$.
\end{proof}\end{comment}

Define 
\[
C_G(e,\lambda) := C_G(e)\cap C_G(\lambda). 
\]

\begin{cor}
\label{cor 3.02} Let $e \in \N$. Then 
\begin{itemize} 
\item[(i)] $\dim C_G(e) =\dim \gg(0)+\dim \gg(1)$; 
\item[(ii)] $\dim R_u(C_G(e))=\dim \gg(1)+\dim\gg(2)$; 
\item[(iii)] $\dim C_G(e,\lambda)=\dim \gg(0)-\dim \gg(2)$.
\end{itemize}
\end{cor}

\begin{proof}
As $C_G(e)=C_P(e)$, part (i) is immediate from Proposition
\ref{prop 4.30}. Using the fact that
$(\Ad(R_u(P)-1)(e)\subseteq \gg_{\geqslant 3}$ (e.g.\ see \cite[\S 5.10]{Ja1})
and Proposition \ref{prop 4.30}, we see that 
$\dim \Ad(R_u(P))(e)=\dim \gg_{\geqslant 3}$ and so $\dim C_{R_u(P)}(e)=\dim
\gg(1)+\dim\gg(2)$. Finally, part (iii) follows from the
first two. 
\end{proof}

The following basic
result regarding the structure of $C_G(e)$ can be found in
\cite[Thm.\  A]{Pr1}.

\begin{prop}
\label{pro 4.35} 
If $\cha k$ is good for $G$, then $C_G(e)$ is the semi-direct product of
$C_G(e,\lambda)$ and $C_G(e)\cap R_u(P)$. Moreover,
$C_G(e,\lambda)^{\circ}$ is reductive and $C_G(e)\cap R_u(P)$ is
the unipotent radical of $C_G(e)$.
\end{prop}

\begin{defn}
\label{defn:height}
Let $e \in \gg$ be nilpotent. 
The \emph{height} of $e$ with respect to an associated
cocharacter $\lambda\in \Omega_G^a(e)$ is defined to be
\[
\htt(e) := \max\limits_{ i\in \mathbb{N}}\{i\mid\gg(i,\lambda)\neq 0\}.
\]
Thanks to Proposition \ref{prop 4.10},  the height of $e$ 
does not depend on the choice of $\lambda\in\Omega_G^a(e)$.
Since conjugate nilpotent elements have the same height, we may
speak of the height of a given nilpotent orbit.
Since $\lambda\in\Omega_G^a(e)$, we have $\htt(e)\geqslant 2$ for any nilpotent
element $e \in \gg$, cf.\ Definition \ref{def 4.30}.
\end{defn}

Let $\gg$ be classical with natural module $V$. 
Set $n = \dim V$. We write a partition $\pi$ of $n$ in one of the following
two ways, either 
$\pi = (d_1, d_2, \ldots , d_r)$ with $d_1 \geqslant d_2 \geqslant \cdots \geqslant d_r \geqslant 0$ and
$\sum_{i=1}^r = n$; or  
$\pi = [1^{r_1} , 2^{r_2} , \ldots]$ with $\sum_i ir_i = n$.
These two notations are related by $r_i = |\{j \mid d_j = i\}|$ for $i \geqslant  1$.

For $\gg$ classical with natural module $V$ 
it is straightforward to determine the height of a nilpotent orbit
from the corresponding partition of $\dim V$.
We leave the proof of the next proposition to the reader.

\begin{prop}
\label{prop 4.50} 
Let $e \in \gg$ be nilpotent with
partition $\pi_e=(d_1,d_2, \ldots, d_r)$.
\begin{itemize}
\item[(i)] If
$\gg= \mathfrak{gl}(V)$, $\mathfrak{sl}(V)$ or $\mathfrak{sp}(V)$,
then $\htt(e)=2(d_1-1)$.  
\item[(ii)] If $\gg=
\mathfrak{so}(V)$, then $\htt(e)=\left\{%
\begin{array}{ll}
    2(d_1-1) & \text{ if } \:\:d_1=d_2, \\
    2d_1-3 & \text{ if } \:\:d_1=d_2+1, \\
    2(d_1-2) & \text{ if } \:\:d_1>d_2+1. \\
\end{array}%
\right.$
\end{itemize}
\end{prop}

\begin{rems}
\label{rem 3.30} 
(i). For $\cha k=0$, Proposition \ref{prop 4.50} 
was proved in \cite[Thm.\  $2.3$]{Pa2}.

(ii). If $e$ is a nilpotent element in $\mathfrak{gl}(V)$,
$\mathfrak{sl}(V)$ or $\mathfrak{sp}(V)$, then $\htt(e)$ is even.
If $e$ is a nilpotent element in $\mathfrak{so}(V)$, then
$\htt(e)$ is odd if and only if %$d_1$ is odd and 
$d_2=d_1-1$.
%Also, if $\htt(e)$ is odd, then $\htt(e) \equiv 3\mod 4$.
\end{rems}

\subsection{Fibre Bundles}
\label{sub:Fibbun}

Let $H$ be a closed subgroup
of $G$. Suppose that $H$ acts on an affine variety $Y$.
Define a morphic action of $H$ on the affine variety
$G\times Y$ by $h\cdot(g,y)=(gh,h^{-1}\cdot y)$ for $h \in H,g
\in G$ and $y \in Y$. 
Since $H$ acts fixed point freely on
$G\times Y$, every $H$-orbit in $G\times Y$ has dimension $\dim H$. 
There exists a surjective quotient morphism 
$\rho : G\times Y \to (G\times Y)/H$, \cite[\S 1.2]{MuFo}, \cite[\S 4.8]{PV}. 
We denote the quotient $(G\times Y)/H$ by $\fib$, 
the \emph{fibre bundle} associated to the \emph{principal bundle} 
$\pi :G\to G/H$ defined by $\pi(g)=gH$ and
\emph{fibre} $Y$. We denote the element $(g,y)H$ of $\fib$
simply by $g\ast y$, see \cite[\S 2]{Ri5}. 
Let $X$ be a $G$-variety and $Y \subseteq X$ be an
$H$-subvariety. The \emph{collapsing} of the fibre bundle $\fib$ is
the morphism $\fib \to G\cdot Y \subseteq X$ defined
by $g\ast y \rightarrow g\cdot y$.

Define an action
of $G$ on $\fib$ by $g\cdot(g'\ast y)=(gg')\ast y$ for $g,g' \in
G$ and $y \in Y$.
We then have a $G$-equivariant surjective morphism $\varphi : \fib \to G/H$ 
by $\varphi(g\ast y)=gH$. Note that
$\varphi^{-1}(gH)\cong Y$ for all $gH \in G/H$. 

\begin{comment}For $x \in G$ set $V_x=\varphi^{-1}(xH)$, so $V_x=\{g\ast v
\in \fib \mid g=xh \text{ for some } h \in H\}$. The morphism
$\phi :V_x\to V$, defined by $\phi(g\ast
v)=x^{-1}g\cdot v$, is clearly well-defined and surjective.
Suppose that $\phi(g\ast v)=\phi(k\ast u)$, so $g\cdot v=k\cdot
u$. Since $g^{-1}k\in H$, we have $g\ast v =g(g^{-1}k)\ast
(g^{-1}k)^{-1}\cdot v=k\ast u$, and so $\phi$ is bijective. Now
define a morphism $\psi : V\to V_x$ by
$\psi(v)=x\ast v$. Again, we have that $\psi$ is bijective. Also
$\psi(\phi(g\ast v))=g\ast v$ and $\phi(\psi(v))=v$. Thus, $\phi$
is an isomorphism of varieties, see \cite[\S 3]{Ha}, (with
inverse $\psi$). This fact also readily implies that

\begin{equation}\label{eq 5.1} \dim \fib = \dim G+\dim V -\dim
H.\end{equation}\end{comment}

\begin{prop}
\label{prop 5.10}
Let $H$ be a closed subgroup of $G$ and let $Y$ be an $H$-variety.
Suppose that $B$ is a Borel subgroup of $G$ such that
$\dim B\cap H$ is minimal 
(among all subgroups of the form $B'\cap H$ for $B'$ ranging over $\B$). 
Then we have
\[
\kappa_G(\fib)=\kappa_G(G/H)+\kappa_{B\cap H}(Y).
\]
\end{prop}

\begin{proof} 
We apply Theorem \ref{thm 1.10} to
the morphism $\varphi :\fib\to G/H$. 
Thus, for a Borel subgroup $B$ of $G$ and
$g\ast y\in \Gamma_{\fib}(B)$,  we have that
$\kappa_G(\fib)=\kappa_G(G/H)+\kappa_K(Z)$, where $Z$ is an
irreducible component of $\varphi^{-1}(\varphi(g\ast y))$ passing
through $g\ast y$, $K = C_B(gH)^{\circ}$. Note that 
$C_B(gH)=B\cap gHg^{-1}$. So, since $g\ast y\in \Gamma_{\fib}(B)$, 
the dimension of $g^{-1}C_B(gH)g = g^{-1}Bg\cap H$ is minimal. 
Now, as $\fib$ is a fibre bundle, for $x \in G$ 
we have $Y_x := \varphi^{-1}(\varphi(x\ast y))\cong Y$. 
Define a morphism $\phi : Y_x\rightarrow Y$ by 
$\phi(g\ast y)=x^{-1}g\cdot y$. 
Clearly, $xhx^{-1} \in B\cap xHx^{-1}$ acts on
$g\ast y \in Y_x$, as $xhx^{-1}\cdot (g\ast y)=xhx^{-1}g\ast y$.
Since $g=xh'$ for some $h' \in H$, 
we have $xhx^{-1}\cdot (g\ast y)=xhh'\ast y$. 
So $\phi(xhh'\ast y)=hh'\cdot y$. Thus, if we
define an action of $B\cap xHx^{-1}$ on $Y$ by 
$xhx^{-1}\cdot y = h\cdot y$, the morphism $\phi : Y_x\to Y$ becomes a
$(B\cap xHx^{-1})$-equivariant isomorphism. 
%By Theorem \ref{thm 1.10}, we have 
It follows that
$\kappa_{B\cap xHx^{-1}}(Y_x)=\kappa_{B\cap xHx^{-1}}(Y)$. 
Since $x^{-1}(B\cap xHx^{-1})x = x^{-1}Bx\cap H$, we
finally get $\kappa_{B\cap xHx^{-1}}(Y)=\kappa_{x^{-1}Bx\cap H}(Y)$. 
The result follows.
\end{proof}

\begin{comment}
Two parabolic subgroups $P$ and $Q$ of $G$ are called
\emph{opposite} if $P\cap Q$ is a Levi subgroup of both $P$ and
$Q$, see \cite[\S 14.20]{Bor}. The following result regarding
opposite parabolic subgroups of $G$ combines Propositions 14.21 and 14.22
from \cite{Bor}.

\begin{prop}
%\label{prop 5.20} 
Let $P$ and $Q$ be parabolic
subgroups of $G$. 
\begin{itemize} 
\item[(i)] For any Levi subgroup
$L$ of $P$ there exists precisely one parabolic subgroup of $G$
opposite to $P$ and containing $L$. 
\item[(ii)] Any two parabolic
subgroups of $G$ opposite to $P$ are conjugate by a unique element
of $R_u(P)$.
\item[(iii)] The group $P\cap Q$ is connected and $(P\cap
Q)R_u(P)$ is a parabolic subgroup of $G$. In particular, $P\cap Q$
contains a maximal torus of $G$.
\end{itemize}
\end{prop}
\end{comment}

Next we need a technical lemma.

\begin{lem}
\label{lem 5.10} 
Let $P$ be a parabolic subgroup of $G$. Then for $B$ ranging over $\B$, 
the intersection $B\cap P$ is minimal if and only if $B\cap P$ is a Borel
subgroup of a Levi subgroup of $P$.
\end{lem}

\begin{proof} 
We may choose a Borel subgroup $B$ of $G$ so that $BP$ is open dense in $G$,
cf.~the proof of Lemma \ref{lem 1.30}.
Then the $P$-orbit of the base point in $G/B \cong \B$ is open dense in $\B$.
Consequently, the stabilizer of this base point in $P$, that is $P \cap B$ 
is minimal among all the isotropy subgroups $P \cap B'$ for $B'$ in $\B$.
Clearly, $B$ is opposite to a Borel subgroup of $G$ contained in $P$.
Thanks to \cite[Cor.~14.13]{Bor}, $P \cap B$ contains a maximal torus $T$ of $G$.
Let $L$ be the unique Levi subgroup of $P$ containing $T$.
Then \cite[Thm.\ 2.8.7]{Ca} implies that 
$P \cap B = 
%(T \cap L)(T\ \cap R_u(P))(R_u(B) \cap L)(R_u(B) \cap R_u(P)) = 
T(R_u(B) \cap L)$.
Clearly, $T(R_u(B) \cap L)$ is solvable and thus lies in a Borel subgroup of $L$.
A simple dimension counting argument, using
Theorem \ref{thm 1.20} applied to the multiplication map 
$B \times P \to BP$ and the fact that $\dim BP = \dim G$,
shows that $P \cap B$ is a Borel subgroup of $L$.

Reversing the argument in the previous paragraph 
shows that if $P \cap B$ is a Borel subgroup of $L$,
then $BP$ is dense in $G$ and thus $P \cap B$ is minimal again in the sense
of the statement.
\end{proof}

Next we consider a special case of Proposition \ref{prop 5.10}.

\begin{lem}
\label{lem 5.20}
Let $P$ be a parabolic subgroup of $G$ and let $Y$ be a $P$-variety. 
Then 
\[\kappa_G(\fibp)=\kappa_L(Y),\] 
where $L$ is a Levi subgroup of $P$.
\end{lem}

\begin{proof} 
Proposition \ref{prop 5.10} implies
that $\kappa_G(\fibp)=\kappa_G(G/P)+\kappa_{B\cap P}(Y)$, where
$\dim B\cap P$ is minimal. Lemmas \ref{lem 1.30} and 
\ref{lem 5.10} imply that $\kappa_G(G/P)=0$ and $B\cap P$ is a
Borel subgroup of a Levi subgroup of $P$. The result follows.
\end{proof}

Let $e \in \N$ be a non-zero nilpotent element, $\lambda \in \Omega_G^a(e)$ 
be an associated cocharacter of $e$ and
$\gg=\bigoplus_{i\in \mathbb{Z}}\gg(i)$ be the grading of
$\gg$ induced by $\lambda$. Also let $P$ be the destabilizing 
parabolic subgroup of $G$ defined by $e$,
cf.\ Subsection \ref{sub:kempf}. In particular, we have 
$\Lie P=\gg_{\geqslant 0}$, see Subsection \ref{sub:assco}.

%Since $\gg_{\geqslant 2} \subseteq \gg$ is $P$-stable, we may
%consider the collapsing $\fibl \to G\cdot
%\gg_{\geqslant 2}$. The next lemma shows how this collapsing is
%related to the nilpotent orbit $G\cdot e$.

\begin{lem}
\label{lem 5.30} 
Let $e \in \N$.
Then $G\cdot\gg_{\geqslant 2}=\overline{G\cdot e}$. 
In particular, $\dim G\cdot \gg_{\geqslant 2} = \dim G\cdot e$.
\end{lem}

\begin{proof} 
Since $\gg_{\geqslant 2}$ is $P$-stable, 
$G\cdot\gg_{\geqslant 2}$ is closed, \cite[Prop.\ 0.15]{Hu2}. 
Thus, since $e \in
\gg(2)\subseteq\gg_{\geqslant 2}$, we have $\overline{G\cdot
e}\subseteq G\cdot\gg_{\geqslant 2}$. By Proposition \ref{prop 4.30}, 
$\overline{P\cdot e}=\gg_{\geqslant 2}$. Since
$\overline{P\cdot e}\subseteq \overline{G\cdot e}$, we thus have
$\gg_{\geqslant 2}\subseteq \overline{G\cdot e}$. Finally, as
$\overline{G\cdot e}$ is $G$-stable, 
$G\cdot \gg_{\geqslant 2}\subseteq \overline{G\cdot e}$. 
The result follows.
\end{proof}

\begin{thm}
\label{thm 5.10} 
Let $e \in \N$.
Then 
\[\kappa_G(G\cdot e) = \kappa_L(\gg_{\geqslant 2}),\]
where $L$ is a Levi subgroup of $P$.
\end{thm}

\begin{proof}
We have $\kappa_G(G\cdot e) = \kappa_G(G/C_G(e)) = \kappa_G(G/C_P(e))$,
thanks to \eqref{eq 1.10} and the fact that $G_G(e) = C_P(e)$.
Moreover, since $G \!\ast_P\! P/C_P(e) \cong G/C_P(e)$,
it follows from Lemma \ref{lem 5.20} that 
$\kappa_G(G/C_P(e)) = \kappa_L(P/C_P(e))$.
Finally, thanks to Proposition \ref{prop 4.30}
and \eqref{eq 1.10}, we obtain 
$\kappa_L(P/C_P(e)) =  \kappa_L(\gg_{\geqslant 2})$.
The result follows.
\begin{comment}
The collapsing $\phi: \fibl \to G\cdot \gg_{\geqslant 2}$ is 
$G$-equivariant and dominant. 
Theorem \ref{thm 1.10} and Lemma \ref{lem 5.30} imply
that for a Borel subgroup $B$ of $G$ and $g \ast x \in \Gamma_{\fibl}(B)$, 
we have $\kappa_G(\fibl) = \kappa_G(\overline{G\cdot e})+\kappa_{H}(Z)$,
where $Z$ is an irreducible component of $\phi^{-1}(\phi(g\ast x))$ 
passing through $g\ast x$ and $H = C_B(g\cdot x)^{\circ}$.
We claim that $\dim Z = 0$. 
Since 
$\dim Z = \dim \phi^{-1}(\phi(g\ast x))=\dim \fibl -\dim \overline{G\cdot e}$, 
it suffices to show that $\dim \fibl =\dim G\cdot e$.
So, as $\dim \fibl = \dim G -\dim P +\dim \gg_{\geqslant 2}$ and
$\dim P = \dim \gg_{\geqslant 0}$, we have $\dim \fibl = \dim G
-\dim \gg(0)-\dim \gg(1)$. 
Finally, by Corollary \ref{cor 3.02}(i), 
we have $\dim \fibl =\dim G - \dim C_G(e)= \dim G\cdot
e$. So the claim follows. Consequently, $Z = \{g\ast x\}$, and thus 
$\kappa_{H}(Z)=0$ and so
$\kappa_G(\fibl)=\kappa_G(\overline{G\cdot e})$. Since 
%$G\cdot e$ is dense in $\overline{G\cdot e}$, 
the inclusion map  $G\cdot e \to \overline{G\cdot e}$ is dominant,
we have $\kappa_G(\overline{G\cdot e}) = \kappa_G(G\cdot e)$, 
by Theorem \ref{thm 1.10}. The result 
follows from Lemma \ref{lem 5.20}.
\end{comment}
\end{proof}

\begin{rem}
For $\cha k=0$, Theorem \ref{thm 5.10} was
proved by Panyushev in \cite[Thm.\  4.2.2]{Pa2}.
\end{rem}

\begin{rem}
\label{rem 5.20} 
Thanks to Theorem \ref{thm 5.10}, 
in order to determine whether a nilpotent orbit is spherical,
it suffices to show that a Borel subgroup of a Levi subgroup of
$P$ acts on $\gg_{\geqslant 2}$ with a dense orbit. 
In our classification we pursue this approach.
\end{rem}

\subsection{Borel Subgroups of Levi Subgroups Acting on Unipotent
Radicals}
\label{sub:sec 6.2}

Let $e \in \gg$ be a non-zero nilpotent element and let 
$\lambda\in \Omega_G^a(e)$ be an associated
cocharacter for $e$.
Let $P=P_{\lambda}$ be the destabilizing parabolic subgroup defined by  
$e$. We denote the Levi
subgroup $C_G(\lambda)$ of $P$ by $L$.  Our next result is
taken from \cite[\S 3]{Ja1}. We only consider the
case when $G$ is simple, the extension to 
the case when $G$ is reductive is straightforward.

\begin{prop}
\label{prop 5.05}
Let $G$ be a simple classical algebraic group 
and $0 \ne e\in \gg$ be nilpotent with corresponding partition
$\pi_e=[1^{r_1},2^{r_2},3^{r_3},\ldots]$. 
Let $a_i,b_i,s,t \in \mathbb{Z}_{\geqslant 0}$ such that 
$a_i+1=\sum_{j \geqslant i}r_{2j+1}$,  
$b_i+1=\sum_{j\geqslant i}r_{2j}$, 
$2s=\sum_{j\geqslant 0} r_{2j+1}$, and
$2t+1=\sum_{j \geqslant 0} r_{2j+1}$. 
Then the structure of $\D L$ is as follows.
\begin{itemize}
\item[(i)] If $G$ is of type $A_n$, then $\D L$
is of type $\prod_{i\geqslant 0} A_{a_i} \times
\prod_{i\geqslant 1} A_{b_i}$.
\item[(ii)] If $G$ is of type $B_n$, then $\D
L$ is of type $\prod_{i\geqslant 1} A_{a_i} \times
\prod_{i \geqslant 1} A_{b_i}\times B_t$.

\item[(iii)] If $G$ is of type $C_n$, then $\D
L$ is of type $\prod_{i\geqslant 1} A_{a_i} \times
\prod_{i\geqslant 1} A_{b_i}\times C_s$.

\item[(iv)] If $G$ is of type $D_n$, then $\D
L$ is of type $\prod_{i\geqslant 1} A_{a_i} \times
\prod_{i\geqslant 1} A_{b_i}\times D_s$.
\end{itemize}
\end{prop}

We use the conventions that $A_0=B_0=C_0=D_0=\{1\}$,
$D_1\cong k^*$ and $D_2 = A_1\times A_1$.

\begin{comment}
\begin{exmp}\label{eg0.5}Let $G$ be of type $C_{10}$ and $e \in \gg$
be nilpotent with corresponding partition $\pi_e=[4,3^2,2^3,1^4]$.
Thus, $a_1=1,b_2=3$ and $s=3$. Proposition \ref{prop 5.05} says
that $\D L$ is of type $A_1\times A_3\times C_3$.
\end{exmp}
\end{comment}

In order to describe the Levi subgroups $C_G(\lambda)$ 
for the exceptional groups
we need to know more about associated cocharacters.
Let $T$ be a maximal torus of $G$ such that $\lambda(k^*)\leqslant T$. 
Now let $G_{\mathbb{C}}$ be the simple, simply connected 
group over $\mathbb{C}$ with the same root system as $G$. Let
$\gg_{\mathbb{C}}$ be the Lie algebra of $G_{\mathbb{C}}$. For a
nilpotent element $e \in \gg_{\mathbb{C}}$ we can find an
$\mathfrak{sl}_2$-triple containing $e$.
Let $h\in\gg_{\mathbb{C}}$ be the semisimple element of this
$\mathfrak{sl}_2$-triple. 
Note that $h$ is the image of $1$ under the differential 
of $\lambda_{\mathbb{C}} \in G_{\mathbb{C}}$ (corresponding to $\lambda$)
at $1$.
Then there exists a set of simple 
roots $\Pi$ of $\Psi$ such that
$\alpha(h) \geqslant 0$ for all $\alpha \in \Psi^+$ and
$\alpha(h) = m_{\alpha}\in\{0,1,2\}$ for all $\alpha \in \Pi$, 
see \cite[\S 5.6]{Ca}. For
each simple root $\alpha\in \Pi$ we attach the numerical label
$m_{\alpha}$ to the corresponding node of the Dynkin diagram.
The resulting labels form the \emph{weighted Dynkin diagram} 
$\Delta(e)$ of $e$. 
We denote the set of weighted Dynkin diagrams of $G$ by $\D(\Pi)$.
For $e,e' \in \gg_{\mathbb{C}}$ nilpotent, we have
that $\Delta(e)=\Delta(e')$ if and only
if $e$ and $e'$ are in the same $G_{\mathbb{C}}$-orbit.
%Following \cite[\S 5.7]{Ca}, we say that $e$ is \emph{even} resp.\ \emph{odd} 
%if all the labels of $\Delta(e)$ are even resp.\ odd.

In order to determine the weighted Dynkin diagram of a given
nilpotent orbit we refer to the method outlined in
\cite[\S 13]{Ca} for the classical groups, and to the tables
in \emph{loc.\ cit.} for the exceptional groups.

We return to the case when the characteristic of $k$ is good for $G$. 
In this case the classification of the 
nilpotent orbits does not depend on the field $k$.
\cite[\S 5.11]{Ca}. 
Recently, in  \cite{Pr1}
Premet gave a proof of this fact for the unipotent classes of $G$
which is free from case by case
considerations.
This applies in our case, since the classification
of the unipotent conjugacy classes in $G$ and of the
nilpotent orbits in $\N$ is the same in
good characteristic, \cite[\S 9 and \S 11]{Ca}.
First assume that $G$ is simply connected and 
that $G$ admits a finite-dimensional
rational representation such that the trace form on $\gg$ is
non-degenerate; see \cite[\S 2.3]{Pr1} for the motivation of
these assumptions. Under these assumptions,
given $\Delta \in \D(\Pi)$, there exists
a cocharacter $\lambda=\lambda_{\Delta}$ of $G$ which is associated to $e$, 
where $e$ lies in the dense $L$-orbit in $\gg(2,\lambda)$,
for $L = C_G(\lambda)$, such that
\begin{equation}
\label{eq 6.2}
\Ad(\lambda(t))(e_{\pm \alpha})=t^{\pm m_{\alpha}}e_{\pm \alpha} 
\text{ and } \Ad(\lambda(t))(x) = x
\end{equation} 
for all $\alpha \in \Pi, e_{\pm \alpha}\in \gg_{\pm \alpha}, x \in \ttt$
and $t \in k^*$, \cite[\S 2.4]{Pr1}. 
We extend this action linearly to all of $\gg$. Now
return to the general simple case. Let $\widehat{G}$ be the simple,
simply connected group with the same root datum as $G$. 
Then there exists a surjective
central isogeny $\pi :\widehat{G} \to G$, 
\cite[\S 1.11]{Ca}. Also, %by Lemma \ref{lem 4.70}, 
an associated cocharacter for
$e = d\pi(\widehat{e})$ in $\gg$ is of the form
$\pi\circ\widehat{\lambda}$, where $\widehat{\lambda}$ is a
cocharacter of $\widehat G$ that is associated  to
$\widehat{e}$ in $\widehat{\gg}$. This
implies that \eqref{eq 6.2} holds for an arbitrary simple
algebraic group, when the characteristic of $k$ is good for $G$.
%Also see the proof of Proposition 16 in \cite{Mc2}.

After these deliberations we can use the tables in
\cite[\S 13]{Ca} to determine the structure of the Levi subgroup
$C_G(\lambda)$ for the exceptional groups. 
Recall that
$\Lie C_G(\lambda)=\gg(0)$ and $\gg(0)$ is the sum of the root 
spaces $\gg_{\alpha}$, where $\alpha \in \Psi$ with 
$\langle\alpha,\lambda\rangle=0$. 
Let $\Pi_{0} = \{\alpha \in \Pi \mid m_\alpha = 0\}$,
the set of nodes $\alpha$ of the
corresponding weighted Dynkin diagram with label $m_\alpha = 0$. 
Then $C_G(\lambda)=\langle T,U_{\pm\alpha}\mid \alpha \in \Pi_{0}\rangle$.

\begin{comment}
\begin{exmp}\label{eg 6.20}Using Example \ref{eg 6.10} we have for $G$ of type
$B_5$ that $\D C_G(\lambda)$ is of type $A_1A_2$, and for $G$ of
Type $E_8$ that $\D C_G(\lambda)$ is of type $A_1A_1A_1$.
\end{exmp}
\end{comment}

It is straightforward to determine the height
of a nilpotent orbit from its associated 
weighted Dynkin diagram. 
Let $\tilde\alpha=\sum_{\alpha \in \Pi}c_{\alpha}\alpha$ 
be the highest root of $\Psi$. For each
simple root $\alpha \in \Pi$ we have $\gg_{\alpha}\subseteq
\gg(m_{\alpha})$ where $m_{\alpha}$ is the corresponding numerical
label on the weighted Dynkin diagram, by \eqref{eq 6.2}.

\begin{lem}
\label{lem 6.05} 
Let $\tilde\alpha$ be the highest root of $\Psi$ and set $d = \htt(e)$. 
Then $\gg_{\tilde\alpha}\subseteq \gg(d)$.
\end{lem} 

\begin{proof} 
Clearly, we have $\gg_{\tilde\alpha}\subseteq \gg(i)$
for some $i\geqslant 0$. The lemma is immediate, because 
%then follows from the following trivial observation: If
if $\tilde\alpha=\sum_{\alpha\in\Pi}c_{\alpha}\alpha$ and
$\beta=\sum_{\alpha\in\Pi}d_{\alpha}\alpha$ is any other
root of $\Psi$, then $c_{\alpha}\geqslant d_{\alpha}$ for all
$\alpha \in \Pi$.
\end{proof}

Lemma \ref{lem 6.05} readily implies 
\begin{equation}
\label{eq 6.05} 
\htt(e) = \sum_{\alpha \in \Pi}m_{\alpha}c_{\alpha}.
\end{equation} 
The  identity \eqref{eq 6.05} is also observed in \cite[\S 2.1]{Pa3}.

For the remainder of this section we assume that $G$ is simple. 
The generalization of each of the subsequent results to the case
when $G$ is reductive is straightforward.

For $P$ a parabolic subgroup of $G$ we set $\pp_u=\Lie R_u(P)$.

\begin{prop}
\label{prop 6.10} 
Let $P=LR_u(P)$ be an arbitrary parabolic
subgroup of $G$, where $L$ is a Levi subgroup of $P$. Then
\[
\kappa_G(G/L)=\kappa_L(P/L)=\kappa_L(R_u(P))=\kappa_L(\pp_u).
\]
\end{prop}

\begin{proof}
%We prove each of the equalities in turn.
Thanks to Lemma \ref{lem 5.20}, we have 
$\kappa_G(G/L) = \kappa_G(G\!\ast_P\! P/L) = \kappa_L(P/L)$.
\begin{comment}
%Consider the two quotient morphisms $\pi : G\to G/L$
%and $\psi : G \to G/P$. Since $\psi$ is constant
%on the fibres of $\pi$, applying Theorem \ref{UMP}, we obtain a
There is a canonical $G$-equivariant dominant morphism
$\varphi : G/L \to G/P$, by $\varphi(xL)= xP$ for $x \in G$. 
%Also, since $G$ acts on $G/L$ and
%$G/P$ by left multiplication, $\varphi$ is $G$-equivariant. So
%$\varphi$ is a $G$-equivariant dominant morphism. 
Since $G$ acts transitively on $G/L$ and 
$\varphi^{-1}(\varphi(1L))=P/L$, by Remark \ref{rem1.20}(i)
and Theorem \ref{thm 1.10}, we have
$\kappa_G(G/L)=\kappa_G(G/P)+\kappa_H(P/L)$, where $H$ is a
subgroup of $G$ of the form $B\cap P$, where $B$ is a Borel
subgroup of $G$ and $\dim B\cap P$ is minimal. 
By Lemma \ref{lem 5.10}, 
$B\cap P$ is a Borel subgroup of a Levi subgroup of $P$. The
fact that all Levi subgroups of $P$ are conjugate 
implies that we may assume that
$B\cap P$ is a Borel subgroup of $L$, so
$\kappa_H(P/L)=\kappa_L(P/L)$. Finally, by Lemma \ref{lem 1.30}, 
$\kappa_G(G/P)=0$. It follows that
$\kappa_G(G/L)=\kappa_L(P/L)$.
\end{comment}
%Let $\tau : P\to P/L$ and $\zeta : P \to R_u(P)$ 
%defined by $\psi(x)=\psi(zy)=z$, where $x = zy$ with  $z \in R_u(P)$ and $y \in L$. 
%be the canonical map and projection, respectively.
%The fibres of $\tau$ are of the form $zL$, where $z \in R_u(P)$.
%Since $\zeta$ is constant on the fibres of $\tau$, 

If we write $P = R_u(P)L$, then the bijection
$P/L = R_u(P)L/L \cong R_u(P)$ gives
a canonical $L$-equivariant isomorphism 
$\phi : P/L \to R_u(P)$ defined by
$\phi(xL) = y$, where $x = yz$ with $y\in R_u(P)$ and $z \in L$.
%by Theorem \ref{UMP}. 
%Since $\varphi(lxL)=\varphi(lzyL)=\varphi(lzl^{-1}lyL)=lzl^{-1}$, where
%$z \in R_u(P)$ and $l,y \in L$, we see that $\varphi$ is an
%One checks that $\phi$ is $L$-equivariant (and clearly bijective). 
Thus, we have $\kappa_L(P/L) = \kappa_L(R_u(P))$. %, by Theorem \ref{thm 1.10}.

A Springer isomorphism between 
the unipotent variety of $G$ and $\N$
restricts to an $L$-equivariant
isomorphism $R_u(P) \to \pp_u$, e.g., see \cite[Cor.\ 1.4]{Go3}, 
so that $\kappa_L(R_u(P))=\kappa_L(\pp_u)$. %, by Theorem \ref{thm 1.10}.
\end{proof}

%\begin{cor} If $\htt(e)=2$ and $\lambda$ induces an even grading
%on $\gg$, then $\kappa_G(e)=\kappa_G(G/L)$.\end{cor}

\begin{rems}
\label{rema 5.20}
(i). While the first
two equalities of Proposition \ref{prop 6.10} hold in 
arbitrary characteristic, the third equality requires the characteristic
of the underlying field to be zero or a good prime for $G$;
this assumption is required for the existence of a Springer isomorphism,
cf.\ \cite[Cor.\ 1.4]{Go3}.

(ii). Lemma 4.2 in \cite{Br1} states 
that there is a dense $L$-orbit on $G/B$ if
and only if there is a dense $B_L$-orbit on $R_u(P)$, where $B_L$
is a Borel subgroup of $L$. Notice that there is a dense $L$-orbit
on $G/B$ if and only if there is a dense $B$-orbit on $G/L$.
In other words, $\kappa_G(G/L)=0$ if and only
if $\kappa_L(R_u(P))=0$. Thus, Proposition \ref{prop 6.10}
generalizes \cite[Lem.\ 4.2]{Br1}.
\end{rems}

By Proposition \ref{prop 6.10}, 
the problem of determining $\kappa_L(R_u(P))$ is equivalent
to the problem of determining $\kappa_G(G/L)$. In particular, a
Borel subgroup of $L$ acts on $R_u(P)$ with a dense orbit if and
only if $L$ is a spherical subgroup of $G$. In fact, the latter
have been classified: In characteristic zero this result was
proved by M.\  Kr\"{a}mer in \cite{Kr} and extended to arbitrary
characteristic by J.\ Brundan in \cite[Thm.\ 4.1]{Br1}:

\begin{thm}
\label{thm 6.10} 
Let $L$ be a proper Levi subgroup of a simple
group $G$. Then $L$ is spherical in $G$ if and only
if $(G, \D L)$ is one of
$(A_n, A_{i-1}A_{n-i})$, $(B_n, B_{n-1})$, $(B_n, A_{n-1})$, $(C_n,
C_{n-1})$, $(C_n, A_{n-1})$, $(D_n, D_{n-1})$, $(D_n, A_{n-1})$,
$(E_6, D_5)$,  or $(E_7, E_6)$.
\end{thm}

%If we assume $G$ to be simple in Theorem \ref{thm 5.30}, then $L$
%is spherical if and only if $L=G$ or $(G, L')$ is one of the
%following
% \[(A_n, A_{m}A_{n-m-1}), (B_n, B_{n-1}), (B_n, A_{n-1}), (C_n,
%C_{n-1}),\] \[(C_n, A_{n-1}), (D_n, D_{n-1}), (D_n, A_{n-1}),
%(E_6, D_5) \text{ or } (E_7, E_6).\]
\begin{comment}\begin{rem}\label{rem 6.10}Note if $G$ is of type $E_8$,
$F_4$ or $G_2$, then $G$ has no proper spherical Levi subgroups.
Also note that if $L$ is a proper spherical Levi subgroup of $G$,
then $\rk \D L=\rk G -1$.\end{rem}\end{comment}

We also recall the classification of the parabolic subgroups of $G$ with
an abelian unipotent radical, cf.\ \cite[Lem.\ 2.2]{RiRoSt}.

\begin{lem}
\label{lem 6.20} 
Let $G$ be a simple algebraic group and $P$ be a parabolic
subgroup of $G$. Then $R_u(P)$ is
abelian if and only if $P$ is a maximal parabolic subgroup
of $G$ which is conjugate to the standard parabolic subgroup $P_I$
of $G$, where $I = \Pi\setminus\{\alpha\}$ and $\alpha$
occurs in the highest root $\tilde\alpha$ with coefficient $1$.
\end{lem}
% The Borel subgroup $B$ of $G$ in this
%theorem is the Borel subgroup of $G$ such that $B=P_{\varnothing}$.

Let $\Pi=\{\alpha_1,\alpha_2, \ldots, \alpha_n\}$ be a set of
simple roots of the root system $\Psi$ of $G$.
Using Lemma \ref{lem 6.20}, we can readily determine the standard
parabolic subgroups $P_I$ of $G$ with an abelian unipotent radical.
For $G$ simple we gather this information in Table \ref{Tab 6.1} below 
along with the structure of the corresponding standard Levi subgroup 
$L_I$ of $P_I$.
Set $P_{\alpha_{i}'} = P_{\Pi\setminus \{\alpha_i\}}$.
Here the simple roots are labelled as in \cite[Planches I - IX]{Bo2}. 

\begin{table}[ht]
\renewcommand{\arraystretch}{1.2}
  \centering
\begin{tabular}{|c|c|c|}
  \hline
  % after \\: \hline or \cline{col1-col2} \cline{col3-col4} ...
  Type of $G$ & $P_I$ & Type of $\D L_I$ \\
    \hline
 $A_n$ & $P_{\alpha_{i}'}$ for $1\leqslant i \leqslant n$ & $A_{i-1}A_{n-i}$ \\
  $B_n$ & $P_{\alpha_{1}'}$                         & $B_{n-1}$ \\
  $C_n$ & $P_{\alpha_{n}'}$                         & $A_{n-1}$ \\
  $D_n$ & $P_{\alpha_{1}'},P_{\alpha_{n-1}'}$ and $P_{\alpha_{n}'}$ & $D_{n-1}$ or $A_{n-1}$ \\
  $E_6$ & $P_{\alpha_{1}'}$ and $P_{\alpha_{6}'}$             & $D_5$ \\
  $E_7$ & $P_{\alpha_{7}'}$                                   & $E_6$ \\
%  $E_8$ &-                                                    & - \\
%  $F_4$ & -                                                   & - \\
%  $G_2$ & -                                                   & - \\
  \hline
\end{tabular}
\medskip
  \caption{Parabolic Subgroups with Abelian Unipotent Radical.}
\label{Tab 6.1}
\end{table}

Note that if $G$ is of type $E_8$, $F_4$ or $G_2$, then $G$
does not admit a parabolic subgroup with an abelian unipotent
radical. Also compare the list of pairs $(G,\D L)$ from Table
\ref{Tab 6.1} with the list in Theorem \ref{thm 6.10}.

Our next result is immediate from \cite[Thm.\ 4.1, Lem.\ 4.2]{Br1}.

\begin{prop}
\label{prop 6.20} 
If $P=LR_u(P)$ is a parabolic subgroup of $G$
with $R_u(P)$ abelian, then $\kappa_L(R_u(P))=0$.
\end{prop}

\begin{proof}
If $R_u(P)$ is abelian, then using Table \ref{Tab 6.1} we see that
all the possible pairs $(G,\D L)$ appear in the list of spherical
Levi subgroups given in Theorem \ref{thm 6.10}, that is
$\kappa_G(G/L)=0$. Proposition \ref{prop 6.10} then implies that
$\kappa_L(R_u(P))=0$. 
\end{proof}

\begin{cor}
\label{cor 6.10}
If $P$ is a parabolic subgroup of $G$
with $R_u(P)$ abelian, then
$\kappa_L(\pp_u)=0$.
\end{cor}

%\begin{proof} Simply use Proposition\ref{prop 6.10}. \end{proof}

Let $\Psi$ be the root system of $G$ 
and let $\Pi\subseteq \Psi$ be a set of simple
roots of $\Psi$. 
Let $P = P_I$ ($I\subseteq \Pi$) be a standard parabolic subgroup of $G$.
Let $\Psi_I$ be the
root system of the standard Levi subgroup $L_I$,
i.e., $\Psi_I$ is spanned by $I$.
Define $\Psi_I^+=\Psi_I\cap \Psi^+$. For any root
$\alpha \in \Psi$ we can uniquely write
$\alpha=\alpha_I+\alpha_{I'}$ where $\alpha_I=\sum_{\beta
\in I} c_{\beta}\beta$ and $\alpha_{I'}=\sum_{\beta \in
\Pi\setminus I} d_{\beta}\beta$. 
We define the \emph{level of} $\alpha$ 
(\emph{relative to $P$} or \emph{relative to $I$}) to be
\[
\lv(\alpha) := \sum\limits_{\beta \in \Pi\setminus I} d_{\beta},
\]
cf.\ \cite{AzBaSe}.
%The map $\lv : \Psi \to \mathbb{Z}$ is linear. 
Let $d$ be the maximal level of any root in $\Psi$. If
$2i> d$, then
\[
A_i := \prod\limits_{\lv(\alpha)=i}U_{\alpha}
\] 
is an abelian unipotent subgroup of $G$. Note $A_d$ is the centre of $R_u(P)$.
Since $L$ normalizes each $A_i$, we can consider $\kappa_L(A_i)$.
\begin{comment}
It is sufficient to show that $(U_{\alpha},U_{\beta})=\{e\}$ for
all $\alpha,\beta \in \Psi$ with $\lv(\alpha)=\lv(\beta)=i$. The
Chevalley commutator relations imply that
$(U_{\alpha},U_{\beta})\subseteq\prod_{l,j>0\:;\:l\alpha+j\beta
\in \Psi} U_{l\alpha+j\beta}$, $(U_{\alpha},U_{\beta})=\{e\}$.
Since $\lv(l\alpha+j\beta)=l\lv(\alpha)+j\lv(\beta)\geqslant
2i>d$, we see that $l\alpha+j\beta$ is not a root for any $l,j>0$.
Note, if $i > d=\htt(e)$ then $A_i$ is trivial.
\end{comment}

\begin{prop}
\label{prop 6.30}
If $P$ is a parabolic subgroup of
$G$ and $2i>d$, then $\kappa_L(A_i)=0$.
\end{prop}

\begin{proof}
We maintain the setup from the previous paragraph.
Setting $A_i=\prod_{\lv(\alpha)=i}U_{\alpha}$ and
$A_i^-=\prod_{\lv(\alpha)=-i}U_{\alpha}$, let 
$H$ be the subgroup of $G$ generated by $A_i$, $A_i^-$, and $L$. 
Then $H$ is reductive, with root system 
$\Psi_I\cup\{\alpha\in\Psi \mid \lv(\alpha)=\pm i\}$, 
and $LA_i$ is a parabolic subgroup of $H$. Since $A_i$ is
abelian, we can invoke Proposition \ref{prop 6.20} to deduce
that $\kappa_L(A_i)=0$.
 \end{proof}

There is a natural Lie algebra analogue of Proposition \ref{prop 6.30}: 
Maintaining the setup from above, for $2i> d$, we see that
$\aaa_i := \bigoplus_{\lv(\alpha)=i}\gg_{\alpha}$
is an abelian subalgebra of $\gg$. 
Since $\Lie U_{\alpha} = \gg_{\alpha}$ for all $\alpha \in \Psi$, 
we have $\Lie A_i =\aaa_i$. 
Thanks to \cite[Cor.\ 1.4]{Go3}, we obtain the following consequence
of Proposition \ref{prop 6.30}.

\begin{comment} Since the characteristic of $k$ is good for $G$ we have that
there are no degenerations in the Chevalley commutator relations,
so we have that $\Lie Z(R_u(P))=Z(\nn)$.\end{comment}

\begin{cor}
\label{cor 6.35}
If $P$ is a parabolic subgroup of
$G$ and $2i>d$, then $\kappa_L(\aaa_i)=0$.
\end{cor}

\begin{comment}\begin{cor}\label{cor 6.30}If $P=LR_u(P)$ is a parabolic
subgroup of $G$, then $\kappa_L(Z(\nn))=0$.\end{cor}\end{comment}

\begin{rems}
(i).
Corollary \ref{cor 6.35} was first proved, for a field of
characteristic zero, in \cite[Prop.\ 3.2]{Pa3}, although the proof
there is somewhat different from ours.

(ii).
Propositions \ref{prop 6.20} and \ref{prop 6.30} suggest that 
that if $A$ is an abelian subgroup of $R_u(P)$ which is
normal in $P$, then $\kappa_L(A)=0$. It is indeed the case that $P$
acts on $A$ with a dense orbit, see \cite[Thm.\ 1.1]{Ro}.
However, this is not the case when we consider 
instead the action of a Borel subgroup of
a Levi subgroup of $P$ on $A$. 
For example, it follows from  \cite[Table 1]{Ro} 
that if $G$ is of type $A_n$, then the dimension of a maximal normal
abelian subgroup $A$ of a Borel subgroup $B$ of $G$ is $i(n+1-i)$,
where $1\leqslant i \leqslant n$. Clearly, for $1\neq i\neq n$ we
have $\dim A > \rk G$. Thus, a maximal torus of $B$ cannot act on 
$A$ with a dense orbit. 
Using \cite[Table 1]{Ro}, it is easy to construct further examples. 
\end{rems}

\section{The Classification of the Spherical Nilpotent Orbits}
\label{sect:classification}
\subsection{Height Two Nilpotent Orbits}
\label{sub:ht2}
In this subsection we show that height two nilpotent orbits are
spherical. 
Let $e \in \gg$ be nilpotent and let $\lambda \in \Omega_G^a(e)$
be an associated cocharacter of $G$.
Define the following subalgebra of $\gg$:
\begin{equation}
\label{eq 6.1}
\gg_E := \bigoplus_{i \in \mathbb{Z}}\gg(2i).
\end{equation}

%Next we prove some elementary properties of the subalgebra $\gg_E$.

\begin{prop}
\label{prop 6.40}
Let $e \in \N$, $\lambda \in \Omega_G^a(e)$, and let 
$\gg_E$ be the subalgebra of $\gg$ defined in \eqref{eq 6.1}. 
\begin{itemize}
\item[(i)] 
There exists a connected reductive subgroup $G_E$ of
$G$ such that $\Lie G_E=\gg_E$.
\item[(ii)] 
There exists a parabolic subgroup $Q$ of $G_E$ such that $\Lie
Q=\bigoplus_{i\geqslant 0}\gg(2i)$. Moreover, 
$C_G(\lambda)$ is a Levi subgroup of $Q$ 
%the parabolic subgroup
%$Q$ has Levi decomposition $Q=C_G(\lambda)R_u(Q)$ 
and $\Lie R_u(Q)=\bigoplus_{i \geqslant 1}\gg(2i)$.
\end{itemize}
\end{prop}

\begin{proof}
Fix a maximal torus $T$ of $G$ such that $\lambda(k^*) \leqslant T$. 
Set $\Phi=\{\alpha \in\Psi \mid \langle\alpha,\lambda\rangle  \in 2\ZZ\}$. 
Then $\gg_E = \bigoplus_{\alpha \in \Phi} \gg_{\alpha}$.

Then $\Phi$ is a semisimple subsystem of $\Psi$.
The subgroup $G_E$ generated by $T$ and
all the one-dimensional root subgroups $U_{\alpha}$ with 
$\alpha \in \Phi$ is reductive and has Lie algebra $\gg_E$.

Let $Q=P\cap G_E$, where $P = P_\lambda$.
Since $\lambda(k^*)\leqslant T\leqslant G_E$, we see that
$Q$ is a parabolic subgroup of $G_E$, see the remarks preceding
Theorem \ref{thm 4.15}. Since $\Lie C_G(\lambda)=\gg(0)$, we have
$C_G(\lambda)\leqslant Q$ and so $C_G(\lambda)$ is a Levi subgroup
of $Q$. The remaining claims follow from the fact that $\Lie
P=\gg_{\geqslant 0}$, the parabolic subgroup $P$ has Levi
decomposition $P=C_G(\lambda)R_u(P)$ and $\Lie R_u(P)=\gg_{>0}$.
\end{proof}

The following discussion and Lemma \ref{lem 6.500} allow us to reduce the determination 
of the spherical nilpotent orbits to the case when $G$ is simple.
Since the centre of $G$ acts trivially on $\gg$, we may assume that $G$ is semisimple.
Let $\tilde G$ be semisimple of adjoint type and $\pi : G \to \tilde G$ be the 
corresponding isogeny.
Let  $e \in \gg$ be nilpotent and let  $\tilde e = d\pi_1(e)$.
Consider the restriction of $d\pi_1$ to the nilpotent variety
of $\gg$. Then $d\pi_1 : \N \to \tilde \N$ is a dominant 
$G$-equivariant morphism, where $\tilde \N$ is the 
nilpotent variety of $\Lie \tilde G$ and $G$ acts on $\tilde \N$ via $\tilde \Ad \circ \pi$.
It then follows from Theorem \ref{thm 1.10} that 
$\kappa_G(G\cdot e) = %\kappa_G(G \cdot \tilde e) = 
\kappa_{\tilde G}(\tilde G \cdot \tilde e)$.
We therefore may assume that $G$ is semisimple of adjoint type.

\begin{lem}
\label{lem 6.500} 
Let $G$ be semisimple of adjoint type. Then $G$ is a direct product of simple groups
$G = G_1G_2 \cdots G_r$. If $e \in \gg$ is nilpotent, then $e = e_1 + e_2 + \ldots + e_r$
for $e_i$ nilpotent in $\gg_i = \Lie G_i$ and $\kappa_G(G\cdot e) = \sum_{i=1}^r \kappa_{G_i}(G_i \cdot e_i)$.
\end{lem}

\begin{proof} 
Since $G$ is semisimple of adjoint type, so that $G$ is the direct product 
$G = G_1G_2 \cdots G_r$ of simple groups $G_i$,
we have $\Lie G = \oplus \Lie G_i$.
Let $e \in \gg$ be nilpotent. 
Clearly, any element $x \in C_G(e)$ is of the form $x = x_1 x_2 \cdots x_r$ where
$x_i \in G_i$ and we also have that $e = e_1 + e_2 + \ldots +e_r$, where $e_i \in \gg_i$ 
and each $e_i$ must be nilpotent. 
We know that $\Ad(x)(e) = e$ so 
$\Ad(x_1)\Ad(x_2) \cdots  \Ad(x_r)(e_1 + e_2 + \ldots + e_r) = e_1 + e_2 + \ldots + e_r$. 
For $i \ne j$ we have $\Ad(x_i)(e_j) = e_j$, so $\Ad(x)(e_i) = \Ad(x_i)(e_i)$.
Therefore, as $\Ad(x_i)$ stabilizes $\gg_i$, we have $\Ad(x_i)(e_i) = e_i$. 
Thus, we obtain the following
decomposition $C_G(e) = C_{G_1}(e_1)C_{G_2}(e_2) \cdots C_{G_r}(e_r)$. For $B$ a Borel subgroup of $G$
we have $B = B_1 B_2 \cdots B_r$, where each $B_i$ is a Borel subgroup of $G_i$
and $C_B(e) = C_{B_1}(e_1)C_{B_2}(e_2) \cdots C_{B_r}(e_r)$. In particular, for
$B \in \Gamma_G(e)$ we have that $\dim C_B(e)$ is minimal. 
This implies that $\dim C_{B_i}(e_i)$ is minimal for each $i$
and so $B_i \in \Gamma_{G_i}(e_i)$. 
Therefore, we have
%Using that $B \in \Gamma_G(e)$ and that 
%$B_i \in \Gamma_{G_i}(e_i)$, we have 
\begin{align*}
\kappa_G(G \cdot e) & = \dim G -  \dim C_G(e) - \dim B + \dim C_B(e) \\
& = \sum_{i=1}^r \dim G_i  - \sum_{i=1}^r \dim C_{G_i}(e_i) - \sum_{i=1}^r \dim B_i + \sum_{i=1}^r \dim C_{B_i}(e_i)\\
& = \sum_{i=1}^r (\dim G_i - \dim C_{G_i}(e_i) - \dim B_i + \dim C_{B_i}(e_i))\\
& = \sum_{i=1}^r \kappa_{G_i}(G_i \cdot e_i),
\end{align*}
and the result follows. 
\end{proof}

\begin{lem}
\label{lem 6.50} 
Let $G$ be a connected reductive algebraic
group and $e \in \gg$ be nilpotent. If $\htt(e)=2$, then $e$ is
spherical.
\end{lem}

\begin{proof} 
First we assume that $G$ is simple.
Let $\lambda \in \Omega_G^a(e)$. 
Let $\gg_E$ be the Lie subalgebra of $\gg$ as defined
in \eqref{eq 6.1} and let $Q$ be the parabolic subgroup of $G_E$ as
in Proposition \ref{prop 6.40}(ii). Since $\htt(e)=2$, we have
$\gg_E=\gg(-2)\bigoplus \gg(0)\bigoplus \gg(2)$. Set $L=C_G(\lambda)$.
Then $\kappa_G(G\cdot e)=\kappa_L(\gg(2))$, by Theorem \ref{thm
5.10}. Also, by Proposition \ref{prop 6.40}, $\Lie R_u(Q)=\gg(2)$.
Since $R_u(Q)$ is abelian, Corollary \ref{cor 6.10} implies that
$\kappa_L(\gg(2))=0$.

Now suppose that $G$ is reductive. Let $\D
G=G_1G_2\cdots G_r$ be a commuting product of simple groups. 
For $e \in \gg$ we have $e=e_1+e_2+\ldots +e_r$,
where $e_i \in \gg_i = \Lie G_i$ and each $e_i$ is nilpotent. 
Since $\htt(e)=\max_{1\leqslant i \leqslant r}\htt(e_i)$, 
we have $\htt(e_i)\leqslant \htt(e)=2$ for all $i$. 
Since $\kappa_G(G\cdot e) = \sum_{i=1}^r \kappa_{G_i}(G_i\cdot e_i)$,
by Lemma \ref{lem 6.500}, 
the result follows from the simple case just proved.
\end{proof}

\subsection{Even Gradings}
\label{sub:even}
Suppose that the given nilpotent element $e\in\gg$
satisfies $\htt(e)\geqslant 4$. Also assume
that any $\lambda \in \Omega_G^a(e)$ induces
an \emph{even grading} on $\gg$, that is $\gg(i,\lambda)=\{0\}$
whenever $i$ is odd. As usual we denote
$\gg(i,\lambda)$ simply by $\gg(i)$.

\begin{lem}
\label{lem 7.10} 
Let $e \in \N$ and $\lambda \in \Omega_G^a(e)$ be as above.
Then $\gg_{\geqslant 2}$ is non-abelian.
\end{lem}

\begin{proof} 
Set $\htt(e)=d$. 
%By Lemma \ref{lem 6.05},
For the highest root $\tilde\alpha\in \Psi^+$ we have $\gg_{\tilde\alpha}\subseteq
\gg(d)$. Write $\tilde\alpha=\alpha_1+\alpha_2+\ldots +\alpha_r$ as a sum
of not necessarily distinct simple roots. The sequence of 
simple roots $\alpha_1,\alpha_2,\ldots ,\alpha_r$ can be
chosen so that $\alpha_1+\alpha_2+\ldots +\alpha_s$ is a root for
all $1\leqslant s \leqslant r$, \cite[Cor.\ 10.2.A]{Hu1}. 
Since the grading of $\gg$ induced by $\lambda$ is even, 
for all simple roots $\alpha \in \Pi$, we have 
$\gg_\alpha\subseteq\gg(i)$ with $i \in \{0,2\}$, 
cf.\ \eqref{eq 6.2}. Since $d \geqslant 4$, for at least one $\alpha_i$ we
must have $\gg_{\alpha_i}\subseteq \gg(2)$. Let $\alpha_k$ be the
last simple root in the sequence $\alpha_1,\alpha_2,\ldots
,\alpha_r$ with this property. Thus, for
$\beta=\alpha_1+\alpha_2+\ldots +\alpha_{k-1}$ we have
$\gg_{\beta} \in \gg(d-2)\subseteq\gg_{\geqslant 2}$. 
%Thanks to the Chevalley commutator relations and the fact that 
Since $\cha k$ is good for $G$, we have 
$[\gg_{\beta},\gg_{\alpha_k}]=\gg_{\beta'}$
where $\beta'=\beta +\alpha_{k}$. Therefore,
$\gg_{\geqslant 2}$ is non-abelian. 
\end{proof}

\begin{cor}
\label{cor}
Let $P$ be the destabilizing parabolic subgroup
of $G$ defined by $e \in \N$. Then $R_u(P)$ is
non-abelian.
\end{cor}

\begin{comment}
\begin{proof} Since $\lambda$ induces an even grading on $\gg$, we
have that $\Lie R_u(P)=\gg_{\geqslant 2}$, see Section \ref{sec
3.1}. Since $\gg_{\geqslant 2}$ is non-abelian we have that
$R_u(P)$ is also non-abelian: The Chevalley commutator relations
and the fact that $\cha k$ is good for $G$ imply that, for
$\alpha_k$ and $\beta$ as in Lemma \ref{lem 7.10}, we have
$(U_{\alpha_k},U_{\beta})\neq \{e\}$. \end{proof}\end{comment}

Set $\pp_u=\Lie R_u(P)$. Because the grading of $\gg$ is even, 
$\gg_{\geqslant 2}=\pp_u$. Thus, by Proposition
\ref{prop 6.10} and Theorem \ref{thm 5.10}, we have
$\kappa_G(G\cdot e)=\kappa_G(G/L)$, where $L=C_G(\lambda)$. Using
the classification of the spherical Levi subgroups and 
the classification of the parabolic subgroups of $G$ with abelian
unipotent radical, Theorem \ref{thm 6.10} and Lemma \ref{lem 6.20}, we
see that there are only two cases, for $G$ simple, when $R_u(P)$
is non-abelian and $L$ is spherical, namely when $G$ is of type
$B_n$ and $\D L$ is of type $A_{n-1}$ and when $G$ is of type
$C_n$ and $\D L$ is of type $C_{n-1}$.

%Lemma \ref{Lem 7.20} rules out these two cases.

\begin{lem}
\label{lem 7.15}
Let $G$ be of type $B_n$ or of type $C_n$.
Let $e \in \N$ and $\lambda \in \Omega_G^a(e)$.
Set $L=C_G(\lambda)$. If $\pi_e=[1^{r_1},2^{r_2},\ldots]$ is the
corresponding partition for $e$, then 
$\dim Z(L)=|\{a_i,b_i \in \ZZ_{\geqslant 0} \mid
a_i+1=\sum_{j\geqslant i}r_{2j+1}, 
b_i+1=\sum_{j\geqslant i}r_{2j}\}|$.
\end{lem}

\begin{proof}
Since $L$ is reductive, $L=Z(L)\D L$, and 
$Z(L) \cap \D L$ is finite, we have 
$\dim L = \dim Z(L) + \dim \D L$. The result
follows from Proposition \ref{prop 5.05}. 
\end{proof}

It is straightforward to deduce the following from 
Propositions \ref{prop 4.50} and \ref{prop 5.05}.

\begin{lem}
\label{Lem 7.20}
Let $e \in \N$ and $\lambda \in \Omega_G^a(e)$
with $\htt(e)\geqslant 4$. Set $L=C_G(\lambda)$. 
If $G$ is of type $B_n$, then $\D L$ is not of type $A_{n-1}$ and if 
$G$ is of type $C_n$, then $\D L$ is not of type $C_{n-1}$.
\end{lem}

\begin{lem}
\label{lem 7.100}
Let $e \in \N$ and suppose that $\lambda \in \Omega_G^a(e)$ 
induces an even grading on $\gg$. If $\htt(e)\geqslant 4$, then
$e$ is non-spherical.
\end{lem}

\begin{proof}
First we observe that
if $G$ is simple, then the statement follows from the facts that
$R_u(P)$ is non-abelian (Corollary \ref{cor}) and that $(G,\D
L)$ is not one of the pairs $(B_n,A_{n-1})$ or $(C_n,C_{n-1})$
(Lemma \ref{Lem 7.20}). So by Theorem \ref{thm 6.10} and
Lemma \ref{lem 6.20}, we see that $L$ is a non-spherical subgroup.
Therefore, by Proposition \ref{prop 6.10},
$\kappa_L(\gg_{\geqslant 2})>0$ and $e$ is non-spherical.

In case $G$ is reductive, we argue as in the proof of Lemma \ref{lem 6.50} 
and reduce to the simple case.
%Let $\D G=G_1G_2\cdots G_r$ be a commuting product of simple
%groups. Then $e=e_1+e_2+\ldots +e_r$,
%where $e_i \in \gg_i$ and each $e_i$ is nilpotent. We have
%$\htt(e)=\max_{1\leqslant i \leqslant r}\htt(e_i)$.
%In particular, $\htt(e_i)= \htt(e)\geqslant 4$ for some $i$. Since
%$\kappa_G(G\cdot e) = \sum_{i=1}^r \kappa_{G_i}(G_i\cdot e_i)$,
%the result follows from the simple case considered above.
\end{proof}

\subsection{Nilpotent Orbits of Height at Least Four}
\label{sub:ht4}
Let $e \in \gg$ be nilpotent and let $\lambda \in \Omega_G^a(e)$.
Let $\gg_E$ be the subalgebra of $\gg$ as defined in \eqref{eq 6.1}. 
Also let $G_E$ be the connected reductive algebraic group
such that $\Lie G_E=\gg_E$ and $Q$ be the parabolic subgroup of
$G_E$ as in Proposition \ref{prop 6.40}(ii).

Since $e \in \gg_E$ and $\lambda(k^*) \leqslant G_E$, 
it follows from \cite[Thm.\  1.1]{FoRo} that $\lambda$ is a
cocharacter of $G_E$ which is  associated to $e$, i.e.\ 
$\lambda \in \Omega_{G_E}^a(e)$.
%Thus the grading $\gg_E=\bigoplus_{i \in \mathbb{Z}}\gg(2i,\lambda)$ is
%induced by an associated cocharacter for $e$ in $\gg_E$. 
Moreover, for $P = P_\lambda$, we have $Q = P\cap G_E$
%, where $P$ is the canonical parabolic subgroup of $G$ for $e$, 
is the destabilizing parabolic subgroup of $G_E$ defined by $e$.

Let $\htt_E(e)$ denote the height of $e \in \gg_E$. 
Now if $\htt(e)\geqslant 4$ and $\htt(e)$ is even,
then $\htt_E(e)=\htt(e)$. The case when $\htt(e)\geqslant 4$ and
$\htt(e)$ is odd is slightly more involved. First we need some
preliminary results. A proof of the following can be found in
\cite[Prop.\ 2.4]{Pa2}.

\begin{lem}
\label{lem 7.30}
Suppose that $\cha k=0$. If $e \in \N$ with 
$\htt(e)$ odd, then the weighted Dynkin diagram $\Delta(e)$
contains no ``2'' labels.
\end{lem} 

If $\Pi$ is a set of simple
roots of $\Psi$ relative to a maximal torus $T$ which contains
$\lambda(k^*)$, then for $\alpha \in \Pi$ we have
\begin{equation}
\label{eq 7.10} 
\gg_{\alpha} \subseteq \gg(i) \text{ where }i\in\{0,1\}.
\end{equation}
To see this recall \eqref{eq 6.2}:
$\Ad(\lambda(t))(e_{\alpha})=t^{m_{\alpha}}e_{\alpha}$, for
$e_{\alpha}\in \gg_{\alpha}$ and $m_{\alpha}$ is the corresponding
label of the weighted Dynkin diagram $\Delta(e)$ of $e$. 
Thus, by Lemma \ref{lem 7.30}, we have $m_{\alpha}\in \{0,1\}$.

%This now allows us to prove the following lemma.
\begin{lem}
\label{lem 7.40}
If $\htt(e) = d$ odd, then $\gg(d-1)\neq\{0\}$.
\end{lem}

\begin{proof} 
The result follows easily, 
arguing as in the proof of Lemma \ref{lem 7.10}  and using 
\eqref{eq 7.10}.
%For the highest root $\tilde\alpha$ of $\Psi^+$ we have $\gg_{\tilde\alpha} \subseteq \gg(d)$.
%Write $\tilde\alpha=\alpha_1+\alpha_2+\ldots +\alpha_r$ as a sum of
%not necessarily distinct simple roots. The sequence of simple roots
%$\alpha_1,\alpha_2,\ldots ,\alpha_r$ can be chosen so that
%$\alpha_1+\alpha_2+\ldots +\alpha_s$ is a root for all 
%$s \leqslant r$. By \eqref{eq 7.10}, 
%$\gg_{\alpha_i}$ either lies in $\gg(0)$ or $\gg(1)$. Therefore,
%for at least one $\alpha_i$ we have $\gg_{\alpha_i}\subseteq \gg(1)$. 
%Let $\alpha_k$ be the last root in the sequence
%$\alpha_1,\alpha_2,\ldots ,\alpha_r$ with this property. Thus, for
%$\beta=\alpha_1+\alpha_2+\ldots +\alpha_{k-1}$ we have
%$\gg_{\beta} \subseteq \gg(d-1)$. 
\end{proof}

\begin{cor}
If $e \in \N$ with $\htt(e)$ odd, then
$\htt_E(e)=\htt(e)-1$. 
\end{cor}

In particular, we have the following conclusion.

\begin{cor}
\label{cor 7.20} 
If $e \in \N$ with 
$\htt(e)\geqslant 4$,  then $\htt_E(e)\geqslant 4$.
\end{cor}

Thus, by Lemma \ref{lem 7.100}, Corollary \ref{cor 7.20},  
and the fact that 
$\Omega_G^a(e) \cap Y(G_E) = \Omega_{G_E}^a(e)$
(\cite[Thm.\  1.1]{FoRo}), 
we have $\kappa_L(\gg_{E,\geqslant 2})>0$, where
$\gg_{E,\geqslant 2}=\bigoplus_{i \geqslant 1}\gg(2i)$ and
$L=C_G(\lambda)=C_{G_E}(\lambda)$.

\begin{lem}
\label{lem 7.50}
If a Borel subgroup $B_L$ of $L$ acts
on $\gg_{\geqslant 2}$ with a dense orbit, then $B_L$ acts on
$\gg_{E,\geqslant 2}$ with a dense orbit.
\end{lem}

\begin{proof} 
This follows readily from Theorem \ref{lem 1.60}. 
%Clearly, $L$ and thus $B_L$ stabilizes each $\gg(i)$. Indeed, if
%$\sum_{i=2}e_i$ is a representative of a dense $B_L$-orbit
%on $\gg_{\geqslant 2}$, where $e_i \in \gg(i)$, then
%$\sum_{i=1}e_{2i}$ is a representative of a dense
%$B_L$-orbit on $\gg_{E,\geqslant 2}$. 
\end{proof}

Combining Lemmas \ref{lem 7.100}, \ref{lem 7.50} and 
Corollary \ref{cor 7.20}, we get the main result of this subsection.

\begin{prop}
\label{prop 7.20}
Let $e \in \N$. If $\htt(e)\geqslant 4$, then
$e$ is non-spherical.
\end{prop}

\begin{comment}
As we commented on in Section \ref{chak} Proposition \ref{prop
7.20} was first proved by D.I.\ Panyushev for a field of
characteristic zero in \cite{Pa3}. The proof in \cite{Pa3} is very
different from ours. See Section \ref{chak} for a review of
Panyushev's proof. We note that our methods are valid in
characteristic zero as well.

\chapter{}\label{C8}
Throughout this chapter $G$ is a connected reductive algebraic
group over an algebraically closed field $k$, the Lie algebra of
$G$ is denoted by $\gg$ and the characteristic of $k$ is good for
$G$. Recall the standard setup from Chapter \ref{C4}. For $e \in
\gg$ a non-zero nilpotent element we have $\lambda\in Y(G)$ an
associated cocharacter for $e$ in $\gg$ and the canonical
parabolic subgroup $P_{\lambda}=P=C_G(\lambda)R_u(P)$. Also recall
the definition of height of a nilpotent element $e$, see
Definition \ref{def 4.40}.

In this chapter we prove that if $\htt(e)=3$, then $e$ is
spherical. In order to prove this we make use of Theorem \ref{thm
5.10} to see that $\kappa_G(G\cdot
e)=\kappa_{C_G(\lambda)}(\gg_{\geqslant 2})$. We also prove this
result using case-by-case arguments. \end{comment}

\subsection{Nilpotent Orbits of Height Three}
Let $e \in \N$ and let $\lambda \in \Omega_G^a(e)$.
Let $P = P(e)$ be the destabilizing parabolic subgroup defined by $e$.
Then $P = LR_u(P)$ for $L = C_G(\lambda)$.
Let $B_L$ be a Borel subgroup of $L$ so that
$\lambda(k^*)\leqslant B_L$. Write $B_L=TU_L$ for a Levi
decomposition of $B_L$, where $U_L = R_u(B_L)$ 
and $T$ is a maximal torus of $G$ containing $\lambda(k^*)$. 
Let $\bb_L = \Lie B_L$,  $\nn = \Lie U_L$, and
$\mathfrak{t}=\Lie T$. 
%Since $\htt(e) = 3$, we have $\gg_{\geqslant 2}=\gg(2)\bigoplus \gg(3)$.

\begin{lem}
\label{lem 8.10}
Let $e\in \gg$ be nilpotent and 
$\lambda$ be an associated cocharacter for $e$ in $\gg$.
Then the following are equivalent.
\begin{itemize} 
\item[(i)] The nilpotent
element $e$ is spherical.
\item[(ii)] There exists $e' \in \gg_{\geqslant 2}$ 
such that $\overline{\Ad(B_L)(e')}=\gg_{\geqslant 2}$.
\item[(iii)] There exists $e' \in \gg_{\geqslant 2}$ 
such that $\dim C_{B_L}(e')=\dim B_L - \dim \gg_{\geqslant 2}$.
\end{itemize}
\end{lem}

\begin{proof} 
Thanks to Theorem \ref{thm 5.10}, 
$\kappa_G(G\cdot e) = \kappa_L(\gg_{\geqslant 2})$. Thus 
(i) and (ii) are equivalent. 
The equivalence between (ii) and (iii) is clear.
%Also for $e' \in
%\gg_{\geqslant 2}$ we have $\dim \Ad(B_L)(e')=\dim B_L-\dim
%C_{B_L}(e')$, the equivalence of $(ii)$ and $(iii)$ follows.
\end{proof}

Recall from Subsection \ref{sub:not} the definition of the support 
of a nilpotent element in $\uu$.

\begin{lem} 
\label{lem 8.20}
Let $e \in \gg_{\geqslant 2}$. 
If $\supp(e)$ is linearly independent, 
then $\dim C_{T}(e) = \dim T - |\supp(e)|$.
\end{lem} 

\begin{proof}
Suppose that $\supp(e)$ is linearly independent. Then 
$\dim \Ad(T)(e) = |\supp(e)|$, e.g.\ see \cite[Lem.\ 3.2]{Go2}. 
The desired equality follows.
\end{proof}

The following is a standard consequence of orbit maps.
%cf.\ \cite[Prop.\ 6.7]{Bor}.

\begin{lem}
\label{lem 8.30} 
Let $e' \in \gg_{\geqslant 2}$. Then
$\dim C_{B_L}(e')\leqslant \dim \cc_{\bb_L}(e')$ and 
$\dim C_{U_L}(e')\leqslant \dim \cc_\nn(e')$.
\end{lem}

\begin{comment}
\begin{proof}
Let $P = P(e')$ be the destabilizing parabolic subgroup of $G$ 
defined by $e'$ and set $V = \D R_u(P)$.
Note that, as $\Lie R_u(P)=\gg_{\geqslant 1}$ and 
$\gg_{\geqslant 2}=[\gg_{\geqslant 1},\gg_{\geqslant 1}]$, we have 
$\Lie V = \gg_{\geqslant 2}$.
Define $\hh = \bb_L\bigoplus\gg_{\geqslant 2}$, thus $e' \in \hh$. The
fact that $e' \in \gg_{\geqslant 2} = \gg(2)\bigoplus \gg(3)$ and
$\gg_{\geqslant 2}$ is abelian implies that
$\cc_{\hh}(e') = \cc_{\bb_L}(e')\bigoplus \gg_{\geqslant 2}$. In particular,
$\dim \cc_{\hh}(e')=\dim \cc_{\bb_L}(e')+\dim \gg_{\geqslant 2}$. Let
$H$ be the connected subgroup of $G$ such that $\Lie H=\hh$, so
$H=B_LV$. Next we claim that
$\Ad(V)(e')=\{e'\}$: First write $e'=e'_2+e'_3$ where $e'_i\in \gg(i)$.
Using \cite[5.10(3)]{Ja1} we see that for any $v \in V$, we have
$\Ad(v)(e'_2+e'_3)=e'_2 + d_{3}$, for some $d_{3}\in \gg(3)$. Also observe that
$\Ad(v^{-1})(e'_2+e'_3)=e'_2-d_{3}$. 
The claim then follows from the fact that $\Ad([V,V])(e') = \{e'\}$. 
Thus, $C_H(e')=C_{B_L}(e')V$. In particular, $\dim
C_H(e')=\dim C_{B_L}(e')+\dim V=\dim C_{B_L}(e')+\dim \gg_{\geqslant
2}$. In general, $\Lie C_H(e') \subseteq \cc_{\hh}(e')$ 
(cf.\ \cite[\S 1.14]{Ca}), so $\dim C_{B_L}(e')+\dim \gg_{\geqslant
2}=\dim C_H(e')\leqslant \dim \cc_{\hh}(e')=\dim \cc_{\bb_L}(e')+\dim
\gg_{\geqslant 2}$. Thus $\dim C_{B_L}(e')\leqslant \dim
\cc_{\bb_L}(e')$. The proof of the second claim is similar, with
$\nn\bigoplus\gg_{\geqslant 2}$ replacing $\hh$. 
\end{proof}
\end{comment}

In \cite[Prop.\ 5.4]{Go4}, 
Goodwin showed that each $U$-orbit in $\uu$ admits a 
unique 
so called \emph{minimal} orbit representative, see \cite[Def.\ 5.3]{Go4}.
(This depends on a suitable choice of an ordering of the positive roots 
compatible with the height function, cf.\ \cite[Def.\ 3.1]{Go4}.)
Moreover, a special case of \cite[Prop.\ 7.7]{Go4} gives  
that for $e$ the minimal representative of its $U$-orbit in $\uu$, we have
$C_B(e) = C_T(e)C_U(e)$. As a consequence, we readily obtain the following.

\begin{lem}
\label{lem 8.40}
Let $e' \in \gg_{\geqslant 2}$. 
Suppose that 
$e'$ is the minimal representative of its $U$-orbit in $\uu$.
Then
$C_{B_L}(e') = C_{T}(e')C_{U_L}(e')$. In particular, 
$\dim C_{B_L}(e') = \dim C_{T}(e')+\dim C_{U_L}(e')$.
\end{lem}

\begin{prop}
\label{prop 8.10}
Let $G$ be a simple algebraic group.
Table \ref{t:1} below gives a complete list of 
the height $3$ nilpotent orbits in $\gg$.
\end{prop}

\begin{proof} 
For the classical groups we use Proposition \ref{prop 4.50}. 
By Remark \ref{rem 3.30}, there are no height $3$ nilpotent orbits
in types $A_n$ and $C_n$. 
Using the tables in \cite[\S 13]{Ca} and \eqref{eq 6.05}, 
one readily determines the desired orbits when $G$ is exceptional.
\end{proof}

In Table \ref{t:1} we either give the partition
or the Bala--Carter label of the corresponding orbit, cf.~\cite[\S 13]{Ca}.

\begin{table}[h]
\renewcommand{\arraystretch}{1.5}
\begin{tabular}{|c|c|}
\hline
Type of $G$ & Orbits \\
\hline
$A_n$ & - \\
$B_n$ & $[1^j,2^{2i},3]$ with $i>0$\\
$C_n$ & - \\
$D_n$ & $[1^j,2^{2i},3]$ with $i>0$\\
\hline
$G_2$ & $\tilde{A_1}$ \\
$F_4$ & $A_1 + \tilde{A_1}$ \\
$E_6$ & $3A_1$ \\
$E_7$ & $(3A_1)'$,  $4A_1$\\
$E_8$ & $3A_1$,  $4A_1$\\
\hline
\end{tabular}
\medskip
\caption{The nilpotent orbits of height $3$.} 
\label{t:1}
\end{table}

\bigskip

In the next three
subsections we concentrate on the height $3$ orbits in types
$B_n$, $D_n$, and the exceptional types, respectively. 

\subsection{Height Three Nilpotent Elements of $\mathfrak{so}_{2n+1}(k)$}
\label{subsect:soodd}
In this subsection let $G$ be of type $B_n$ for $n \geqslant 3$, 
so $\gg= \mathfrak{so}_{2n+1}(k)$. The nilpotent
orbits in $\gg$ are classified by the partitions of $2n+1$ with
even parts occurring with even multiplicity, 
see \cite[Thm.\ 1.6]{Ja1}. By Proposition \ref{prop 4.50}, the
height $3$ nilpotent orbits correspond to partitions of $2n+1$
of the form $\pi_{r,s}=[1^s,2^{2r},3]$, 
where $r\geqslant 1, s\geqslant 0$ and $2r+s+1=n$. 
Denote the corresponding nilpotent
orbit by $\mathcal{O}_{r,s}$ and a representative of such an orbit
by $e_{r,s}$.

\begin{lem}
\label{lem 8.50} 
There are precisely $\left [\frac{n-1}{2} \right ]$ 
distinct height $3$ nilpotent orbits in $\gg$.
\end{lem}

\begin{proof} 
By our comments above, we need to show 
that there are precisely $\left [ \frac{n-1}{2} \right ]$
partitions of $2n+1$ of the form $\pi_{r,s}$. This is equivalent
to finding all partitions of $n-1$ of the form $[1^{s/2},2^r]$. Thus 
$r$ satisfies $1\leqslant r\leqslant \frac{n-1}{2}$.
Since $r$ is an integer, the result follows. 
\end{proof}

Since the number $2r+1$ appears frequently in the sequel, we set
$\widehat{r}=2r+1$. 
%Also recall that we associate to any nilpotent
%orbit a weighted Dynkin diagram. 
Using \cite[\S 13]{Ca}, we readily see that 
that $e_{r,s}$ has the following weighted Dynkin diagram:

\vfill
\eject

\begin{figure}[ht]
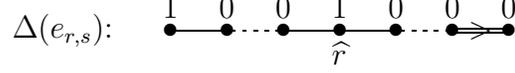

\beginpicture

\setcoordinatesystem units <1.5cm,1.5cm> point at 0 0     % sets scale
\setplotarea x from -.5 to 4, y from 1.5 to 2     % sets plot size up
%

%B n%
\put {$\Delta(e_{r,s})$:} [l] at 2.1 2 \multiput {$\bullet$} at
3.5 2 *2 .5 0 / \put {$>$} at 6.25 2 \multiput {$\bullet$} at 5 2
*3 .5 0 /

\putrule from 3.54 2 to 3.96 2

\putrule from 4.54 2 to 5.54 2

\putrule from 6.04 2.02 to 6.46 2.02

\putrule from 6.03 1.98 to  6.47 1.98

\put {1} at 3.5 2.2 \put {0} at 4 2.2 \put {1} at 5 2.2 \put {0}
at 5.5 2.2 \put {0} at 4.5 2.2

 \put {0} at  6 2.2 \put {0} at 6.5 2.2
 \put {$\widehat{r}$} at 5.0 1.8

\setdashes <.5mm, 1mm>  \putrule from 4.04 2 to 5.96 2

\endpicture 
\caption{Labeling of $\Delta(e_{r,s})$.}
\label{Dynkin3B}

\end{figure}

\begin{rem}
\label{rem 8.10} 
Note that in $\Delta(e_{r,s})$ there are precisely two simple roots, 
$\alpha_1$ and $\alpha_{\widehat{r}}$ that are labeled with a ``1''
and that there is an odd number of simple roots between 
$\alpha_1$ and $\alpha_{\widehat{r}}$.
Also, the short simple root is labeled with a ``1'' 
if and only if $s = 0$, and this can only happen when $n$ is odd.
\end{rem}

We refer to \cite[Planche II]{Bo2} for information regarding the
root system of type $B_n$. 
%We use the following notation for the positive roots in $\Psi^+$
Let $\alpha_1,\ldots,\alpha_n$
be the simple roots of $\Psi^+$ and let 
\begin{align*} 
\beta_{j,k}&=\alpha_j+\ldots +\alpha_k \text{ for }
1\leqslant j \leqslant k\leqslant n, \\ 
\gamma_{j,k}&=\alpha_j+
\ldots +\alpha_{k-1}+2\alpha_k+\ldots + 2\alpha_n \text{ for }
1\leqslant j< k<n,
\end{align*} 
where $\beta_{j,j} = \alpha_j$. Note that all the possible $\beta$'s and
$\gamma$'s exhaust $\Psi^+$.

%\begin{exmp}\label{ex 8.10} The positive roots of $B_4$ are
%\[\Psi^+=\{\alpha_1,\alpha_2,\alpha_3,\alpha_4,\beta_{1,2},\beta_{2,3},
%\beta_{3,4},\beta_{1,3},\beta_{2,4},\beta_{1,4},\gamma_{3,4},\gamma_{2,4},
%\gamma_{1,4},\gamma_{2,3},\gamma_{1,3},\gamma_{1,2}\}.\]Note that
%$\gamma_{1,2}$ is the highest root of $B_4$. In fact,
%$\gamma_{1,2}$ is always the highest root in a root system of type
%$B_n$.
%\end{exmp}

%Let us now look at the structure of the abelian Lie subalgebra
%$\gg_{\geqslant 2}=\gg(2)\bigoplus\gg(3)$.

For a $T$-stable Lie subalgebra $\mm$ of $\uu$ recall the 
definition of the set of roots $\Psi(\mm)$ of $\mm$ with respect to $T$
from Subsection \ref{sub:not}.

\begin{lem}
\label{lem 8.60} 
For an associated cocharacter of $e_{r,s}$ in $\gg$ we have
\begin{itemize}
\item[(i)] $\Psi(\gg(2))=\{ \beta_{1,j},\gamma_{i,m},\gamma_{l,k} 
\mid 1< l < k\leqslant \widehat{r}\leqslant j
\text{ and } 1<i<m\leqslant \widehat{r}\}$, and so
$\dim \gg(2)=2r^2-r+2s+1$;

\item[(ii)] $\Psi(\gg(3))=\{\gamma_{1,k}
\mid  k\leqslant \widehat{r}\}$, and so $\dim \gg(3)=2r$.
\end{itemize}
\end{lem}

\begin{proof} 
%Let $\lambda$ be a cocharacter of $G$ that is  associated to 
%$e_{r,s}$ and let $T$ be a maximal torus of $G$ such that
%$\lambda(k^*)\leqslant T$. 
For every $\delta \in \Psi$ we have
that $\gg_{\delta}\subseteq \gg(i)$ for some $i \in \{0,\pm 1,\pm 2,\pm 3\}$. 
For the simple roots this information can be read off from 
$\Delta(e_{r,s})$, see \eqref{eq 6.2}. 
Let $\delta=\sum_{\alpha\in \Pi}c_{\delta,\alpha}\alpha$ be a
positive root.

Now $\gg_{\delta} \subseteq \gg(2)$ if and only if
$c_{\delta,\alpha_1}+c_{\delta,\alpha_{\widehat{r}}}=2$. 
All of the roots listed above satisfy this
condition, and no others do. 
%A counting exercise gives the dimension.
Finally, $\gg_{\delta} \subseteq \gg(3)$ if and only if
$c_{\delta,\alpha_1}+c_{\delta,\alpha_{\widehat{r}}}=3$. 
All of the roots listed above satisfy this
condition, and no others do. 
%A simple counting exercise gives the dimension.
\end{proof}

%Next we look at the structure of the Lie subalgebra $\bb_L$ of $\gg(0)$.  

\begin{lem}
\label{lem 8.70} 
For an associated cocharacter of $e_{r,s}$ in $\gg$ we have
\begin{itemize}
\item[(i)] $\Psi(\bb_L)=\{\beta_{j,k},\gamma_{l,m} 
\mid \widehat{r} < j \text{ or } 1<j\leqslant k < \widehat{r}, 
\widehat{r}< l<m\}$.
\item[(ii)] $\dim \bb_L=2r^2+s^2+s+r+1$.
\end{itemize}
\end{lem} 

\begin{proof} 
%Let $\lambda$ be a cocharacter of $G$ that is  associated to 
%$e_{r,s}$ and let $T$ be a maximal torus of $G$ such that
%$\lambda(k^*)\leqslant T$. 
For every $\delta \in \Psi$ we have
that $\gg_{\delta}\subseteq \gg(i)$ for some $i \in \{0,\pm 1,\pm
2,\pm 3\}$. 
As mentioned above, for the simple roots this information can be read off from 
$\Delta(e_{r,s})$, see \eqref{eq 6.2}. 
Let $\delta=\sum_{\alpha\in \Pi}c_{\delta,\alpha}\alpha \in \Psi^+$.
Then $\gg_{\delta} \subseteq \bb_L$ if and
only if $c_{\delta,\alpha_1}+c_{\delta,\alpha_{\widehat{r}}}=0$.
All of the roots listed above satisfy this
condition, and no others do. Consequently, 
$\dim \nn = 2r^2+s^2-r$. Since $\dim \mathfrak{t}=n$, we get 
$\dim \bb_L = 2r^2+s^2+s+r+1$. 
%Alternatively we note, thanks to the
%weighted Dynkin diagram, that the Lie algebra $\gg(0)$ is
%isomorphic to $\mathfrak{sl}_{\widehat{r}-1}(k)\times
%\mathfrak{so}_s(k)\times k^2$. The structure and dimension of
%$\bb_L$ follows.
 \end{proof}

%\begin{exmp}
%\label{ex 8.20}
%Let $G$ be of type $B_5$. 
%It follows from Example \ref{eg 6.10} that for
%a nilpotent element $e_{1,2}\in \gg$ we have
%$\Psi(\bb_L) = \{\alpha_2,\alpha_4,\alpha_5,\beta_{4,5},\gamma_{1,2}\}$.
%\end{exmp}

It follows from Figure \ref{Dynkin3B} that $L$ is of 
Dynkin type $A_{\widehat{r}-1} \times B_s$.
Accordingly, there is a natural partition of the roots of $\bb_L$
into a union of two subsets, 
namely the positive roots of the $A_{\widehat{r}-1}$ and $B_s$ subsystems,
respectively.
Thus, we have $\Psi(\bb_L)=\Psi_1(\bb_L)\cup\Psi_2(\bb_L)$, where
\begin{align*}
\Psi_1(\bb_L)&=\{ \beta_{j,k} \mid 1<j\leqslant k < \widehat{r}\}, \\
\Psi_2(\bb_L)&=\{\beta_{j,k},\gamma_{l,m} \mid
\widehat{r}<j\leqslant k , \widehat{r}< l < m\}.
\end{align*}
Similarly, we can decompose the roots of $\gg_{\geqslant 2}$ into two
sets as follows:
$\Psi(\gg_{\geqslant 2}) 
= \Psi_1(\gg_{\geqslant 2}) \cup \Psi_2(\gg_{\geqslant 2})$, where
\begin{align*}
\Psi_1(\gg_{\geqslant 2})&=\{\gamma_{j,k} \mid 1\leqslant j <
k\leqslant \widehat{r} \}, \\
\Psi_2(\gg_{\geqslant 2})&=\{\beta_{1,j},\gamma_{1,k} \mid
\widehat{r} \leqslant j, \widehat{r} < k\}.
\end{align*}

The sets $\Psi_i(\bb_L)$ and $\Psi_i(\gg_{\geqslant 2})$ 
satisfy the following property:
\begin{equation}
\label{eq 8.10}
\delta \in \Psi_i(\bb_L), \eta \in \Psi_{3-i}(\gg_{\geqslant 2}) 
 \ \Rightarrow \ \delta + \eta \notin \Psi,\ i \in \{ 1,2\}.
\end{equation} 
Denote by $\bb_L^i$ the Lie subalgebras of
$\bb_L$ such that $\Psi(\bb_L^i)=\Psi_i(\bb_L)$ for $i = 1,2$. 
For the rest of
this subsection we show that the following element is a
representative of the dense $B_L$-orbit in $\gg_{\geqslant 2}$; set:
\[
e'_{r,s} := 
\sum_{j,k=0}^{r-1}(e_{\gamma_{\widehat{r}-2j-1,\widehat{r}-2j}}
+e_{\gamma_{1,\widehat{r}-2k}})
+e_{\gamma_{1,\widehat{r}+1}}+e_{\beta_{1,\widehat{r}}},\]
where $e_{\delta}\in \gg_{\delta}\setminus\{0\}$ 
for $\delta \in \Psi(\gg_{\geqslant 2})$.

Recall from the paragraph before Lemma \ref{lem 8.40} the 
notion of minimal $U$-orbit representatives in $\uu$ from 
\cite{Go4}. 

\begin{lem}
\label{lem 8.80} 
Each  $e'_{r,s}$
is the minimal representative of its $U$-orbit in $\uu$,
$\supp(e'_{r,s})$ is linearly independent, and  
$|\supp(e'_{r,s})|=\left\{ \begin{array}{ll}
    2r+2 & \text{ if } \:\:s>0; \\
    2r+1 & \text{ if } \:\:s=0. \\
    \end{array} \right.$ 
\end{lem}

\begin{proof} 
It is straightforward to check that $e'_{r,s}$
is the minimal representative of its $U$-orbit in $\uu$ in the sense of
\cite{Go4} and one easily computes $|\supp(e'_{r,s})|$. Note that
the root $\gamma_{1,\widehat{r}+1}$ only occurs if $s > 0$.
%Clearly, we have that $\{\alpha_i\}_{1\leqslant i\leqslant n}$ is
%a basis for $\Psi$. 

Suppose there exist scalars $\tau_j,\xi_k,
\mu$ and $\nu$ such that
\[
\sum_{j=0}^{r-1}\tau_j\gamma_{\widehat{r}-2j-1,\widehat{r}-2j}
+ \sum_{k=0}^{r-1}\xi_k\gamma_{1,\widehat{r}-2k}
+\mu\gamma_{1,\widehat{r}+1}
+\nu\beta_{1,\widehat{r}}=0.
\]
Since the coefficients of $\alpha_1,\alpha_2$,  and $\alpha_3$ must
be zero, we have
\[
\sum_{k=0}^{r-1}\xi_k+\mu+\nu=0,\ 
\tau_{r-1}+\sum_{k=0}^{r-1}\xi_k+\mu+\nu=0, \ \text{ and }\ 
\xi_{r-1}+2\tau_{r-1}+\sum_{k=0}^{r-1}\xi_k+\mu+\nu=0.
\]
These three equations imply that $\tau_{r-1}=0=\xi_{r-1}$. Continuing
in this way, we see that $\tau_j=0=\xi_j$ for all $j$. Thus we
are left to show that 
$\gamma_{1,\widehat{r}+1}$ and $\beta_{1,\widehat{r}}$ are 
linearly independent; but this is obvious. 
\end{proof}

Thanks to Lemma \ref{lem 8.80} it is harmless to assume that $\supp(e'_{r,s})$ is
part of a Chevalley basis of $\gg$.

%For our purposes we only require part of Chevalley's Theorem
%***MOVE TO INTRO***. Assume that $\beta,\gamma\in \Psi$ are
%linearly independent, that is $\beta \neq \pm \gamma$. Let
%$t_{\beta,\gamma} \in \Rn_{0}$ be the largest integer such that
%$\beta-t_{\beta,\gamma}\gamma \in \Psi$. Then:
%\[[e_{\beta},e_{\gamma}]=\left\{
%\begin{array}{ll}
%    \pm(t_{\beta,\gamma}\!+\!1)e_{\beta+\gamma} & \text{ if } \:\:\beta+\gamma\in \Psi; \\
%    0 & \text{ if } \:\:\beta+\gamma\notin \Psi. \\
%    \end{array}
%\right.\]It is well known that $t_{\beta,\gamma}+1\in \{1,2,3\}$
%for all $\beta,\gamma \in \Psi$. A prime $p$ is called
%\emph{very bad for G} if $p$ divides $t_{\beta,\gamma}+1$ for
%some $\beta,\gamma\in \Psi$.
%
%\begin{table}[ht]
%  \centering
%
%
%\begin{tabular}{|c|c||c|c|}
%  \hline
%  % after \\: \hline or \cline{col1-col2} \cline{col3-col4} ...
%  Type of $G$ & Very Bad Primes & Type of $G$ & Very Bad Primes \\
%  \hline
%  $A_n$ & None    & $F_4$   & $2$ \\
%  $B_n$ & $2$     & $e_6$ & None \\
%  $C_n$ & $2$     & $e_7$ & None \\
%  $D_n$ & None     & $e_8$ & None \\
%  $G_2$ & $2,3$   & - & - \\
%
%  \hline
%\end{tabular} \caption{Very Bad Primes}\label{tab B1}
%\end{table}
%Note that if a prime $p$ is good for $G$, then it is not very bad
%for $G$.

\begin{lem}
\label{lem 8.111} 
$\dim \cc_\nn(e'_{r,s})=\left\{
\begin{array}{ll}
    (s-1)^2 & \text{ if } \:\:s>0; \\
    0 & \text{ if } \:\:s=0. \\
    \end{array}
\right.$
\end{lem}

\begin{proof}
Thanks to \eqref{eq 8.10}, we may consider the two summands
$\sum_{j,k=0}^{r-1}(e_{\gamma_{\widehat{r}-2j-1,\widehat{r}-2j}}+e_{\gamma_{1,\widehat{r}-2k}})$
and $e_{\gamma_{1,\widehat{r}+1}}+e_{\beta_{1,\widehat{r}}}$ of
$e'_{r,s}$ separately. Since
$\gamma_{\widehat{r}-2j-1,\widehat{r}-2j}+\gamma_{1,\widehat{r}-2k}
\in \Psi_1(\gg_{\geqslant 2})$, we need only consider the root
spaces $\gg_{\delta}$ for $\delta \in \Psi_1(\bb_L)$. So let
$\beta_{i,m}\in \Psi_1(\bb_L)$. If $m=\widehat{r}-2l$ for some
$0\leqslant l<r$, then, by the Chevalley commutator relations,
$[e_{\gamma_{\widehat{r}-2l+1,\widehat{r}-2(l-1)}},
\gg_{\beta_{i,\widehat{r}-2l}}]=\gg_{\gamma_{i,\widehat{r}-2(l-1)}}$,
since $\cha k$ is good for $G$.
If $m=\widehat{r}-2l-1$ for some $0\leqslant l<r$, then
%by the Chevalley commutator relations
$[e_{\gamma_{1,\widehat{r}-2l}},\gg_{\beta_{i,\widehat{r}-2l-1}}]=\gg_{\gamma_{1,i}}$.
Next we observe that all the $\beta$'s above exhaust the set
$\Psi_1(\bb_L)$. Consequently,
$\cc_{\bb_L^1}(\sum_{j,k=0}^{r-1}
(e_{\gamma_{\widehat{r}-2j-1,\widehat{r}-2j}}+e_{\gamma_{1,\widehat{r}-2k}}))
=\{0\}$.

Next we consider the summand
$e_{\gamma_{1,\widehat{r}+1}}+e_{\beta_{1,\widehat{r}}}$. First
observe that $[\nn,e_{\gamma_{1,\widehat{r}+1}}]=\{0\}$, so
$\cc_\nn(e_{\gamma_{1,\widehat{r}+1}})=\nn$. Secondly, the root
$\beta_{1,\widehat{r}}$ lies in $\Psi_2(\gg_{\geqslant 2})$.
Thanks to property \eqref{eq 8.10}, we need only consider roots $\delta \in
\Psi_2(\bb_L)$. We see that the only roots $\delta\in
\Psi_2(\bb_L)$ with $\delta+\beta_{1,\widehat{r}}\in
\Psi(\gg_{\geqslant 2})$ are of the form $\beta_{\widehat{r}+1,j}$
or $\gamma_{\widehat{r}+1,k}$ where $\widehat{r}+1\leqslant
j\leqslant n$ and $\widehat{r}+1 < k\leqslant n$. Again the
Chevalley commutator relations imply
$[\gg_{\beta_{\widehat{r}+1,j}},e_{\beta_{1,\widehat{r}}}]=
\gg_{\beta_{1,j}}$ and
$[\gg_{\gamma_{\widehat{r}+1,k}},e_{\beta_{1,\widehat{r}}}]=
\gg_{\gamma_{1,k}}$. We also observe that $\beta_{j,k}$ and
$\gamma_{l,m}$ for $\widehat{r}+1<j,l$ have the property that
$\beta_{1,\widehat{r}+1}+\gamma_{l,m},\beta_{1,\widehat{r}+1}+\beta_{j,k}\notin
\Psi_2(\gg_{\geqslant 2})$. All the roots above exhaust
$\Psi_2(\bb_L)$, so we conclude that all the roots $\beta_{j,k}$
and $\gamma_{l,m}$ for $\widehat{r}+1<j,l$ of $\Psi_2(\bb_L)$ are
all contained in $\Psi(\cc_\nn(e_{\beta_{1,\widehat{r}}}))$. If
$s>0$, these roots form the set of positive roots of a root system
of type $B_{s-1}$, there are exactly $(s-1)^2$ positive roots in a
root system of type $B_{s-1}$ and so
$|\Psi(\cc_\nn(e_{\beta_{1,\widehat{r}}}))|=(s-1)^2$. Therefore,
$\dim \cc_\nn(e'_{r,s})=(s-1)^2$, clearly, if $s=0$ then, $\dim
\cc_\nn(e'_{r,s})=0$. \end{proof}

\begin{prop}
\label{prop 8.20} 
The $B_L$-orbit of $e'_{r,s}$ is dense in $\gg_{\geqslant 2}$.
\end{prop}

\begin{proof} 
Thanks to Lemma \ref{lem 8.10}, it is sufficient to show that $\dim B_L=\dim
C_{B_L}(e'_{r,s})+\dim \gg_{\geqslant 2}$. Lemma \ref{lem 8.60}
implies that $\dim \gg_{\geqslant 2}=2r^2+2s+r+1$ and Lemma
\ref{lem 8.70} implies that $\dim B_L = 2r^2+s^2+s+r+1$. 
By Lemma \ref{lem 8.80}, $e'_{r,s}$
is the minimal representative of its $U$-orbit in $\uu$.
Thus, by Lemma \ref{lem 8.40}, 
we have $\dim C_{B_L}(e'_{r,s})=\dim C_{T}(e'_{r,s})+\dim C_{U}(e'_{r,s})$. 
Consequently, Lemmas \ref{lem 8.30},
\ref{lem 8.80}, and \ref{lem 8.111} imply that, for
$s>0$, $\dim C_{B_L}(e'_{r,s})\leqslant n-2r-2+(s-1)^2=s^2-s$. So
\[
\dim C_{B_L}(e'_{r,s})+\dim \gg_{\geqslant 2}\leqslant s^2-s+2r^2+r+2s+1
%=2r^2+s^2+r+s+1
= \dim B_L.
\]
This clearly implies
$\dim B_L=\dim C_{B_L}(e'_{r,s})+\dim \gg_{\geqslant 2}$.
Similarly, if $s=0$, we get $\dim B_L=\dim C_{B_L}(e'_{r,s})+\dim
\gg_{\geqslant 2}$. 
\end{proof}

\begin{cor}
\label{cor 8.10} 
$\dim C_{B_L}(e'_{r,s})=s(s-1)$.
\end{cor}

Finally, from Lemma \ref{lem 8.10} we obtain

\begin{cor}
\label{cor 8.20}
If $G$ is of type $B_n$ and $e\in \N$ with $\htt(e)=3$, 
then $e$ is spherical.
\end{cor}

\subsection{Height Three Nilpotent Elements of $\mathfrak{so}_{2n}(k)$}
\label{subsect:soeven}
Assume for this subsection that $G$ is of type $D_n$ for $n\geqslant
4$, so $\gg=\mathfrak{so}_{2n}$. We know that the nilpotent orbits in $\gg$ are
classified by the partitions of $2n$ with even parts occurring
with even parity, see \cite[Thm.\  1.6]{Ja1}. We showed 
that the height three nilpotent
orbits correspond to partitions of $2n$ of the form
$\pi_{r,s}=[1^{2s+1},2^{2r},3]$ where $r\geqslant 1, s\geqslant 0$
and $2r+s+2=n$, see Proposition \ref{prop 4.50}. 
Similarly to the $B_n$ case, we denote the
corresponding orbit by $\mathcal{O}_{r,s}$ and a representative of
such an orbit by $e_{r,s}$. Because the proofs of the results in this
subsection are virtually identical to the ones in 
Subsection \ref{subsect:soodd}, they are omitted.

\begin{lem}
\label{lem 8.90} 
There are precisely $\left [\frac{n-2}{2} \right ]$ 
distinct height 3 nilpotent orbits in
$\gg$.
\end{lem}

\begin{comment}\begin{proof} It is clearly sufficient to show that there are
precisely $\left [ \frac{n-2}{2} \right ]$ partitions of $2n$ of
the form $\pi_{r,s}$. This is equivalent to finding partitions of
$n-2$ of the form $[1^s,2^r]$. Thus we find that $r$ must satisfy
$1\leqslant r\leqslant \frac{n-2}{2}$. Since $r$ is an integer,
the result follows. \end{proof}\end{comment}

Using \cite[\S 13]{Ca}, we can easily calculate that 
for $s>0$, $e_{r,s}$ has the weighted Dynkin diagram $\Delta(e_{r,s})$
as shown in Figure \ref{Dynkin3D} below.

\begin{figure}[ht]
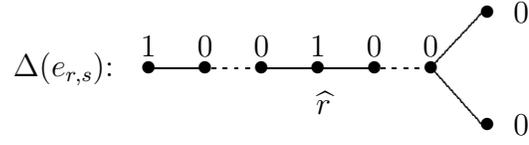

\beginpicture

\setcoordinatesystem units <1.5cm,1.5cm> point at 0 0     % sets scale
\setplotarea x from -5.1 to 3, y from -1 to 1     % sets plot size up

%D n%
\put {$\Delta(e_{r,s})$:} [l] at -2.2 0 \multiput {$\bullet$} at
-1 0 *4 .5 0 / \put {$\bullet$} at 1.5 0

 \putrule from -.96 0 to -.54 0 \putrule from 0.04 0 to .96 0
 \put {$\bullet$} at 2 0.5 \put {$\bullet$} at 2 -0.5

\plot 1.53 0.04 1.96 0.5 /  \plot 1.53 -0.02 1.96 -0.5 /

\put {$0$} at 0 .2 \put {$1$} at .5 .2 \put {$0$} at  1.5 .2 \put
{$0$} at 2.3 0.5 \put {$0$} at 2.3 -0.5 \put {$0$} at 1 .2 \put
{$0$} at -.5 .2 \put {$1$} at -1 .2

\put {$\widehat{r}$} at .55 -.3

\setdashes <.5mm, 1mm> \putrule from -.54 0 to 1.46 0

\endpicture 
\caption{Labelling of $\Delta(e_{r,s})$ for $s>0$.}
\label{Dynkin3D}

\end{figure}

Similarly, when $s=0$, the labelling of $\Delta(e_{r,0})$
is shown in Figure \ref{Dynkin3D0} below.
\begin{figure}[ht]
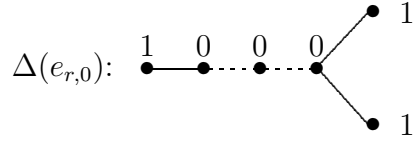


\beginpicture

\setcoordinatesystem units <1.5cm,1.5cm> point at 0 0     % sets scale
\setplotarea x from -4.7 to 3, y from -1 to 1     % sets plot size up

%D n%
\put {$\Delta(e_{r,0})$:} [l] at -1.2 0 \multiput {$\bullet$} at 0
0 *2 .5 0 / \put {$\bullet$} at 1.5 0

 \putrule from 0.04 0 to .46
0 \put {$\bullet$} at 2 0.5 \put {$\bullet$} at 2 -0.5 \plot 1.53
0.04 1.96 0.5 /  \plot 1.53 -0.02 1.96 -0.5 /

\put {$1$} at 0 0.2 \put {$0$} at .5 0.2 \put {$0$} at  1.5 0.2
\put {$1$} at 2.3 0.5 \put {$1$} at 2.3 -0.5 \put {$0$} at 1 0.2

\setdashes <.5mm, 1mm> \putrule from .54 0 to 1.46 0

\endpicture 
\caption{Labelling of $\Delta(e_{r,0})$.}
\label{Dynkin3D0}

\end{figure}

\begin{rem}
\label{rem 8.20} 
Note that there is
always an odd number of ``0'' labels between the first and second
``1'' labels in $\Delta(e_{r,s})$. 
If $s>0$, then there are $s+1$ ``0'' labels to the
right of the second ``1'' label. Finally, $s=0$ only if $n$ is
even.
\end{rem}

We refer to \cite[Planche IV]{Bo2} for information regarding the
root system of type $D_n$. We use the following notation for the
positive roots $\Psi^+$. Let $\alpha_1,\ldots,\alpha_n$ be the set
of simple roots of $\Psi^+$ and let 

\begin{align*}
\beta_{j,k}&=\alpha_j+\ldots +\alpha_k \text{ for } 1\leqslant j
\leqslant k\leqslant n\ \text{ not } j=n-1, k=n,\\
\beta_{j}&=\alpha_j+\ldots +\alpha_{n-2}+\alpha_n \text{ for }
1\leqslant j\leqslant n-2,\\
\gamma_{j,k}&=\alpha_j+ \ldots
+\alpha_{k-1}+2\alpha_k+\ldots +
2\alpha_{n-2}+\alpha_{n-1}+\alpha_n \text{ for } 1\leqslant j<
k<n-2.
\end{align*}
Here we again use the convention $\beta_{j,j} = \alpha_j$. Note that all
the possible $\beta$'s and $\gamma$'s exhaust $\Psi^+$.
\begin{comment}
\begin{exmp}\label{ex 8.30} The positive roots of $D_4$ are:
\[\Psi^+=\{\alpha_1,\alpha_2,\alpha_3,\alpha_4,\beta_{1,2},\beta_{2,3},\beta_{1,3},\beta_{1},\beta_{2,4},\beta_{1,4},\beta_2,\gamma_{1,2}\}.\]Note
that $\gamma_{1,2}$ is the highest root of $D_4$.  In fact,
$\gamma_{1,2}$ is always the highest root in a root system of type
$D_n$.\end{exmp}\end{comment}

Next we consider the structure of the abelian Lie subalgebra
$\gg_{\geqslant 2}=\gg(2)\bigoplus\gg(3)$.

\begin{lem}
\label{lem 8.100} 
An associated cocharacter for $e_{r,s}$ 
affords the following.
\begin{itemize}
\item[(i)] 
  $\Psi(\gg(2))=\left\{
\begin{array}{ll}
  \{ \beta_{1,j},\beta_{1},\gamma_{l,k},\gamma_{1,m} \mid 1< l
< k\leqslant \widehat{r}\leqslant j, \widehat{r}<m\} & \text{ if } s>0; \\
   \{ \beta_{1,n-1},\beta_{1},\beta_{i,n},\gamma_{j,k}\mid 2\leqslant
   i< \widehat{r}, 1<j<k< \widehat{r}\}& \text{ if } s=0. \\
    \end{array}\right.$\\
In particular, $\dim \gg(2)=2r^2-r+2s+2$.
\item[(ii)]  
$\Psi(\gg(3))=\left\{
\begin{array}{ll}
  \{\gamma_{1,k}\mid  k\leqslant \widehat{r}\}& \text{ if } s>0; \\
   \{\beta_{1,n},\gamma_{1,k}\mid 2\leqslant
   k< \widehat{r}\}& \text{ if } s=0. \\
    \end{array}\right.$\\
In particular, $\dim \gg(3)=2r$.
\end{itemize}
\end{lem}

\begin{comment}
\begin{proof}

Let $\lambda$ be an associated cocharacter for $e_{r,s}$ in $\gg$
and $T$ be a maximal torus of $G$ such that $\lambda(k^*)\leqslant
T$. For every $\delta \in \Psi$ we have that
$\gg_{\delta}\subseteq \gg(i)$ for some $i \in \{0,\pm 1,\pm 2,\pm
3\}$. The weighted Dynkin diagram of $e_{r,s}$ tells us where the
simple root spaces lie, see \eqref{eq 6.2}. Let
$\delta=\sum_{\alpha\in \Pi}c_{\delta,\alpha}\alpha$ be a
positive root.\begin{itemize}

\item A root space $\gg_{\delta}$ lies in $\gg(2)$ if and only if
$c_{\delta,\alpha_1}+c_{\delta,\alpha_{\widehat{r}}}=2$. We can
easily see that all of the roots listed above satisfy this
condition, and no others do. A counting exercise gives the
dimension.

\item A root space $\gg_{\delta}$ lies in $\gg(3)$ if and only if
$c_{\delta,\alpha_1}+c_{\delta,\alpha_{\widehat{r}}}=3$. Again we
can easily see that all of the roots listed above satisfy this
condition, and no others do. A simple counting exercise gives the
dimension.
\end{itemize} This completes the proof of the lemma.
\end{proof}\end{comment}

Next we look at the structure of the Lie subalgebra $\bb_L$ of $\gg(0)$.

\begin{lem}
\label{lem 8.110} 
An associated cocharacter for $e_{r,s}$ 
affords the following.

$\Psi(\bb_L)=\left\{
\begin{array}{ll}
\{\beta_{i},\beta_{j,k},\gamma_{l,m} \mid \widehat{r} < j \text{
or } 1<j\leqslant k < \widehat{r}, \widehat{r}<i\:,\: \widehat{r}< l<m\}  
& \text{ if } s>0; \\
 \{\beta_{j,k}\mid 1<j\leqslant k < \widehat{r}\} & \text{ if } s=0. \\
    \end{array}
\right.$ \\
In particular, $\dim \bb_L=2r^2+s^2+r+2s+2$. 
\end{lem}

\begin{comment}\begin{proof}Let $\lambda$ be an associated cocharacter for
$e_{r,s}$ in $\gg$ and $T$ be a maximal torus of $G$ such that
$\lambda(k^*)\leqslant T$. For every $\delta \in \Psi$ we have
that $\gg_{\delta}\subseteq \gg(i)$ for some $i \in \{0,\pm 1,\pm
2,\pm 3\}$. The weighted Dynkin diagram of $e_{r,s}$ tells us
where the simple root spaces lie, see \eqref{eq 6.2}. Let
$\delta=\sum_{\alpha\in \Pi}c_{\delta,\alpha}\alpha$ be a
positive root. A root space $\gg_{\delta}$ lies in $\bb_L$ if and
only if $c_{\delta,\alpha_1}+c_{\delta,\alpha_{\widehat{r}}}=0$.
We can easily see that all of the roots listed above satisfy this
condition, and no others do. A counting exercise gives $\dim \nn =
2r^2+s^2-r+s$. Since $\dim \mathfrak{t}=n$, we get that $\dim
\bb_L=2r^2+s^2+r+2s+2$. Alternatively, we note that the Lie
algebra $\gg(0)$ is isomorphic to
$\mathfrak{sl}_{\widehat{r}-1}(k)\times
\mathfrak{so}_{s+1}(k)\times k^2$, thanks to the weighted Dynkin
diagram. The structure and dimension of $\bb_L$ follows.
 \end{proof}\end{comment}

Similarly to the $B_n$ case, the roots of $\bb_L$ naturally form
two distinct subsets, namely the roots whose support 
lies strictly to the left of the second ``1'' label of the weighted
Dynkin diagram and those whose support lies strictly
to the right of the second ``1'' label on the weighted Dynkin
diagram. More precisely, we have
$\Psi(\bb_L)=\Psi_1(\bb_L)\cup\Psi_2(\bb_L)$ where
\begin{align*}
\Psi_1(\bb_L)&=\{ \beta_{j,k} \mid 1<j\leqslant k < \widehat{r}\},
\\
\Psi_2(\bb_L)&=\{\beta_{j,k},\beta_i,\gamma_{l,m} \mid
\widehat{r}<j\leqslant k , \widehat{r}< i, \widehat{r}< l <
m\}.
\end{align*}

Again we partition the roots of $\gg_{\geqslant 2}$ into two distinct
subsets. More precisely, we write 
$\Psi(\gg_{\geqslant 2}) = 
\Psi_1(\gg_{\geqslant 2})\cup\Psi_2(\gg_{\geqslant 2})$, where
for $s\geqslant 1$, we define

\begin{align*}
\Psi_1(\gg_{\geqslant
2})&=\{\gamma_{j,k} \mid 1\leqslant j < k\leqslant \widehat{r} \},\\
\Psi_2(\gg_{\geqslant 2})&=\{\beta_1,\beta_{1,j},\gamma_{1,k} \mid
\widehat{r} \leqslant j, \widehat{r} < k\},
\end{align*}
and for $s=0$, we define
\begin{align*}
\Psi_1(\gg_{\geqslant 2})&=\{\gamma_{j,k}
\mid 1\leqslant j < k\leqslant \widehat{r} \},\\
\Psi_2(\gg_{\geqslant 2})&=\{\beta_1,\beta_{1,n-1},\beta_{j,n}
\gamma_{1,k} \mid j\leqslant\widehat{r} < k\}.
\end{align*} 
Again, we have the following property of these sets:
\begin{equation}
\label{eq 8.20} 
\delta \in \Psi_i(\bb_L), \eta \in
\Psi_{3-i}(\gg_{\geqslant 2}) \ \Rightarrow\  \delta + \eta \notin \Psi,
\ i \in \{ 1,2\}.
\end{equation} 

For $s>1$, set
\[
e'_{r,s} :=
\sum_{j,k=0}^{r-1}(e_{\gamma_{\widehat{r}-2j-1,\widehat{r}-2j}}
+e_{\gamma_{1,\widehat{r}-2k}})+e_{\gamma_{1,\widehat{r}+1}}
+e_{\beta_{1,\widehat{r}}} \in \gg_{\geqslant 2},
\] 
for $s=1$, set
\[
e'_{r,1} := \sum_{j,k=0}^{r-1}(e_{\gamma_{\widehat{r}-2j-1,\widehat{r}-2j}}
+e_{\gamma_{1,\widehat{r}-2k}})+e_{\beta_{1,n}}+e_{\beta_{1,\widehat{r}}}
\in \gg_{\geqslant 2},
\] 
and for $s=0$, set
\[
e'_{r,0} := \sum_{j,k=1}^{r-1}(e_{\gamma_{\widehat{r}-2j-1,\widehat{r}-2j}}
+e_{\gamma_{1,\widehat{r}-2k}})+e_{\beta_{1,n}}+e_{\beta_{1,n-1}}
+e_{\beta_{n-2,n}}+e_{\beta_{1}} \in \gg_{\geqslant 2}.
\]

\begin{lem}
\label{lem 8.120}
With the notation as above, we have  
$|\supp(e'_{r,s})|=2r+2$, $\supp(e'_{r,s})$ is linearly
independent, and $\dim \cc_\nn(e'_{r,s})=s(s-1)$.
\end{lem}

\begin{prop}
\label{prop 8.33} 
The $B_L$-orbit of $e'_{r,s}$ is dense in $\gg_{\geqslant 2}$.
\end{prop}

\begin{comment}\begin{proof} By Lemma \ref{lem 8.10} it is sufficient to
show that $\dim B_L=\dim C_{B_L}(e'_{r,s})+\dim \gg_{\geqslant 2}$.
Lemma \ref{lem 8.100} implies that $\dim \gg_{\geqslant
2}=2r^2+r+2s+2$ and Lemma \ref{lem 8.110} implies that $\dim B_L =
2r^2+s^2+r+2s+2$. 
One readily checks that $e'_{r,s}$
is the minimal representative of its $U$-orbit in $\uu$ in the sense of
\cite{Go4}. Thus, by Lemma \ref{lem 8.40} we have $\dim
C_{B_L}(e'_{r,s})=\dim C_{T}(e'_{r,s})+\dim C_{U}(e'_{r,s})$. Lemmas
\ref{lem 8.30} and \ref{lem 8.120} and Proposition \ref{prop 8.30}
imply that $\dim C_{B_L}(e'_{r,s})\leqslant n-2r-2+s(s-1)=s^2$. So
$\dim C_{B_L}(e'_{r,s})+\dim \gg_{\geqslant 2}\leqslant
s^2+2r^2+r+2s+2=\dim B_L$, this clearly implies $\dim B_L=\dim
C_{B_L}(e'_{r,s})+\dim \gg_{\geqslant 2}$. Similarly, if $s=0$, we
get $\dim B_L=\dim C_{B_L}(e'_{r,s})+\dim \gg_{\geqslant 2}$.
\end{proof}\end{comment}

%\begin{cor}\label{cor 8.30} $\dim
%C_{B_L}(e'_{r,s})=s^2$.\end{cor}

\begin{cor}
\label{cor 8.40}
If $G$ is of type $D_n$ and $e\in \N$ with $\htt(e)=3$, then $e$ is
spherical.
\end{cor}

\subsection{Height Three Nilpotent Elements of the Exceptional Lie Algebras}
\label{sub:ex}
We fix an ordering of the roots $\alpha_1,\ldots,\alpha_r$ of
$\Psi(\gg_{\geqslant 2})$ such that $\alpha_i\prec \alpha_j$ for
$i<j$. Define the subalgebra $\mm_i$ of $\gg_{\geqslant 2}$ by
setting $\mm_i=\bigoplus_{j=i+1}^r\gg_{\alpha_j}$ 
and the quotient $\qq_i$ by $\qq_i=\gg_{\geqslant 2}/\mm_i$
for $0 \leqslant i \leqslant r$.
Let $B$ be a Borel subgroup of $G$ such that $\gg_{\geqslant
2}\subseteq \Lie R_u(B)=\uu$. 
Note that each $\qq_i$ is a $B$-module. 

The computer programme, $\mathsf{DOOBS}$, devised by S.M.\ Goodwin  
allows us to
determine whether $B$ acts on $\gg_{\geqslant 2}$ with a dense
orbit. For details of the $\mathsf{GAP4}$ (\cite{Ga}) 
computer algebra program of Goodwin, we refer the reader to 
\cite{Go2} and \cite{Go2.5}. 
%The $\mathsf{DOOBS}$ algorithm runs in the 
%For the convenience of the
%reader, we sketch how 
Working inductively, starting with $i = 0$, at each stage of the algorithm, 
$\mathsf{DOOBS}$ determines a representative $x_i+\mm_i$,
with $\supp(x_i)$ linearly independent of a dense $B$-orbit on
$\qq_i$ or decides that $B$ does not act on $\qq_i$ with a dense
orbit.

$\mathsf{DOOBS}$ also keeps a record of the primes for which
$\dim_p \cc_{\uu}(x_i+\mm_{i+1})>\dim_0 \cc_{\uu}(x_i+\mm_{i+1})$,
where $\dim_p\cc_{\uu}(x_i+\mm_{i+1})$ and
$\dim_0\cc_{\uu}(x_i+\mm_{i+1})$ denote the dimension of
$\cc_{\uu}(x_i+\mm_{i+1})$ over a field of characteristic $p$ and
characteristic $0$ respectively, see Remark 3.2 in \cite{Go2}. For
these primes we cannot conclude that $B$ acts on $\gg_{\geqslant
2}$ with a dense orbit. 
If $\mathsf{DOOBS}$ determines that $B$ acts on $\gg_{\geqslant 2}$ 
with a dense orbit, then it calculates a representative of
the dense orbit and a list of primes for which the result is
not necessarily valid.

There is a variant of $\mathsf{DOOBS}$ called
$\mathsf{DOOBSLevi}$, see \cite{Go2.5}. This program considers a
parabolic subgroup $P=LR_u(P)$ and determines whether a Borel subgroup $B_L$
of $L$ acts on an ideal of $\Lie R_u(P)$ with a dense orbit. The
algorithm used to determine whether $B_L$ acts on an
ideal with a dense orbit is essentially the same as the
$\mathsf{DOOBS}$ algorithm, with $B_L$ replacing $B$.
$\mathsf{DOOBSLevi}$ also records the primes for which its
conclusions are not necessarily valid.

\begin{comment}The input required for $\mathsf{DOOBSLevi}$ is the type of
$G$, the rank of $G$, a generating set for the ideal and a
generating set for the Levi subgroup. For this we number the roots
of $\Psi^+$ as $\alpha_1,\ldots,\alpha_r$ and note that if $\rk
G=n$, then $\alpha_1, \ldots,\alpha_n$ are the simple roots of
$\Psi^+$. For $\mathsf{DOOBSLevi}$ to generate the ideal we input
a subset $\{i_1,\ldots,i_k\}$ of $\{1,\ldots,r\}$ such that the
roots corresponding to $\{i_1,\ldots,i_k\}$ generate the ideal.
For $\mathsf{DOOBSLevi}$ to generate a Levi subgroup we input a
subset $\{i_1,\ldots,i_s\}$ of $\{1,\ldots,n\}$ such that $L$ is
generated by the simple roots corresponding to
$\{1,\ldots,n\}\setminus\{i_1,\ldots,i_s\}$.\end{comment}

Let $e \in \N$ of height $3$ and let 
$\lambda$ be a cocharacter of $G$ that is associated to $e$.
We use the same numbering of the positive roots as in $\mathsf{GAP4}$.
Table \ref{tab 8.10} 
below lists the roots whose root subgroups together with $T$ generate the Levi
subgroup $C_G(\lambda)$ and we also list the roots whose root subspaces 
generate $\gg_{\geqslant 2}$ (as a $B$-submodule of $\gg$) 
for the $7$ cases of 
height three nilpotent orbits for the simple exceptional groups,
see Proposition \ref{prop 8.10}. These are determined by means of the
weighted Dynkin diagrams.

\begin{comment} Indeed, see Section \ref{levi} for how to determine the Levi
subgroup from the weighted Dynkin diagram. A similar method can be
used for $\gg_{\geqslant 2}$, or we can note that generators for
$\gg_{\geqslant 2}$ can be read off Table 1 in
\cite{Ro}.\end{comment}

\begin{table}[ht]
\renewcommand{\arraystretch}{1.5}
  \centering
\begin{tabular}{|c|c|c|c|}
  \hline
  % after \\: \hline or \cline{col1-col2} \cline{col3-col4} ...
  Type of $G$& Bala--Carter Label & Generators for $L$ & Generators for $\gg_{\geqslant 2}$ \\
  \hline
  $G_2$ & $\tilde{A_1}$&$\alpha_2$ & $\alpha_4$ \\
  $F_4$ & $A_1 + \tilde{A_1}$&$\Pi\setminus\{\alpha_4\}$ & $\alpha_{16}$ \\
  $E_6$ & $3A_1$&$\Pi\setminus\{\alpha_4\}$ & $\alpha_{24}$ \\
  $E_7$ & $(3A_1)'$&$\Pi\setminus\{\alpha_3\}$ & $\alpha_{37}$ \\
  $E_7$ & $4A_1$&$\Pi\setminus\{\alpha_{2},\alpha_{7}\}$ & $\alpha_{30},\alpha_{53}$
  \\
  $E_8$ & $3A_1$ &$\Pi\setminus\{\alpha_7\}$ & $\alpha_{74}$ \\
  $E_8$ & $4A_1$&$\Pi\setminus\{\alpha_2\}$ & $\alpha_{69}$ \\
  \hline
\end{tabular}
\bigskip
\caption{Height Three Nilpotent Orbits in the Exceptional Lie Algebras.}
\label{tab 8.10}
\end{table}

The height $3$ cases for the exceptional groups were analyzed
using the $\mathsf{DOOBSLevi}$ algorithm. 
It turns out that there are no characteristic restrictions in these cases:
%For further details on these calculations, see \cite{Fo}.

\begin{lem}
\label{lem 8.130} 
If $G$ is simple of exceptional type and $e\in \N$ 
with $\htt(e)=3$, then $e$ is spherical.
\end{lem}

Corollaries \ref{cor 8.20} and \ref{cor 8.40} combined with 
Lemma \ref{lem 8.130} give the following result.

\begin{prop}
\label{prop 8.40}
Let $G$ be a connected reductive algebraic group and let
$e \in \N$. If $\htt(e)=3$, then $e$ is spherical.
\end{prop}

\begin{proof} 
If $G$ is simple, then the statement follows from
Corollaries \ref{cor 8.20} and \ref{cor 8.40} and
Lemma \ref{lem 8.130}.
In the general case we argue as in Lemma \ref{lem 6.50} to reduce to the 
simple case.
%let $\D G=G_1G_2\cdots G_r$ be a
%commuting product of simple groups. We have 
%that $e = e_1+e_2+\ldots +e_r$ where $e_i \in \gg_i = \Lie G_i$ and each $e_i$
%is nilpotent. We have  
%$\htt(e) = \max_{1\leqslant i \leqslant r}\htt(e_i)$.
%%, see the proof of Lemma \ref{lem 6.50}. 
%In particular, $\htt(e_i)\leqslant \htt(e)=3$ for all $i$.
%We have $\kappa_G(G\cdot e)=\sum_{i=1}^r
%\kappa_{G_i}(G_i\cdot e_i)$. By the simple case just discussed and 
%by Lemma \ref{lem 6.50}, the desired result follows.
\end{proof}

\subsection{The Classification}
\label{sub:class}
Our main classification theorem now follows readily 
from Lemma \ref{lem 6.50} and Propositions \ref{prop 7.20} 
and \ref{prop 8.40}.

\begin{thm}
\label{thm 8.40}
Let $G$ be a connected reductive
algebraic group. Suppose that $\cha k$ is a good prime for $G$. 
Then a nilpotent
element $e \in \gg$ is spherical if and only if $\htt(e)\leqslant 3$.
\end{thm}

\begin{rem}
\label{rem 8.50}
Let $G$ be a simple 
algebraic group and let $\cha k$ be a good prime for $G$. 
Then the spherical nilpotent orbits are given in 
Table \ref{t:sherical}. 
We present the orbits 
by listing the corresponding partition 
in the classical cases or by giving the corresponding 
Bala--Carter label for the exceptional groups.
\end{rem}

\begin{table}[ht]
\renewcommand{\arraystretch}{1.5}
  \centering
\begin{tabular}{|c|l|}
  \hline
  Type of $G$ & Spherical Orbits \\
  \hline
$A_n$ &  $[1^j,2^i]$\\
$B_n$ & $[1^j,2^{2i}]$, or $[1^j,2^{2i},3]$ with $i \geqslant 0$\\
$C_n$ & $[1^{2j},2^i]$ \\
$D_n$ & $[1^j,2^{2i}]$, or $[1^j,2^{2i},3]$ with $i \geqslant 0$\\
\hline
$G_2$ & $A_1$ or $\tilde{A_1}$\\
$F_4$ & $A_1$, $\tilde{A_1}$, or $A_1 + \tilde{A_1}$\\
$E_6$ & $A_1$, $2A_1$, or $3A_1$\\
$E_7$ & $A_1$, $2A_1$, $(3A_1)'$, $(3A_1)''$, or $4A_1$\\
$E_8$ & $A_1$, $2A_1$, $3A_1$, or $4A_1$\\
  \hline
\end{tabular}
\bigskip
\caption{The spherical nilpotent Orbits for $G$ simple.}
\label{t:sherical}
\end{table}

\begin{rem}
\label{rem 8.55}
Using the fact that in good characteristic 
%we have a Springer map 
%between the unipotent variety of $G$ and $\N \subset \gg$, 
%and the fact that such 
a Springer map affords a bijection between the set of unipotent
$G$-conjugacy classes and the set of nilpotent $G$-orbits (see \cite{serre}),
Theorem \ref{thm 8.40} also gives a classification of the spherical unipotent
classes in $G$. Here we define the height of a unipotent element $u$ of $G$
as the height of the image of $u$ in $\N$
under a Springer isomorphism. 
\end{rem}

\begin{comment}\begin{proof} Use Theorems \ref{thm 6.30} and \ref{thm 8.10}
\end{proof}

We have now classified the spherical nilpotent orbits, over a
field of good characteristic for $G$. The classification is the
same as the classification over a field of characteristic zero,
see \cite{Pa3}. In \emph{loc.\!\! cit.\!} the proof that height
three implies spherical was also shown by case-by-case
considerations, however a general proof was later found, see
\cite{Pa1}. In the final chapter we discuss some applications of
this result.

\chapter{}\label{C9}

In this chapter we discuss some applications of the main result,
Theorem \ref{thm 8.40}. Throughout this chapter $G$ is a connected
reductive algebraic group over an algebraically closed field $k$,
the Lie algebra of $G$ is denoted by $\gg$ and the characteristic
of $k$ is good for $G$, except in the final section.

\end{comment}

\section{Applications and Complements}
\label{sec:appl}
Here we discuss some applications of the main result and some further 
consequences.

\subsection{Spherical Distinguished Nilpotent Elements}
\label{sub:dist}
Recall that a nilpotent element $e\in \N$ 
is distinguished in $\gg$ if every torus contained in
$C_G(e)$ is contained in the centre of $G$. For now we assume that
$G$ is simple, so $e$ is distinguished in $\gg$ if and only if any torus
contained in $C_G(e)$ is trivial and hence $C_G(e)^{\circ}$
is unipotent. Further recall that 
$\kappa_G(G\cdot e)=\kappa_G(G/C_G(e)^{\circ})$, 
cf.\ equation \eqref{eq 1.10}. 
Since $C_G(e)^{\circ}$ is
connected and unipotent, it is contained in the unipotent radical $U$
of a Borel subgroup $B = TU$ of $G$. 
Let $B^-=TU^-$ be the unique opposite Borel subgroup to $B=TU$ relative 
to $T$, see \cite[\S 26.2]{Hu}. 
Consequently, 
$B^-\cap C_G(e)^{\circ} \subseteq B^-\cap U = \{1\}$. 
Thus, by equation \eqref{SubCom}, we have 
$\kappa_G(G/C_G(e)^{\circ}) = \dim G-\dim C_G(e)^{\circ} -\dim B^- = 
\dim U -\dim C_G(e)$,
or equivalently, $\kappa_G(G\cdot e) = |\Psi^+| - \dim C_G(e)$.
We summarize what we have just shown.

\begin{prop}
\label{prop 20.11}
Let $e \in \N$ be a distinguished
nilpotent element. Then 
\[
\kappa_G(G\cdot e)=|\Psi^+|-\dim C_G(e).
\]
\end{prop}

\begin{rem}
\label{rem 20.10}
Proposition \ref{prop 20.11} was first observed by Panyushev for a field of
characteristic zero in \cite[Cor.\ 2.4]{Pa3}.
\end{rem} 

If $G$ is a
simple classical group, then the distinguished nilpotent elements
are given as follows, see Lemmas 4.1 and 4.2 in \cite{Ja1}.

\begin{lem}
\label{lem 20.01} 
Let $e \in \N$ and let $\pi_e$ be the corresponding partition of $\dim V$. 
\begin{itemize}
\item[(i)] If $G=\SL(V)$, 
then $e$ is distinguished if and only if $\pi_e=[\dim V]$.
\item[(ii)] If $G=\SP(V)$, then $e$ is
distinguished if and only if $\pi_e$ consists only of distinct
even parts.
\item[(iii)]  If $G=\SO(V)$, then $e$ is
distinguished if and only if $\pi_e$ consists only of distinct odd
parts.
\end{itemize}
\end{lem}

%\begin{proof}In \cite[Lem.\ 4.2]{Ja1} the above is proved
%for $G=\SP(V)$ and $\SO(V)$. In \cite[Lem.\ 4.1]{Ja1} it is proved
%that if $G=\GL(V)$, then $e$ is distinguished if and only if $e$
%has partition $\pi=[\dim V]$. The result follows readily for
%$\SL(V)$. \end{proof}

\begin{cor}
\label{cor 20.1} 
If $G=\SO(V)$ and $e \in \N$ is spherical and distinguished, then
$\htt(e)=2$.
\end{cor}

\begin{proof} 
Thanks to Proposition \ref{prop 8.10}, 
the height $3$ nilpotent elements have partitions of the
form $\pi=[1^s,2^{2r},3]$, where $r>0$. Thus such a partition has
even parts and so is not distinguished. So if $e$ is spherical and
distinguished, then $\htt(e)=2$. 
\end{proof} 

Proposition \ref{prop 4.50} and Lemma \ref{lem 20.01} 
imply the following result.

\begin{prop}
\label{prop 20.02} 
Let $e \in \N$ be distinguished
and $\pi_e$ be the corresponding partition of $\dim V$.
\begin{itemize}
\item[(i)] If $G=\SL(V)$, then $\htt(e)=2$ if and only if $\pi_e=[2]$.
\item[(ii)] If $G=\SP(V)$, then $\htt(e)=2$  if and only if $\pi_e=[2]$.
\item[(iii)]  If $G=\SO(V)$, then $\htt(e)=2$ if and only if $\pi_e=[3]$ or
$\pi_e=[1,3]$.
\end{itemize}
\end{prop}  

%It follows from Proposition \ref{prop 20.02} that
%if $G$ is reductive and admits a height $2$ distinguished
%nilpotent orbit, then $\D G$ is of type $A_1^t$. 
%In particular, it follows from Proposition \ref{prop 20.02} and
%Theorem \ref{thm 8.40} that for $G$ simple
%the only occurrence of a spherical distinguished
%nilpotent orbit is when $G=\SL_2(k)$.

\begin{thm}
\label{thm 20.02} 
If $G$ is a simple algebraic group and
$e \in \N$ is spherical and distinguished, then $G$ is
of type $A_1$.
\end{thm}

\begin{proof}
For $G$ simple classical, 
Proposition \ref{prop 20.02} implies that $G$ is
of type $A_1$. For $G$ of exceptional type it follows from 
Remark \ref{rem 8.50} and the tables in
\cite[\S 13]{Ca} that there are no nilpotent orbits in $\gg$
that are both spherical and distinguished. 
\end{proof}

\subsection{Orthogonal Simple Roots and Spherical Nilpotent Orbits}
\label{sub:orth}
In \cite[Thm.\  3.4]{Pa2}, Panyushev proved that if the
characteristic of $k$ is zero, then $e\in \N$ is spherical if and
only if there exist pairwise orthogonal simple roots
$\alpha_1,\alpha_2,\ldots ,\alpha_t$ in $\Pi$
such that $G\cdot e$ contains an element of the form
$\sum_{i=1}^{t}e_{\alpha_i}$ where
$e_{\alpha_i}\in\gg_{\alpha_i}\setminus\{0\}$. By pairwise
orthogonal we mean that
$\langle\alpha_i,\alpha_j\rangle=0$ for $i\neq j$. 
In this subsection we show that
this is also the case if the characteristic of $k$ is good for
$G$.

%Also $\Psi$ is the root system of $G$ with respect to a
%maximal torus $T$ of $G$. 
%First we need some information regarding
%the nilpotent orbits of $G$ when $\D G$ is of type $A_1^{t}$.

\begin{lem}
\label{lem 20.15} 
Let $\D G$ be of type $A_1^{t}$ for some $t \geqslant 1$. Then
there is precisely one distinguished nilpotent orbit in $\N$.
\end{lem}

\begin{proof} 
Since the nilpotent orbits of $G$ in $\gg$
are precisely the nilpotent orbits of $\D G$ in $\Lie \D G$, we
may assume that $G$ is semisimple. Thus, $G=G_1 G_2 \cdots G_r$ and
each $G_i$ is of type $A_{1}$. %By Lemma \ref{lem 20.01}
There is precisely one distinguished nilpotent orbit when $G_i$ is
of type $A_1$: the unique non-zero nilpotent orbit. Also $G\cdot e$
is distinguished in $\gg$ if and only if 
$G_i\cdot e_i$ is distinguished in $\gg_i = \Lie G_i$ for all $i$, where
$e = e_1+\ldots +e_r$ and $e_i \in \gg_i$ is nilpotent. 
\end{proof}

%Let $e \in \N$, $\lambda \in \Omega_G^a(e)$, 
%and $S$ be a maximal torus of $C_G(e)$. 

\begin{lem}
\label{lem 20.111}
Let $e \in \N$ and $S$
be a maximal torus of $C_G(e)$. 
Then $\D C_G(S)$ is of type $A_{1}^{t}$  for some $t \geqslant 1$
if and only if there exist
pairwise orthogonal simple roots $\alpha_1$, $\alpha_2$, $\ldots$,
$\alpha_t$ in $\Pi$ such that $G\cdot e$ contains an element
of the form $\sum_{i=1}^{t}e_{\alpha_i}$, where
$e_{\alpha_i}\in\gg_{\alpha_i}\setminus \{0\}$.
\end{lem}

\begin{proof}
Suppose that $\D C_G(S)$ is of type $A_{1}^{t}$. Let
$\alpha_1,\ldots,\alpha_t$ be simple roots of $\Phi$, where $\Phi$
is the root system of $C_G(S)$ relative to a maximal torus $T$ of
$C_G(S)$. As $\D C_G(S)$ is of type $A_{1}^t$, the roots
$\alpha_1,\ldots,\alpha_t$ are pairwise orthogonal. Clearly, $e\in
\Lie C_G(S) = \cc_{\gg}(S)$ and $e$ is distinguished in $\cc_{\gg}(S)$,
see Proposition \ref{prop4.20}. By Lemma \ref{lem 20.15}, an
element of the form $\sum_{i=1}^{t}e_{\alpha_i}$ is also
distinguished in $\cc_\gg(S)$ and there is precisely one
distinguished nilpotent orbit in $\cc_\gg(S)$. Thus, 
$e$ and $\sum_{i=1}^{t}e_{\alpha_i}$ are in the same 
$C_G(S)$-orbit, hence they are in the same
$G$-orbit. So $G\cdot e$ contains an element of the desired form.
%$\sum_{i=1}^{t}e_{\alpha_i}$. 
%As $C_G(S)$ is a Levi
%subgroup of $G$, the simple roots $\alpha_1,\ldots,\alpha_r$ are
%also simple roots of $\Psi$.

Conversely, suppose that there exist pairwise orthogonal simple roots
$\alpha_1,\alpha_2,\ldots , \alpha_t\in \Psi$ such that $G\cdot e$
contains an element of the form
$e' = \sum_{i=1}^{t}e_{\alpha_i}$. Let $H$ be the subgroup of
$G$ generated by $\{T,U_{\pm \alpha_i}\mid 1\leqslant i\leqslant
t\}$, where $T$ is as in the previous paragraph. Then $\D H$ is of
type $A_1^t$. By construction, $e'$ is distinguished in
$\hh$. By Proposition \ref{prop4.20}, $H$ is of the
form $C_G(S')$, where $S'$ is a maximal torus of $C_G(e')$. Thus,
$\D C_G(S')$ is of type $A_{1}^t$. Since $e$ and $e'$ 
are $G$-conjugate, so are $C_G(e)$ and $C_G(e')$, as well as $S$ and $S'$.
Finally, we get that $C_G(S)$ and $C_G(S')$ are
$G$-conjugate. The result follows.
 \end{proof}

\begin{lem}
\label{lem 20.1} 
If $e \in \N$ is spherical, then $\D
C_G(S)$ is of type $A_{1}^{t}$ for some $t \geqslant 1$.
\end{lem}

\begin{proof} 
Let $\lambda$ be a cocharacter of $G_G(S)$ that is associated  to $e$,
i.e.\ $\lambda \in \Omega_{C_G(S)}^a(e)$. 
Then, since $\Lie C_G(S)=\cc_\gg(S)$, it follows from 
\cite[Cor.\ 3.21]{FoRo} that $\lambda \in \Omega_G^a(e)$. 
As $e$ is spherical in $\gg$, we have $\htt(e)\leqslant 3$, by 
Theorem \ref{thm 8.40}. As $\lambda \in \Omega_{C_G(S)}^a(e)$, 
we also have $\htt(e)\leqslant 3$ when we
regard $e$ as an element of $\cc_\gg(S)$. Thus, again by Theorem
\ref{thm 8.40}, $e$ is spherical in $\cc_\gg(S)$. So $e$ is
distinguished and spherical in $\cc_\gg(S)$ and so  
$\D C_{G}(S)$ is of type $A_{1}^t$, by Theorem \ref{thm 20.02}.
\end{proof} 

In order to prove the reverse implication of Lemma
\ref{lem 20.1} we first need to consider the group $C_G(S)$. If
$G$ is classical, then the structure of $C_G(S)$ can be determined
from the partition $\pi_e$ corresponding to $e$; see
\cite[\S 4.8]{Ja1} for the following result.

\begin{lem}
\label{lem 20.11}
Let $G$ be simple classical and $e\in \N$ 
with  corresponding partition $\pi_e$. 
%Then we have the following.
\begin{itemize}
\item[(i)]
If $G$ is of type $A_n$ and $\pi_e=[1^{r_1},2^{r_2},\ldots]$,
then $\D C_G(S)$ is of
type $\prod_{i\geqslant 1}A_{i-1}^{r_i}$.
\item[(ii)] 
If $G$ is of type $B_n$ and
$\pi_e=[1^{2s_1+\epsilon_1},2^{2s_2},3^{2s_3+\epsilon_3},\ldots]$,
where $s_i\geqslant  0$ and
$\epsilon_i \in \{0,1\}$, then $\D C_G(S)$ is of type
$\prod_{i\geqslant 1}A_{i-1}^{s_i} \times B_{m}$, where
$2m+1=\sum_{\epsilon_i \neq 0} i$.
\item[(iii)] 
If $G$ is of type $C_n$ and
$\pi_e=[1^{2s_1},2^{2s_2+\epsilon_2},3^{2s_3},4^{2s_4+\epsilon_4},\ldots]$,
where $s_i\geqslant  0$ and
$\epsilon_i \in \{0,1\}$, then $\D C_G(S)$ is of type
$\prod_{i\geqslant 1} A_{i-1}^{s_i} \times C_{m}$, where
$2m=\sum_{\epsilon_i \neq 0} i$.
\item[(iv)] 
If $G$ is of type $D_n$ and
$\pi_e=[1^{2s_1+\epsilon_1},2^{2s_2},3^{2s_3+\epsilon_3},\ldots]$,
where $s_i\geqslant 0$ and
$\epsilon_i \in \{0,1\}$, then $\D C_G(S)$ is of type
$\prod_{i\geqslant 1} A_{i-1}^{s_i} \times D_{m}$, where
$2m=\sum_{\epsilon_i \neq 0} i$.
\end{itemize}
\end{lem}

\begin{lem}
\label{lem 20.2} 
If $G$ is simple classical and $\D C_G(S)$ is of type
$A_{1}^{t}$, then $e$ is spherical. 
\end{lem} 

\begin{proof} 
%We prove this case-by-case. 
First suppose that $G$ is of type $A_n$. 
Since $\D C_G(S)$ is of type $A_{1}^{t}$, it follows from 
Lemma \ref{lem 20.11} that 
$r_i=0$ for all $i \geqslant 3$.
Thus $\pi_e=[1^{r_1},2^{r_2}]$ and
so $e$ is spherical, by Remark \ref{rem 8.50}.

Let $G$ be of type $B_n$. 
Since $\D C_G(S)$ is of type $A_{1}^{t}$, it follows from Lemma
\ref{lem 20.11} that $s_i=0$ for $i\geqslant 3$ and
$m\leqslant 1$, so $2m+1 \leqslant 3$. Since $2m+1$ is a sum of
distinct odd integers, we  either have $2m+1=1$ or $2m+1=3$. Thus
$\pi_e=[1^{2s_1+1},2^{2s_2}]$ or
$\pi_e=[1^{2s_1},2^{2s_2},3]$ and 
so $e$ is spherical, again by Remark \ref{rem 8.50}.

Let $G$ be of type $C_n$. Since $\D C_G(S)$ is of type $A_{1}^{t}$, 
it follows from Lemma
\ref{lem 20.11} that $s_i=0$ for $i\geqslant 3$ and
$m\leqslant 1$, so $2m\leqslant 2$. Since $2m$ is a sum of
distinct even integers, we either have $2m=0$ or $2m=2$. 
Thus $\pi_e=[1^{2s_1},2^{2s_2}]$ or $\pi_e=[1^{2s_1},2^{2s_2+1}]$ 
and so, by Remark \ref{rem 8.50}, $e$ is spherical.

Finally, let $G$ be of type $D_n$. 
Since $\D C_G(S)$ is of type $A_{1}^{t}$, it again follows from Lemma
\ref{lem 20.11} that $s_i=0$ for $i\geqslant 3$ and
$m\leqslant 2$, so $2m\leqslant 4$. Since $2m$ is a sum of
distinct odd integers, we either have $2m = 0$ or $2m = 1+3$. Thus
$\pi_e = [1^{2s_1},2^{2s_2}]$ or
$\pi_e = [1^{2s_1+1},2^{2s_2},3]$ and so, by Remark \ref{rem 8.50},
$e$ is spherical.
\end{proof}

All that remains is to check the exceptional cases. 
%Again this is done case-by-case. 
The Bala--Carter label of $e \in \N$ gives the Dynkin type of a Levi
subgroup $L$ of $G$ such that $e$ is distinguished in $\Lie \D L$.
By Proposition \ref{prop4.20}, such a Levi subgroup is
the centralizer of a maximal torus of $C_G(e)$. Thus, the
Bala--Carter label gives the type of $\D C_G(S)$. It follows from the
tables in \cite[\S 13]{Ca} and Remark \ref{rem 8.50} that 
any nilpotent orbit with Bala--Carter label $A_1^t$ is spherical. 
We summarize this in Table \ref{tabel20.1} below.

\begin{table}[ht]
\renewcommand{\arraystretch}{1.2}
  \centering
\begin{tabular}{|c|c|c||c|c|c|}
  \hline
  % after \\: \hline or \cline{col1-col2} \cline{col3-col4} ...
  Type & Bala--Carter Label & Height &Type & Bala--Carter Label & Height\\
  \hline
  $G_2$ & $A_1$ & 2 &$E_7$ & $A_1$ & 2 \\
  $G_2$ & $\tilde{A_1}$ & 3& $E_7$ & $2A_1$ & 2
  \\\cline{1-3}
  $F_4$ & $A_1$ & 2 &$E_7$ & $(3A_1)''$ & 2 \\
  $F_4$ & $\tilde{A_1}$ & 2 & $E_7$ & $(3A_1)'$ & 3 \\
  $F_4$ & $A_1 + \tilde{A_1}$ & 3&$E_7$ & $4A_1$ & 3
  \\\hline
  $E_6$ & $A_1$ & 2 &$E_8$ & $A_1$ & 2 \\
  $E_6$ & $2A_1$ & 2 &$E_8$ & $2A_1$ & 2 \\
  $E_6$ & $3A_1$ & 3 &$E_8$ & $3A_1$ & 3 \\
 - &- &- & $E_8$ & $4A_1$ & 3 \\
  \hline
\end{tabular}
\bigskip
\caption{Orbits in Exceptional Lie Algebras with $\D C_G(S)$ 
of Type $A_{1}^t$.}
\label{tabel20.1}
\end{table}

%This leads us to the following.
%We summarize:

\begin{lem}
\label{lem 20.25}
If $G$ is a simple exceptional algebraic group and $\D C_G(S)$ is of type
$A_{1}^{t}$, then $e$ is spherical.
\end{lem}

\begin{lem}
\label{lem 20.3}
Let $e\in \N$. If $\D
C_G(S)$ is of type $A_{1}^{t}$, then $e \in \gg$ is spherical.
\end{lem}

\begin{proof} 
For $G$ simple, 
the result follows from Lemmas \ref{lem 20.2} and \ref{lem 20.25}. 
In the general case let $\D G=G_1G_2\cdots G_r$ be a commuting product of
simple groups and $e = e_1+e_2+\ldots +e_r$, where $e_i \in \gg_i = \Lie G_i$
and each $e_i$ is nilpotent. %, see Proposition \ref{prop 3.90}.
A maximal torus $S$ of $C_G(e)$ is of the
form $S_1S_2\cdots S_r$, where $S_i$ is a maximal torus of
$C_{G_i}(e_i)$. The simple case implies that $\D C_{G_i}(S_i)$ is of
type $A_{1}^{t}$. 
\end{proof}

Lemmas \ref{lem 20.3} and \ref{lem 20.111} now imply the main 
result of this subsection.

\begin{thm}
\label{thm 20.2} 
Let $e \in \N$ and let $S$
be a maximal torus of $C_G(e)$. Then the following are
equivalent.
\begin{itemize}
\item[(i)] $e$ is spherical;
\item[(ii)] $\D C_G(S)$ is of type $A_{1}^t$;
\item[(iii)] there exist pairwise orthogonal simple roots 
$\alpha_1,\alpha_2,\ldots ,\alpha_t\in \Pi$ such
that $G\cdot e$ contains an element of the form
$\sum_{i=1}^{t}e_{\alpha_i}$, where
$e_{\alpha_i}\in\gg_{\alpha_i}\setminus\{0\}$.
\end{itemize}
\end{thm}

\begin{comment}
\begin{exmp}\label{ex 20.05}\begin{itemize}\item Let $G$ be of type $A_3$
and let $\alpha_1,\alpha_2$ and $\alpha_3$ be simple roots of
$\Psi$. Clearly, the following sets contain pairwise orthogonal
simple roots $\{\alpha_1\}\:,\:\{\alpha_2\}\:,\:\{\alpha_3\}$ and
$\{\alpha_1,\alpha_3\}$. If $e_{\alpha_i}\in
\gg_{\alpha_i}\setminus\{0\}$ for $i=1,2,3$, then $e_{\alpha_i}$
and $e_{\alpha_j}$ lie in the same $G$-orbit, since the roots
$\alpha_i$ and $\alpha_j$ are conjugate under the Weyl group of
$G$. So by Theorem \ref{thm 20.2} we have two distinct
non-zero spherical nilpotent orbits in $\gg$, cf.\ Theorem \ref{thm 8.50}.

\item Let $G$ be of type $G_2$ and let $\alpha$ and $\beta$ be
simple roots of $\Psi$, with $\alpha$ the long root. Trivially the
following sets contain pairwise orthogonal simple roots
$\{\alpha\}$ and $\{\beta\}$. Since $\alpha$ and $\beta$ have
different lengths, they are not conjugate under the Weyl group of
$G$. So by Theorem \ref{thm 20.2} we have two distinct
non-zero spherical nilpotent orbits in $\gg$, cf.\! Theorem
\ref{thm 8.50}.
\end{itemize}\end{exmp}\end{comment}

\subsection{Spherical Orbits and ad-Nilpotent Ideals}
\label{sub:ideals}

In this section we generalize some results from \cite{PaRo1} and \cite{PaRo2} 
to a field of good characteristic.

When $G$ is simple and classical, Panyushev gave 
simple algebraic criteria for
a nilpotent element $e \in \N$ to be spherical in 
\cite[\S 4]{Pa3}. We show that 
these criteria are still valid for a field of good characteristic.

\begin{lem}
\label{lem 20.3a} 
Let $G$ be a simple classical
algebraic group and $e \in \N$.
\begin{itemize}
\item[(i)] 
Let $e$ be a nilpotent matrix in $\mathfrak{sl}_n$ or $\mathfrak{sp}_n$. 
Then $e$ is spherical if
and only if  $e^2=0$.
\item[(ii)] 
Let $e$ be a nilpotent matrix in $\mathfrak{so}_n$. 
Then $e$ is spherical if and only if the rank of $e^2$ is at most one.
\end{itemize}
\end{lem}

\begin{proof}
Let $e$ be a nilpotent matrix in $\mathfrak{sl}_n$ or $\mathfrak{sp}_n$. 
If $e$ is spherical, then 
$\pi_e=[1^j,2^i]$, for appropriate $i$ and $j$, see Remark \ref{rem 8.50}. 
By considering the corresponding Jordan blocks
for $\pi_e$, we see that $e^2=0$. Conversely, if
$e^2=0$, then $e$ is conjugate to an element $e'$ with partition
$\pi_{e'}=[1^j,2^i]$ and so $e$ is spherical, again by Remark \ref{rem 8.50}.

Let $e$ be a nilpotent matrix in $\mathfrak{so}_n$. 
If $e$ is spherical, then 
$\pi_e=[1^j,2^i]$ or $\pi_e=[1^j,2^i,3]$, for appropriate $i$ and
$j$, see Remark \ref{rem 8.50}. By considering the corresponding
Jordan blocks for $\pi_e$, we see that either
$e^2 = 0$ or $e^2$ has partition $\pi_{e^2}=[1^k,2]$. Thus the rank
of $e^2$ is either $0$ or $1$. Conversely, if the rank of $e^2$
is at most $1$, then $e$ is conjugate to an element $e'$ with partition
$\pi_{e'} = [1^j,2^i]$ or $\pi_{e'} = [1^j,2^i,3]$ and so $e$ is spherical,
again by Remark \ref{rem 8.50}.
\end{proof}

In \cite{PaRo1} and \cite{PaRo2}, 
D.I.\ Panyushev and the second author gave a 
classification of the spherical ideals of $\bb=\Lie B$ contained 
in $\bb_u=\Lie R_u(B)$, where $B$ is a Borel subgroup of $G$ 
in characteristic $0$. 
An ideal $\cc$ of
$\bb$ is \emph{ad-nilpotent} if $\cc$ is
contained in $\bb_u$. An ad-nilpotent ideal $\cc$ of $\bb$ is
called \emph{spherical} if its $G$-saturation $G\cdot
\cc=\{x\cdot e \mid x \in G, e \in \cc\}$ is a spherical
$G$-variety. First in \cite[Cor.\ 2.4]{PaRo1} it is proved that if
$\aaa$ is an Abelian ideal of $\bb$, then $\aaa$ is spherical. In 
\cite[Prop.\ 4.1 and Thm.\  4.2]{PaRo2} it is proved that 
there are non-abelian spherical ideals 
only if $G$ is not simply-laced, that is 
if the Dynkin diagram of $G$ has a multiple bond.

Theorem 2.3 in \cite{PaRo1} states that any $G$-orbit meeting an
abelian ad-nilpotent ideal $\aaa$ is spherical. This is proved by
means of the fact that an orbit $G\cdot e$ is spherical if and
only if $\ad(e)^4=0$, see \cite[Cor.\ 2.2]{Pa3}. Unfortunately, 
this equivalence  is no longer true in positive characteristic, see
Example \ref{ex 20.10}. However, the forward implication of this
equivalence is still valid in good characteristic.

\begin{lem}
\label{lem 20.30} 
If $e \in \N$ is spherical, then $\ad(e)^4=0$.
\end{lem}

\begin{proof} 
If $e$ is spherical, then by Theorem \ref{thm 8.40}
$\htt(e)\leqslant 3$. Let $\gg=\bigoplus_{i=-3}^{3}\gg(i)$
be the grading of $\gg$ afforded by an associated cocharacter 
in $\Omega_G^a(e)$.  We have that $e \in \gg(2)$. 
Consequently,  $\ad(e)^4(\gg(i)) \subseteq \gg(i+8) = \{ 0 \}$ for
any $-3 \leqslant i \leqslant 3$.
Consequently, $\ad(e)^4=0$ on all of $\gg$.
\end{proof}

The next example shows that the converse of Lemma \ref{lem 20.30}
is not true in general in positive characteristic.

\begin{exmp}
\label{ex 20.10}
Let $G=\SL_3(k)$ and $\cha k=3$. So
$\gg=\mathfrak{sl}_3(k)$. Set $e=e_{2,1}+e_{3,2}$, where $e_{i,j}$
is the elementary matrix with a 1 in the $(i,j)$ position and 0's
elsewhere. So $e$ is a regular nilpotent element in $\gg$.
Consider the grading of $\gg$ 
afforded by an associated cocharacter in $\Omega_G^a(e)$. 
We have $\gg=\bigoplus_{i=-2}^{2}\gg(2i)$. In order to prove
$\ad(e)^4=0$, it is sufficient to show that
$\ad(e)^4(\gg(-4))=\{0\}$. Clearly, $\gg(-4)=ke_{1,3}$. Now
$\ad(e)(e_{1,3})=e_{2,3}-e_{1,2}$ and
$\ad(e)(e_{2,3}-e_{1,2})=e_{1,1}-2e_{2,2}+e_{3,3}$. Since $\cha
k=3$, we have $e_{1,1}-2e_{2,2}+e_{3,3}=e_{1,1}+e_{2,2}+e_{3,3}$
and $e_{1,1}+e_{2,2}+e_{3,3} \in Z(\gg)$. Thus, $\ad(e)^4=0$.
However, $e$ is not spherical, as %it has corresponding partition
$\pi_e=[3]$, see Remark \ref{rem 8.50}.
\end{exmp}

We note that Proposition 4.1 and Theorem 4.2 in \cite{PaRo2} both 
also hold in good characteristic, as their proofs only require
properties of the underlying root system $\Psi$ and the results established in 
Lemmas \ref{lem 20.3a} and \ref{lem 20.30}.

So we are left to show that if $\aaa$ is an abelian ad-nilpotent
ideal, then $\aaa$ is spherical. Since $G\cdot \aaa$ is irreducible, it is
the closure of some nilpotent orbit, 
say $\overline{G\cdot e} = G\cdot \aaa$. 
The maximal abelian ad-nilpotent ideals of $\bb$
are the same in good characteristic as in 
characteristic zero, see Table 1 in \cite{Ro} and
Tables I and II in \cite[\S 4]{PaRo1}.
Using the description of the orbits in Tables I and II in \cite[\S 4]{PaRo1},
we infer that the Bala--Carter
label of $G\cdot e$ is of the form $A_1^t$, so $G\cdot e$
is spherical, thanks to Theorem \ref{thm 20.2}.
Since $G\cdot e$ is open in $G\cdot \aaa$, it follows that  
$G\cdot \aaa$ is spherical.
It is straightforward to get the 
sphericity of $G\cdot\aaa$ for any abelian ideal $\aaa$ of $\bb$
from the sphericity  result of the maximal abelian ideals.
Thus we have established the following.

\begin{thm}
\label{thm 20.04} 
Let $\aaa$ be an abelian ad-nilpotent ideal of
$\bb$. Then $\aaa$ is spherical. 
\end{thm}

As a corollary of Theorem \ref{thm 20.04} we get \cite[Thm.\ 1.1]{Ro}
in good characteristic. 

\begin{cor}
\label{cor 20.2} 
Let $P$ be a parabolic subgroup of $G$
and let $\aaa$ be an abelian ideal of $\Lie P$ in $\Lie R_u(P)$.
Then $P$ acts on $\aaa$ with finitely many orbits.
\end{cor}

\begin{rem}
We note that Theorem \ref{thm 20.04} and Corollary \ref{cor 20.2} 
do in fact hold in arbitrary characteristic, cf.\ \cite[Thm.\ 1.1]{Ro}.
%However, in contrast to the arguments in \cite{Ro}, 
%our proofs of these facts in good characteristic presented here are free 
%of case by case considerations 
%(of course, modulo the classification of the spherical nilpotent orbits).
\end{rem}

\begin{rem}
If $\cc$ is a spherical ideal of $\bb$, then clearly $B$ acts on 
$\cc$ with a finite number of orbits.
However, the converse does not hold. There are many additional instances
when $B$ acts on a given ideal $\cc$ of $\bb$ only with a finite number 
of orbits, e.g. see the results in 
\cite{hilleroehrle} and \cite{juergensroehrle}.
\end{rem}

\subsection{A Geometric Characterization of Spherical Orbits}
\label{sub:geom}
In this subsection we describe a formula characterizing spherical $G$-orbits
in a simple algebraic group $G$ 
in terms of elements of the Weyl group $W$ of $G$
that is proved in \cite[Thm.\ 1]{CaCaCo}. 
For $x \in G$ the conjugacy class $G\cdot x$ is spherical if  
$G\cdot x$ is a spherical variety. 
While this characterization in {\it loc.\ cit.} is
based on case by case arguments, 
recently, G.\ Carnovale \cite[Thm.\ 2]{Carno} gave a proof of this result
which is free of 
case by case considerations and applies in good odd characteristic.
Using the arguments from \cite{CaCaCo} combined with our classification
of the spherical unipotent nilpotent orbits, Remark \ref{rem 8.55}, 
we can generalize this formula to good characteristic. 

Let $G$ be simple and suppose that $p$ is good for $G$.
Fix a Borel subgroup $B$ of $G$. Let $W$ be the Weyl group of $G$ and
let $BwB$ be the $(B,B)$-double coset of $G$ containing $w \in W$. 
The following was shown in \cite{CaCaCo} in an
argument independent of the  characteristic of the underlying field:
Suppose that $\CO$ is a conjugacy class in $G$ which intersects 
the double coset $BwB$ so that 
$\dim \CO = \ell(w) + \rk(1 - w)$ holds. 
Then  $\CO$ is spherical.
Here $\rk(1 - w)$ denotes the rank of the linear map $1 - w$ in the
standard representation of $W$ and $\ell$ is the usual length function 
of $W$ with respect to a distinguished set of generators of $W$.
Conversely, let $\CO$ be a spherical conjugacy class in $G$ and let
$BwB$ be the $(B,B)$-double coset containing the dense $B$-orbit 
in $\CO$. Then $\dim \CO = \ell(w) + \rk(1 - w)$, see \cite[Thm.\ 2]{Carno}.
Consequently, this gives a geometric characterization of the spherical 
conjugacy classes in $G$.
For proofs we refer the reader to \cite{CaCaCo} and \cite{Carno}. 
Observe that as a consequence of the finiteness of the 
Bruhat decomposition of $G$ and the fact that  
any $(B,B)$-double coset and any conjugacy class of $G$ are
irreducible subvarieties of $G$, 
for a given conjugacy class
$\CO$ in $G$ there is a unique $w \in W$
such that $\CO \cap BwB$ is dense in $\CO$.

\begin{thm} (\cite[Thm.\ 1]{CaCaCo})
\label{thm:CaCaCo}
Let $\CO$ be a conjugacy class in $G$ 
and let $w \in W$ be such that $\CO \cap BwB$ is dense in $\CO$.
Then $\CO$ is spherical if and only if 
$\dim \CO = \ell(w) + \rk(1 - w)$.
\end{thm}

\subsection{Bad Primes and Spherical Nilpotent Orbits}
\label{sub:bad}
Finally, we 
briefly discuss the situation when the characteristic of $k$ is bad
for $G$. In this case the classification of the nilpotent orbits in $\N$ is
different from that in good characteristic, see \cite[\S 5.11]{Ca}. 
However, there is still only a finite number of nilpotent orbits,
\cite{HoltSpaltenstein}. 
Unfortunately, our methods do not allow us to give a classification 
of the spherical nilpotent orbits in this case.
For, in our classification we made use of 
the height of a nilpotent orbit,
where the height is defined via an associated
cocharacter. However, it is not known whether 
associated cocharacters always exist for all nilpotent elements 
in bad characteristic, cf.\ \cite[\S 5.14, \S 5.15]{Ja1}.

In principle one can still determine whether a given nilpotent orbit is
spherical by a case by case analysis.
Next we give two examples of this.
In particular, we show that Theorem \ref{thm 20.2} fails in 
bad characteristic in general.
These examples show that there can be additional spherical nilpotent orbits 
in bad characteristic.

\begin{exmps}
\label{eg 20.20}
(i).
Let $G$ be of type $B_2$ and $\cha k=2$. 
Let $\alpha$ and $\beta$ be the simple roots of $\Psi$ with
$\alpha$ the long root. Let
$e=e_{\alpha+\beta}+e_{\alpha+2\beta}$. According to
\cite[\S 5.14]{Ja1} the centralizer $C_G(e)$ is the unipotent
radical of a Borel subgroup of $G$. Thus, by Lemma \ref{lem 1.30}, 
$C_G(e)$ is a spherical subgroup of $G$ and so $e$ is
spherical. Note that the $G$-orbit of $e$ does not contain an
element of the form $e_{\alpha}$ or $e_{\beta}$, but $e$ is
still spherical. Thus, Theorem \ref{thm 20.2} is no longer
true in bad characteristic.
Moreover, $e$ is distinguished in $\gg$, \cite[\S 5.14]{Ja1}.
This shows that Theorem \ref{thm 20.02} can also fail for bad characteristic.

(ii).  
Let $G$ be of type $G_2$ and $\cha k=3$. 
Let $\alpha$ and $\beta$ be the simple roots of $\Psi$ with
$\alpha$ the long root. 
Let $e=e_{\alpha+2\beta}+e_{2\alpha+3\beta}$. According to
\cite[\S 5.15]{Ja1}, the centralizer $C_G(e)$ is the unipotent
radical of a Borel subgroup of $G$. Thus, by Lemma \ref{lem 1.30}, 
$C_G(e)$ is a spherical subgroup of $G$ and so $e$ is
spherical. Again, the $G$-orbit of $e$ does not
contain an element of the form $e_{\alpha}$ or $e_{\beta}$, but
$e$ is spherical. 
Again, $e$ is distinguished in $\gg$, \cite[\S 5.15]{Ja1}.
\end{exmps}

%%%%%%%%%%%%%%%%%%%%%%%%%%%%%%%%%%%%%%%%%%%%%%%%%%%%%%%%%%%%%%%%%%%%%%
%%%%%%%%%%%%% Acknowledgements
%%%%%%%%%%%%%%%%%%%%%%%%%%%%%%%%%%%%%%%%%%%%%%%%%%%%%%%%%%%%%%%%%%%%%%

\bigskip {\bf Acknowledgements}:
The first author acknowledges funding by the EPSRC.
We are grateful to S.M.~Goodwin for providing
the relative version $\mathsf{DOOBSLevi}$ of his program
$\mathsf{DOOBS}$ that was used in Subsection \ref{sub:ex}
to determine the sphericity of the nilpotent orbits of height $3$ 
for the exceptional cases and for very helpful discussions and
improvements of the paper. 
We would also like to thank the referee for suggesting some
improvements.

\bigskip

%%%%%%%%%%%%%%%%%%%%%%%%%%%%%%%%%%%%%%%%%%%%%%%%%%%%%%%%%%%%%%%%%%%%%%
%%%%%%%%%%%%% bibliography
%%%%%%%%%%%%%%%%%%%%%%%%%%%%%%%%%%%%%%%%%%%%%%%%%%%%%%%%%%%%%%%%%%%%%%

\end{document}